\documentclass{amsart}
\usepackage{amsmath,amsthm}
\usepackage {latexsym}
\usepackage{amssymb}

\newcommand{\less}{\lesssim}

\newcommand{\beal}{\begin{align}}
\newcommand{\enal}{\end{align}}
\newcommand{\bealn}{\begin{align*}}
\newcommand{\enaln}{\end{align*}}
\newcommand{\bear}{\begin{eqnarray}}
\newcommand{\eear}{\end{eqnarray}}
\newcommand{\beeq}{\begin{equation}}
\newcommand{\eneq}{\end{equation}}

\newcommand{\const}{\mbox{\rm const}}

\newcommand{\spec}{{\rm spec}}

\newcommand{\eps}{{\varepsilon}}
\newcommand{\R}{{\mathbb R}}
\newcommand{\Compl}{{\mathbb C}}

\newcommand{\Z}{{\mathbb Z}}
\newcommand{\Nat}{{\mathbb N}}
\newcommand{\calC}{{\mathcal C}}

\newcommand{\calF}{{\mathcal F}}
\newcommand{\calL}{{\mathcal L}}
\newcommand{\calQ}{{\mathcal Q}}
\newcommand{\calS}{{\mathcal S}}

\newcommand{\calK}{{\mathcal K}}

\newcommand{\tileps}{{\tilde{\epsilon}}}

\newcommand{\la}{\langle}
\newcommand{\ra}{\rangle}

\renewcommand{\ln}{\log}

\newcommand{\IS}{IS}

\def\nn{\nonumber}
\def\bm{\left[ \begin{array}{cc}}
\def\endm{\end{array}\right]}

\newtheorem{theorem}{Theorem}
\newtheorem{lemma}[theorem]{Lemma}
\newtheorem{defi}[theorem]{Definition}
\newtheorem{cor}[theorem]{Corollary}

\newtheorem{proposition}[theorem]{Proposition}

\theoremstyle{remark}
\newtheorem{remark}[theorem]{Remark}

\def\calR{{\mathcal R}}

\def\il{\int\limits}

\def\half{\frac12}

\def\cQ{{\mathcal Q}}
\def\calE{{\mathcal E}}
\def\dom{{\rm Dom}}

\def\cL{\mathcal L}
\def\tilphi{\phi}
\def\tilm{m}

\def\tiltheta{\theta}

\renewcommand{\Im}{\,{\rm Im}\,}
\renewcommand{\Re}{\,{\rm Re}\,}

\def\Bla{\Big\langle}
\def\Bra{\Big\rangle}

\renewcommand{\hat}{\widehat}
\renewcommand{\epsilon}{\eps}
\numberwithin{equation}{section}
\numberwithin{theorem}{section}

\begin{document}

\title{Renormalization and blow up for charge one equivariant critical wave maps.}

\author{J.\ Krieger}
\address{Harvard University, Dept. of Mathematics, Science Center, 1 Oxford Street, Cambridge, MA 02138, U.S.A.}
\email{jkrieger@math.harvard.edu}

\author{W.\ Schlag}
\address{Department of Mathematics, The University of Chicago, 5734 South University Avenue, Chicago, IL 60637, U.S.A.}
  \email{schlag@math.uchicago.edu}

\author{D.\ Tataru}
\address{Department of Mathematics, The University of California at Berkeley, Evans Hall, Berkeley, CA 94720, U.S.A.}
\email{tataru@math.berkeley.edu}

\thanks{The authors were partially supported by the National Science Foundation, J.\ K.\ by DMS-0401177, W.\ S.\ by
DMS-0617854, D.~T.\ by DMS-0354539, and DMS-0301122. The first
author thanks UC Berkeley and the University of Chicago for their
hospitality. The second author thanks Harvard University for its
hospitality, and Fritz Gesztesy for helpful discussions.}

\subjclass{Primary 35L05, 35Q75, 35P25} \keywords{Nonlinear
hyperbolic equations, wave maps, spectral theory of strongly
singular Sturm Liouville operators, zero energy resonances}

\begin{abstract} We prove the existence of equivariant finite time
blow-up solutions for the wave map problem from $\R^{2+1}\to S^2$ of
the form $u(t,r)= Q(\lambda(t) r) + \calR(t,r)$ where $u$ is the
polar angle on the sphere, $Q(r)=2\arctan r$ is the ground state
harmonic map, $\lambda(t)=t^{-1-\nu}$, and $\calR(t,r)$ is a
radiative error with local energy going to zero as $t\to0$. The
number $\nu>\frac12$ can be described arbitrarily. This is
accomplished by first "renormalizing" the blow-up profile, followed
by a perturbative analysis.
\end{abstract}

\maketitle

\section{Introduction}
\label{sec:intro}

We consider Wave Maps $U:\R^{2+1}\rightarrow S^{2}$ which are
equivariant with co-rotation index~$1$. Specifically, they satisfy
$U(t,\omega x)=\omega U(t,x)$ for $\omega\in SO(2,\R)$, where the
latter group acts in standard fashion on $\R^{2}$, and the action on
$S^{2}$ is induced from that on $\R^{2}$ via stereographic projection.
Wave maps are characterized by being critical with respect to the
functional
\begin{equation}\nonumber
  U\rightarrow \int_{\R^{2+1}}\la \partial_{\alpha}U, \partial^{\alpha}U\ra \,d\sigma,\,\alpha=0,1,2
\end{equation}
with Einstein's summation convention being in force,
$\partial^{\alpha}=m^{\alpha\beta}\partial_{\beta}$,
$m_{\alpha\beta}=(m^{\alpha\beta})^{-1}$ the Minkowski metric on
$\R^{2+1}$, and $\,d\sigma$ the associated volume element. Also,
$\la\cdot,\cdot\ra$ refers to the standard inner product on $\R^{3}$
if we use ambient coordinates to describe $u$, $\partial_{\alpha}u$
etc. Recall that the energy is preserved:
\[
\calE(u) = \int_{\R^2} \la DU(\cdot,t),DU(\cdot,t)\ra\, dx =\const
\]
If one instead uses spherical coordinates, and lets $u$ stand for the
longitudinal angle, and similarly use polar coordinates $r,\theta$ on
$\R^{2}$, we describe the Wave Map by $(t,r,\theta)\longrightarrow
(u(t,r), \theta)$, where now $u(t,r)$, a scalar function, satisfies
the equation
\begin{equation}\label{TheEquation}
  -u_{tt}+u_{rr}+\frac{u_{r}}{r}=\frac{\sin(2u)}{2r^{2}}
\end{equation}
The problem at hand is {\it{energy critical}}, meaning that the
conserved energy is invariant under the natural re-scaling
$U\rightarrow U(\lambda t, \lambda x)$ (using the original coordinates
and meaning of $U$). By contrast, the analogous wave map problem on
$\R^{n+1}$, $n\geq 3$ is energy-supercritical in the sense that the
natural scale-invariant Sobolev space is then $\dot{H}^{\frac{n}{2}}$,
and the corresponding norm $\|u\|_{\dot{H}^{\frac{n}{2}}}$ is not
expected to be controlled globally-in-time for general initial data,
which leads to the general belief that in this case, there should not
be a good well-posedness theoory for general initial data,
{\it{irrespective of the target.}} Indeed, singular wave maps stemming
from $C^{\infty}$-data have been constructed on background $\R^{3+1}$
with target $S^{3}$ in \cite{Sha1}, and with origin $\R^{n+1}$, $n\geq
4$ and for more general targets in \cite{Sha2}.
\\

In the critical case, global well-posedness is expected for hyperbolic
targets, while singularity development is expected for certain
positively curved targets, such as $S^{2}$. More precisely, numerical
evidence in \cite{Bi}, \cite{Li} strongly suggests singularity
development for equivariant wave maps of co-rotation index $1$ from
$\R^{2+1}$ to $S^{2}$ with smooth data, while wave maps from
$\R^{2+1}$ to ${\mathbf{H}}^{2}$, and more generally
${\mathbf{H}}^{k}$, $k\geq 2$, are expected to preserve the
regularity\footnote{By this we mean that if initial data have
  regularity $H^{1+\delta}$, $\delta>0$, the Wave Map can be uniquely
  globally extended in this class.}  of the initial data. Further
evidence for possible singularity development up to this date
specifically, in the co-rotation~$1$ equivariant case has been
recently provided in \cite{Co}. We note that a fairly satisfactory
understanding has been achieved for small-energy wave maps from
$\R^{2+1}$ to general targets \cite{Tao}, \cite{Tat}, \cite{Kri}, as
well as for rotationally invariant wave maps and general initial
data \cite{Chr-Tah}, \cite{Str1}. In particular, {\it{it is known
that the
    latter never develop singularities}} \cite{Str1}, and that for
equivariant wave maps of co-rotation index~$1$, regularity breakdown
can only occur in an {\it{energy concentration scenario}}~\cite{Str2}.
For equivariant wave maps, it is known that regularity of the initial
data is preserved (see previous footnote) provided the target
satisfies a geodesic convexity condition~\cite{Sh-Tah}.
\\
Our objective in this paper is to rigorously demonstrate regularity
breakdown for equivariant wave maps $u: \R^{2+1}\longrightarrow S^{2}$
of co-rotation index $1$ with certain $H^{1+}$ regular initial data.
More precisely, the data $(u,\,u_{t})$ will be of class
$H^{1+\delta}\times H^{\delta}$ for some $\delta>0$. It is well-known
that such data result in unique local solutions of the same
regularity until possible breakdown occurs via an energy-concentration
scenario. We note that a result of Struwe shows that if the solution
is indeed $C^{\infty}$--smooth before breakdown\footnote{This result
  most likely can be adapted to solutions of lesser smoothness}, such
a scenario can only happen by the bubbling off of a harmonic
map~\cite{Str2}: specifically, let $Q(r): \R^{2}\longrightarrow S^{2}$
be an equivariant hamonic map, which can be constructed for every
co-rotation index $k\in\Z$ (for example, for $k=1$ stereographic
projection will do). We shall identify $Q(r)$ with the longitudinal
angle, as above. Then according to \cite{Str2}, if an equivariant wave
map $u$ of co-rotation index $k=1$, again identified with the
longitudinal angle, with smooth initial data at some time $t_{0}>0$
breaks down at time $T=0$, then energy focuses at the origin, and
there is a decomposition
\begin{equation}\nonumber
  u(t,r)=Q(\lambda(t) r) +\epsilon(t,r),\,\text{$Q(r)$ a co-rotation $k=1$ index equivariant harmonic map}
\end{equation}
where there is a sequence of times $t_{i}\rightarrow 0$, $t_{i}<0$,
$i=1,2,\ldots$, with $\lambda(t_{i})|t_{i}|\rightarrow \infty$, such
that the rescaled functions $u(t_{i},\frac{r}{\lambda(t_{i})})$
converge to $Q(r)$ in the strong energy topology.
\\

This is borne out by our main theorem. We let $Q(r)$ represent the
standard harmonic map of co-rotation $k=1$, i.e., $Q(r)=2\arctan r$.
Recall that in the equivariant formulation the energy is
\[
\calE(u) = \int_{\R^2} \Big[\frac12 ( u_t^2 + u_r^2) +
\frac{\sin^2(u)}{2r^2}\Big] \,r\, dr
\]
The {\em local} energy relative to the origin is defined as
\[
\calE_{\rm loc}(u) = \int_{r<t} \Big[\frac12 ( u_t^2 + u_r^2) +
\frac{\sin^2(u)}{2r^2}\Big] r\, dr
\]
It is well-known that for equivariant wave-maps singularities can only
develop at the origin and this happens at time zero iff
\[
\liminf_{t\to0} \calE_{\rm loc}(u)(t)>0
\]
The following theorem is the main result of this paper. Note that we
need to "renormalize" the profile $Q(r\lambda(t))$ by means of a large
perturbation (denoted $u^e$ below). We find it convenient to solve backwards in time,
with blow-up as $t\to0+$.

\begin{theorem}\label{Main}
  Let $\nu >\frac12$ be arbitrary and  $t_{0} > 0$ be sufficiently small.
  Define $\lambda(t)=t^{-1-\nu}$ and fix a large integer $N$. Then there exists a
  function\footnote{We refer to this as an "elliptic profile
    modifier"; see Section~\ref{sec:mod} for a detailed explanation of this notion. Also, $C^\beta$ for noninteger $\beta$
means $C^{[\beta],\beta-[\beta]}$}
  $u^e$ satisfying
  \[
  u^e\in C^{\nu+1/2-}(\{t_0>t>0, \;|x|\le t\}),\qquad \calE_{\rm
    loc}(u^e)(t)\lesssim (t\lambda(t))^{-2}\,|\log t|^2 \text{\  \ as\  \ }t\to0
  \]
  and a blow-up solution $u$ to \eqref{TheEquation} in $[0,t_0]$ which has the form
  \[
  u(r,t) = Q(\lambda(t) r) + u^e(r,t) + \epsilon(r,t), \qquad 0\le r\le t
  \]
  where $\epsilon$ decays at $t=0$; more precisely,
\[
\epsilon \in t^N H_{\rm loc}^{1+\nu-}(\R^2), \qquad  \epsilon_t \in
t^{N-1} H_{\rm loc}^{\nu-}(\R^2), \qquad \calE_{\rm
    loc}(\eps)(t)\lesssim t^N \text{\  \ as\  \ }t\to0
\]
with spatial norms that are uniformly controlled as $t\to0$. Also,
$u(0,t)=0$ for all $0<t<t_0$. The solution $u(r,t)$ extends as an
$H^{1+\nu-}$ solution to all of $\R^2$ and the energy of $u$
concentrates in the cuspidal region $0\leq r\lesssim
\frac{1}{\lambda(t)}$ leading to blow-up at $r=t=0$.
\end{theorem}

We remark that a somewhat surprising feature of our theorem is that
the blow-up rate is prescribed. This is in stark contrast to the usual
modulation theoretic approach where the rate function is used to
achieve orthogonality to unstable modes of the linearized problem.
Heuristically speaking, there are two types of instabilities which
typically arise in linearized problems: those due to symmetries of the
nonlinear equation (typically leading to algebraic growth of the
linear evolution) and those that produce exponential growth in the
linear flow (due to some kind of discrete spectrum). For example, the
latter arises in the recent work on ``center-stable manifolds'',
see~\cite{Sch}, \cite{KS1}, \cite{KS2} whereas for the former
see~\cite{KrSch}. Both types can lead to blow up. Here we do not have
any discrete spectrum in the linearized equation, but rather a
zero-energy resonance which is due to the scaling symmetry. It is unclear
at this point which role (other than a technical one) the resonance
plays in the formation of the blow-up. Indeed, our approach is really
non-perturbative as the crucial elliptic profile modifier produces a
{\em large} perturbation of the basic profile~$Q$. The perturbative
component of our proof deals with the removal of errors produced
by the elliptic profile modifier (it is crucial that these errors decay
rapidly in time).

A recent preprint\footnote{The conclusions of our paper were reached
  before the appearance of this preprint} by I.~Rodnianski and
J.~Sterbenz \cite{Ro-St} details the construction of generic sets of
initial data (including smooth data) resulting in blow-up with a rate
$\lambda(t)\sim \frac{\sqrt{|\log t|}}{t}$ for equivariant wave maps
from $\R^{2+1}$ to $S^{2}$ with co-rotation index $k\geq 4$.  These
data can be chosen arbitrarily close to the corresponding co-rotation
$k$ harmonic map with respect to a suitable norm stronger than
$\|.\|_{H^{1}}$. This latter behavior appears specific for
sufficiently large co-rotation indices {\it{but, due to numerical
    experiments (e.g., \cite{Bi}) is not expected for the case of
    co-rotation index one.}}  More precisely, numerical experiments suggest
that perturbations of the ``ground state'' $Q(r)=2\arctan(r)$ need to be
of a certain size (depending on their profile) to result in blow-up.
Indeed, our theorem, which is partly based on perturbative
techniques, has a non-perturbative flavor in that the {\it{elliptic
    profile modifier cannot be made small at fixed time $t=t_{0}>0$
    and with a  fixed profile $\lambda(t)$.}} On a technical level, we remark that
the corresponding linearlized operator has zero energy as an eigenvalue
for $k>1$ but for $k=1$ zero energy becomes a {\em resonance} (indeed,
$\partial_\lambda Q_k(\lambda r)|_{\lambda=1}\in L^2(0,\infty)$
iff $k>1$ where $Q_k(r)=2\arctan(r^k)$).

Our argument is correspondingly divided into two parts: first, we use
a direct method, exploiting the algebraic fine structure of the
system, to find an approximate solution $Q(\lambda(t)r)+u^{e}(t,r)$.
Roughly speaking, one may think of $u_{e}(.,.)$ as being obtained by a
finite sequence of approximations which alternately improve the
accuracy near the light cone and near the origin. To model the
solution near the light cone, one introduces the coordinates $(a,t)$
where $a=\frac{r}{t}$ and reduces to solving an elliptic problem in
$a$ by neglecting time derivatives. More precisely, one treats time
derivatives as error source terms, which get decimated by iterating
the elliptic construction. Similarly, one improves accuracy near the
origin $r=0$ by working with the coordinates $(R, t)$ where
$R=\lambda(t)r$, again reducing to an elliptic problem by neglecting
time derivatives. This process does not lead to an actual solution, as
one ``keeps losing time derivatives'', which leads to worse and worse
implicit constants. Thus in a second stage, we construct a
parametrix for the wave equation which is obtained by passing
to coordinates $(R, \tau)$ where $R=\lambda(t)r$,
$\tau=\frac{1}{\nu}t^{-\nu}$. This in turn relies on a careful
analysis of the spectral and scattering theory of the Schr\"odinger operator
which arises by linearizing around $Q(r)$. The remaining error is then iterated away by
continued application of the wave parametrix.

A surprising feature of Theorem~\ref{Main} is the fact that
{\it{blow-up may be arbitrarily slow, as we can prescribe $\nu$
    arbitrarily large.}} However, the data leading to this blow-up are
not generic, and indeed rather difficult to describe. We observe
that in~\cite{Bi} solutions blowing up with $\lambda(t)\sim
t^{-2.3}$ were observed numerically, corresponding to a transient
regime dividing blow-up from global smoothness and scattering. In
particular, this blow-up rate appears to correspond to a set of
initial data of co-dimension one. Our initial data sets seem to lie
 on manifolds of very large co-dimension, which increases with $\nu$ (we plan to
 return to a rigorous treatment of this conditional stability issue in a later paper). In
particular, it appears unlikely that such solutions corresponding to
$\nu>>1$ would be detected numerically.

Finally, we note that our technique quite likely lends itself to
constructing blow-up solutions for higher homotopy indices, too, as
well as to problems of a similar nature, such as the critical
Yang-Mills equation.

\section{Approximate solutions}
\label{sec:mod}

\subsection{The elliptic profile modifier}

In this section we show how to construct an arbitrarily good
approximate solution to the wave map equation as a perturbation of a
time-dependent harmonic map profile
\[
u_0 = Q(R), \qquad R =r \lambda(t)
\]
with the polynomial timescale
\[
\lambda(t) = t^{-1-\nu}
\]
To describe the approximate solution we use the time variable, the
variable $R$ which corresponds to the harmonic map scale, and the
self-similar variable $a=r/t$ which is useful in analyzing the
behavior near the cone. The only trade off in this construction is
that we need to allow singularities of the form
\[
(1-a^2)^\nu (\ln(1-a^2))^k
\]
as we approach the cone. Thus the larger the parameter $\nu$, the
better the regularity of the approximate solutions.

\begin{theorem}
  \label{thm:sec2} Let $k \in \Nat$. There exists an approximate
  solution $u_{2k}$ for \eqref{TheEquation} of the form
  \[
  u_{2k-1}(r,t) = Q(\lambda(t) r) + \frac{c_k}{(t\lambda)^2} R {\ln(1+ R^2)} +
  O\left( \frac{R^{-1} (\ln (1+R^2))^2}{(t\lambda)^2}\right)
  \]
  so that the corresponding error has size
  \[
  e_{2k-1} = O\left(\frac{R(\log (2+R))^{2k-1}}{t^2 (t \lambda)^{2k}}\right)
  \]
 Here the $O(\cdot)$ terms are uniform in $0\le r\le t$ and
  $0<t<t_0$ where $t_0$ is a fixed small constant.
\end{theorem}

\begin{remark}
  In the proof we obtain $u_{2k-1}$ and $e_{2k-1}$ which are analytic
  inside the cone and $C^{\frac12+\nu-}$, respectively
  $C^{-\frac12+\nu-}$ on the cone, with a good asymptotic expansion
  both on the $R$ scale and near the cone.

More precisely, using our notations
  defined below we have
  \[
  u_{2k-1} \in Q(\lambda(t) r)+ \frac{1}{(t\lambda)^2} \IS^3(R {\ln
    R},\cQ_k)
  \]
  while the error satisfies
  \[
  t^2 e_{2k-1} \in \frac{1}{(t \lambda)^{2k}} \IS^1(R (\ln
  R)^{2k-1},\cQ'_{k-1})
  \]
\end{remark}

\begin{proof}
  We iteratively construct a sequence $u_k$ of better approximate
  solutions by adding corrections $v_k$,
  \[
  u_k = v_{k} + u_{k-1}
  \]
  The error at step $k$ is
  \[
  e_k = (-\partial_t^2 + \partial_r^2 +\frac{1}r \partial_r) u_k
  -\frac{\sin(2u_k)}{2r^{2}}
  \]
  To construct the increments $v_k$ we first make a heuristic
  analysis. If $u$ were an exact solution, then the difference
  \[
  \epsilon = u-u_{k-1}
  \]
  would solve the equation
  \begin{equation}\label{eq:gleich}
  (-\partial_t^2 + \partial_r^2 +\frac{1}r \partial_r) \epsilon -
  \frac{\cos(2u_{k-1})}{2r^{2}} \sin (2\epsilon) +
  \frac{\sin(2u_{k-1})}{2r^{2}} (1-\cos(2\epsilon)) = e_{k-1}
  \end{equation}
  In a first approximation we linearize this equation around $\epsilon
  =0$ and substitute $u_{k-1}$ by $u_0$. Then we obtain the linear
  approximate equation
  \begin{equation}
    \left(-\partial_t^2 + \partial_r^2 +\frac{1}r \partial_r -
      \frac{\cos(2u_0)}{r^{2}}\right) \epsilon  \approx e_{k-1}
    \label{eps}  \end{equation} For $r \ll t$ we expect the time
  derivative to play a lesser role  so we neglect it and we are left
  with an elliptic equation with respect to the variable $r$,
  \begin{equation}
    \left(\partial_r^2 +\frac{1}r \partial_r -
      \frac{\cos(2u_0)}{r^{2}}\right) \epsilon  \approx e_{k-1}, \qquad r
    \ll t \label{epsodd}\end{equation} For $r \approx t$ we can
  approximate $ \cos(2u_0)$ by $1$ and rewrite \eqref{eps} in the form
  \[
  \left(-\partial_t^2 + \partial_r^2 +\frac{1}r \partial_r -
    \frac{1}{r^{2}}\right) \epsilon \approx e_{k-1}
  \]
  Here the time and spatial derivatives have the same strength.
  However, we can identify another principal variable, namely $a =
  r/t$ and think of $\epsilon$ as a function of $(t,a)$.  As it turns
  out, neglecting a "higher order" part of $e_{k-1}$ which can be directly included
  in $e_k$, we are able to use scaling and the exact structure of the
  principal part of $e_{k-1}$ to reduce the above equation to a
  Sturm-Liouville problem in~$a$
  which becomes singular at~$a=1$.

  The above heuristics lead us to a two step iterative construction of
  the $v_k$'s. The two steps successively improve the error in the two
  regions $r \ll t$, respectively $r \approx t$. To be precise, we
  define $v_k$ by
  \begin{equation}
    \left(\partial_r^2 +\frac{1}r \partial_r - \frac{\cos(2u_0)}{r^{2}}\right) v_{2k+1}  = e_{2k}^0
    \label{vkodd}\end{equation} respectively
  \begin{equation}
    \left(-\partial_t^2+\partial_r^2 +\frac{1}r \partial_r - \frac1{r^{2}}\right) v_{2k}  = e_{2k-1}^0
    \label{vkeven}\end{equation} both equations having zero Cauchy
  data\footnote{The coefficients are singular at $r=0$, therefore this
    has to be given a suitable interpretation} at $r=0$.  Here at each
  stage the error term $e_k$ is split into a principal part and a
  higher order term (to be made precise below),
  \[
  e_k = e_k^0 + e_k^1
  \]
  The successive errors are then computed as
  \[
  e_{2k} = e_{2k-1}^1 + N_{2k} (v_{2k}), \qquad e_{2k+1} = e_{2k}^1 -
  \partial_t^2 v_{2k+1} + N_{2k+1} (v_{2k+1})
  \]
  where
  \begin{equation}
    -N_{2k+1}(v) =  \frac{\cos(2u_0)
      -\cos(2 u_{2k})}{r^2} v +  \frac{\sin(2u_{2k})}{2r^{2}}
    (1-\cos(2v))+ \frac{  \cos(2u_{2k})}{2r^{2}}(2
    v - \sin( 2 v))
    \label{eodd}\end{equation} respectively
  \begin{equation}
    -N_{2k}(v) =   \frac{1
      -\cos(2 u_{2k-1})}{r^2} v  +  \frac{\sin(2u_{2k-1})}{2r^{2}}
    (1-\cos(2v)) + \frac{  \cos(2u_{2k-1})}{2r^{2}}(2
    v - \sin( 2 v))
    \label{eeven}\end{equation}

  To formalize this scheme we need to introduce suitable function
  spaces in the cone \[\calC_0=\{(r,t)\::\:0\le r<t, 0<t<t_0\}\] for
  the successive corrections and errors. We first consider the $a$
  dependence. For the corrections $v_k$ we use

\begin{defi}\label{def:Q}
  For $i \in \Nat$ we let $j(i)=i$ if $\nu$ is irrational,
  respectively $j(i) = 2i^2$ if $\nu$ is rational.

  a) $\mathcal Q$ is the algebra of continuous functions $q:[0,1] \to
  \R$ with the following properties:

  (i) $q$ is analytic in $[0,1)$ with an even expansion at $0$

  (ii) Near $a=1$ we have an absolutely convergent expansion of the
  form
  \[
  q = q_0(a) + \sum_{i =1}^\infty \left(  (1-a)^{(2i-1) \nu+\frac12}\sum_{ j=0}^{ j(2i-1)}
    q_{2i-1,j}(a) (\ln (1-a))^j + (1-a)^{2i \nu+1}
    \sum_{ j=0}^{ j(2i)} q_{2i,j}(a)  (\ln
    (1-a))^j\right)
  \]
  with analytic coefficients $q_0$, $q_{ij}$

  b) $\mathcal Q_m$ is the algebra which is defined similarly, with
  the additional requirement that
  \[
  q_{ij}(1) =0 \text{\ \ if\ }\ \ i \geq 2m+1, \text{odd}
  \]
\end{defi}

For the errors $e_k$ we introduce

\begin{defi}
  a)  With $j(i)$ as above, $\mathcal Q'$ is the space of continuous functions $q:[0,1] \to
  \R$ with the following properties:

  (i) $q$ is analytic in $[0,1)$ with an even expansion at $0$

  (ii) Near $a=1$ we have a convergent expansion of the form
  \[
  q = q_0(a) + \sum_{i =1}^\infty \left( (1-a)^{(2i-1) \nu-\frac12} \sum_{ j=0}^{ j(2i-1)}
    q_{2i-1,j}(a) (\ln (1-a))^j + (1-a)^{2i \nu}
    \sum_{ j=0}^{ j(2i)} q_{2i,j}(a) (\ln
    (1-a))^j\right)
  \]
  with analytic coefficients $q_0$, $q_{ij}$

  b) $\mathcal Q'_m$ is the space which is defined similarly, with
  the additional requirement that
  \[
  q_{ij}(1) =0 \text{\ \ if\ }\ \ i \geq 2m+1, \ \text{odd}
  \]
\end{defi}

Next we define the class of functions of $R$:

\begin{defi}
  $S^m(R^k (\ln R)^\ell)$ is the class of analytic functions
  $v:[0,\infty) \to \R$ with the following properties:

  (i) $v$ vanishes of order $m$ at $R=0$

  (ii) $v$ has a convergent expansion near $R=\infty$,
  \[
  v = \sum_{0 \leq j \leq \ell+i} c_{ij}\, R^{k-2i} (\ln R)^{j}
  \]
\end{defi}

We also introduce another auxiliary variable,
\begin{equation}\label{eq:bdef} b = \frac{(\ln(2+ R^2))^2}{(t \lambda)^2}
\end{equation}
Since we seek solutions inside the cone we can restrict $b$ to a
small interval $[0,b_0]$. We combine these three components in order
to obtain the full function class which we need:

\begin{defi}
  a) $S^m(R^k (\ln R)^\ell,\mathcal Q_n)$ is the class of analytic
  functions $v:[0,\infty) \times [0,1]\times [0,b_0] \to \R$ so that

  (i) $v$ is analytic as a function of $R,b$,
  \[
  v: [0,\infty) \times [0,b_0] \to \mathcal Q_n
  \]

  (ii) $v$ vanishes of order $m$ at $R=0$

  (iii) $v$ has a convergent expansion at $R=\infty$,
  \[
  v(R,\cdot,b) = \sum_{0 \leq j \leq \ell+i} c_{ij}(\cdot,b) R^{k-2i}
  (\ln R)^{j}
  \]
  where the coefficients $c_{ij}: [0,b_0] \to \cQ_m$ are analytic with
  respect to $b$

  b) $\IS^m(R^k (\ln R)^\ell,\mathcal Q_n)$ is the class of analytic
  functions $w$ on the cone $\calC_0$ which can be represented as
  \[
  w(r,t) = v(R,a,b), \qquad v \in S^m(R^k (\ln R)^\ell,\mathcal Q_n)
  \]
\end{defi}

We note that the representation of functions on the cone as in part
(b) is in general not unique since $R,a,b$ are dependent variables.
Later we shall exploit this fact and switch from one representation
to another as needed. We shall prove by induction that the
successive corrections $v_k$ and the corresponding error terms $e_k$
can be chosen with the following properties: For each $k\ge1$,

\begin{equation}
  v_{2k-1} \in \frac{1}{(t \lambda)^{2k}} \IS^3(R (\ln
  R)^{2k-1},\cQ_{k-1}) \label{v2k-1}\end{equation}
\begin{equation}
  t^2 e_{2k-1} \in \frac{1}{(t \lambda)^{2k}} \IS^1(R (\ln
  R)^{2k-1},\cQ'_{k-1}) \label{e2k-1}\end{equation}
\begin{equation}
  v_{2k} \in \frac{1}{(t \lambda)^{2k+2}} \IS^3(R^3 (\ln
  R)^{2k-1},\cQ_{k}) \label{v2k}\end{equation}
\begin{equation}
  t^2 e_{2k} \in \frac{1}{(t \lambda)^{2k}} \big[\IS^1(R^{-1} (\ln
  R)^{2k},\cQ_k)   + b \IS^1(R (\ln R)^{2k-1},\cQ'_{k})    \big]
  \label{e2k}\end{equation} Moreover, for $k=0$,
\begin{equation}\label{e20}
  t^2 e_0 \in S^1(R^{-1})
\end{equation}

\medskip

{\bf Step 0:} {\em We begin the analysis at $k=0$, where we explicitly
  compute $e_0$.}

\noindent We have
\begin{eqnarray*}
  e_0 &=& -u_{0tt}
  \\ &=& -  |\lambda'(t)|^2 r^2  Q''(R) - \lambda''(t)r Q'(R)
  \\ &=& -\left(\frac{\lambda'}{\lambda}\right)^2 R^2 Q''(R) -
  \frac{\lambda''}{\lambda} RQ'(R)
  \\ &=&   \frac{1}{t^2}\left((\nu+1)^2 \frac{4R^3}{(1+R^2)^2} -
    (\nu+1)(\nu+2) \frac{2R}{1+R^2}\right)
  \\ &=&   \frac{1}{t^2}\left(-(\nu+1)^2 \frac{4R}{(1+R^2)^2} +
    \nu(\nu+1) \frac{2R}{1+R^2}\right)
\end{eqnarray*}
With our notations
\[
t^2 e_0 \in S^1(R^{-1})
\]
as claimed. It remains to complete the induction step. Hence we
assume we know the above relations hold up to $k-1$ with $k\ge1$,
and construct $v_{2k-1}$, respectively $v_{2k}$, so that they hold
for the index $k$.

\medskip

{\bf Step 1:} {\em Begin with $e_{2k-2}$ satisfying \eqref{e2k}
  or~\eqref{e20} and choose $v_{2k-1}$ so that \eqref{v2k-1} holds.}

\smallskip
\noindent If $k=1$, then define $e_0^0:= e_0$. If $k>1$, we define the
principal part $e_{2k-2}^0$ of $e_{2k-2}$ by setting $b=0$, i.e.,
\[
e_{2k-2}^0(R,a) := e_{2k-2}(R,a,0)
\]
For the difference we can pull out a factor of $b$ and conclude that
\begin{align*}
  t^2e_{2k-2}^1 &\in \frac{b}{(t \lambda)^{2k-2}}\big[\IS^1(R^{-1} (\ln
  R)^{2k-2},\cQ_{k-1})+  \IS^1(R (\ln R)^{2k-3},\cQ_{k-1}')\big]\\
  &\subset \frac{1}{(t \lambda)^{2k}} \IS^1(R (\ln
  R)^{2k-1},\cQ_{k-1}')
\end{align*}
which can be included in $e_{2k-1}$, cf.~\eqref{e2k-1}.

We define $v_{2k-1}$ as in \eqref{vkodd} neglecting the $a$ dependence
of $e_{2k-2}^0$. In other words, $a$ is treated as a parameter.
Changing variables to $R$ in \eqref{vkodd} we need to solve the
equation
\[
(t \lambda)^2 L v_{2k-1} = t^2 e_{2k-2}^0 \in \frac{1}{(t
  \lambda)^{2k}} \IS^1(R^{-1} (\ln R)^{2k-2} , \cQ_{k-1} )
\]
where the operator $L$ is given by
\[
L = \partial_R^2 +\frac{1}R \partial_R - \frac{\cos(2u_0)}{R^{2}} =
\partial_R^2 +\frac{1}R \partial_R - \frac{1}{R^2} \frac{1-6R^2
  +R^4}{(1+R^2)^2}
\]
Then \eqref{v2k-1} is a consequence of the following ODE lemma.

\begin{lemma} Let $k\ge1$. Then the solution $v$ to the equation
  \[
  L v= f \in S^1(R^{-1} (\ln R)^{2k-2}), \qquad v(0) = v'(0)=0
  \]
  has the regularity
  \[
  v \in S^3(R (\ln R)^{2k-1})
  \]
\end{lemma}

\begin{proof}
  Since $f$ is analytic at $0$ with a linear leading term, one can
  easily write down a Taylor series for $v$ at $0$ with a cubic
  leading term.

  It remains to determine the asymptotic behavior of $v$ at infinity.
  For this it is convenient to remove the first order derivative in
  $L$ (to achieve constancy of the Wronskian). Thus, we seek a
  solution of
  \[ \tilde L \sqrt{R}\, v= \sqrt{R} f, \qquad \tilde L= \partial_R^2 -
  \frac{3}{4R^2} + \frac{8}{(1+R^2)^2}
  \]
  We use this fundamental system of solutions for $\tilde L$:
  \[
  \phi(R) = \frac{R^{\frac32}}{1+R^2}, \qquad \theta(R) =
  \frac{-1+4R^2 \ln R +R^4}{\sqrt{R}\,(1+R^2)}
  \]
  Their Wronskian is $W(\phi,\theta)=2$.  This allows us to obtain an
  integral representation for $v$ using the variation of parameters
  formula, which gives
  \[
  v = \frac12 R^{-\frac12}\theta(R)\int_0^R \phi(R')\sqrt{R'} f(R')\,
  \,dR' - \frac12 R^{-\frac12} \phi(R) \int_0^R \theta(R') \sqrt{R'}
  f(R')\, \,dR'
  \]
  Carrying out the integration shows that the right-hand side grows
  like $R(\ln R)^{2k-1}$ as claimed.
\end{proof}

As a special case of the above computation we also note the
representation for $v_1$,
\begin{equation}\label{eq:v1}
  v_1 = \frac{1}{(t \lambda)^{2}} V(R), \qquad V \in S^3(R \ln R)
\end{equation}

\smallskip

{\bf Step 2:} {\em Show that if $ v_{2k-1}$ is chosen as above then
  \eqref{e2k-1} holds.}

\medskip \noindent Thinking of $v_{2k-1}$ as a function of $t$, $R$
and $a$ we can write $e_{2k-1}$ in the form
\[
e_{2k-1} = N_{2k-1}(v_{2k-1}) + E^t v_{2k-1} + E^a v_{2k-1}
\]
Here $N_{2k-1}(v_{2k-1})$ accounts for the contribution from the
nonlinearity and is given by \eqref{eodd}. $E^t v_{2k-1}$ contains the
terms in
\[
\partial_t^2 v_{2k-1} (t,R,a)
\]
where no derivative applies to the variable $a$, while $ E^a v_{2k-1}$
contains the terms in
\[
(-\partial_t^2 + \partial_r^2 +\frac{1}r \partial_r) v_{2k-1} (t,R,a)
\]
where at least one derivative applies to the variable $a$. We begin
with the terms in $ N_{2k-1}$. We first note that, by summing the
$v_j$ over $1\le j\le 2k-2$,
\begin{equation}\label{eq:uentw1}
  u_{2k-2} - u_0 \in \frac{1}{(t \lambda)^{2}} \IS^3(R \ln R,\mathcal
  Q_{k-1})\end{equation}
To switch to trigonometric functions we need

\begin{lemma}
  Let
  \[
  v \in \frac{1}{(t \lambda)^{2}} \IS^3(R \ln R,\mathcal Q_{k-1})
  \]
  Then
  \[
  \sin v \in \frac{1}{(t \lambda)^{2}} \IS^3(R \ln R,\mathcal
  Q_{k-1}), \qquad \cos v \in \IS^0(1,\mathcal Q_{k-1})
  \]
  \label{sincosv}\end{lemma}

\begin{proof}
  We write
  \[
  \sin v = v g(v^2)
  \]
  with $g$ an entire function. Then it suffices to show that $g(v^2)
  \in \IS^0(1,\mathcal Q_{k-1})$. We begin with
  \[
  v^2 \in \frac{1}{(t \lambda)^{4}} \IS^6(R^2 (\ln R)^2,\mathcal
  Q_{k-1}) \subset \frac{R^2}{(t \lambda)^{2}} \frac{1}{(t
    \lambda)^{2}} \IS^4( (\ln R)^2,\mathcal Q_{k-1})
  \]
  But $a^2 = R^2 (t \lambda)^{-2}$ and the remaining ${(t
    \lambda)^{-2}}$ together with up to two $\log R$ factors combines
  to give one $b$ factor. We conclude that
  \[
  v^2 \in a^2 b\IS^2(1,\mathcal Q_{k-1}) \subset \IS^2(1,\mathcal
  Q_{k-1})
  \]
  For $w \in S^2(1,\mathcal Q_{k-1})$ we evaluate $g(w)$. Since $g$ is
  analytic we conclude that $g(w)$ is analytic in $R,b$ when
  interpreted as
  \[
  g(w)\;:\: [0,\infty) \times [0,b_0] \to \mathcal Q_{k-1}
  \]
  We consider the asymptotic expansion at $R=\infty$. Since we have an
  absolutely convergent asymptotic expansion for $w$ and a convergent
  Taylor series for $g$ at $0$, we obtain an absolutely convergent
  asymptotic expansion for $g(w)$. This gives
  \[
  g(w) \in S^0(1, \mathcal Q_{k-1})
  \]
  and concludes the proof of the lemma.
\end{proof}

Using Lemma~\ref{sincosv} and~\eqref{eq:uentw1} we compute
\begin{eqnarray*}
  \cos(2u_0) -\cos(2u_{2k-2}) &=& 2   \cos(2u_0)
  \sin^2(u_{2k-2} - u_0) + 2 \sin(2u_0) \sin(u_{2k-2} - u_0)
  \cos(u_{2k-2} - u_0)
  \\ &\subset& \frac{1}{(t \lambda)^{4}}
  \IS^6(R^2 (\ln R)^2,\mathcal Q_{k-1})+ \frac{1}{(t \lambda)^{2}}
  \IS^4(\ln R,\mathcal Q_{k-1})
\end{eqnarray*}
Hence
\begin{eqnarray*}
  && t^2 \frac{\cos(2u_0) -\cos(2 u_{2k-2})}{r^2} v_{2k-1} \\ &\in
  &\frac{(t \lambda)^{2}}{R^2}
  \left(  \frac{1}{(t \lambda)^{2}}\IS^4( \ln R,\mathcal Q_{k-1}) +
    \frac{1}{(t \lambda)^{4}} \IS^6( R^2(\ln R)^2,\mathcal
    Q_{k-1})\right) \frac{1}{(t \lambda)^{2k}} \IS^3(R (\ln
  R)^{2k-1},\mathcal Q_{k-1})
  \\ &\subset &
  \frac{1}{(t \lambda)^{2k}}   \left( \IS^5(R^{-1} (\ln R)^{2k},\mathcal Q_{k-1}) +
    \frac{1}{(t \lambda)^{2}}
    \IS^7(R (\ln R)^{2k+1},\mathcal Q_{k-1})\right)   \\
    &\subset &  \frac{1}{(t \lambda)^{2k}} \IS^5(R (\ln R)^{2k-1},\mathcal Q_{k-1})
\end{eqnarray*}
where at the last step we have pulled a $b$ factor out of the second
term.  Similarly we have
\begin{eqnarray*}
  && t^2 \frac{\sin(2u_{2k-2})}{2r^{2}} (1-\cos(2v_{2k-1})) \\  &\in &
  \frac{(t \lambda)^{2}}{R^2} \left(\IS^1(R^{-1},\cQ_{k-1}) + \frac{1}{(t
      \lambda)^{2}} \IS^3(R \ln R,\mathcal Q_{k-1})\right)
  \left(\frac{1}{(t \lambda)^{2k}} \IS^3(R (\ln R)^{2k-1},\mathcal Q_{k-1})   \right)^2
  \\ &=& \frac{1}{(t \lambda)^{2k}} \left(\frac{1}{(t \lambda)^{2k-2}}
    \IS^5(R^{-1}( \ln R)^{4k-2},\mathcal Q_{k-1})
    + \frac{1}{(t \lambda)^{2k}}
    \IS^7(R (\ln R)^{4k-1},\mathcal Q_{k-1})   \right)
  \\ & \subset& \frac{1}{(t \lambda)^{2k}}
  \IS^5(R (\ln R)^{2k-1},\mathcal Q_{k-1})
\end{eqnarray*}
where we have used a power of $b^k$ to pass to the final inclusion.
Finally,
\begin{eqnarray*}
  t^2  \frac{  \cos(2u_{2k-2})}{r^{2}}(2 v_{2k-1} - \sin( 2 v_{2k-1}))
  &\in& \frac{(t \lambda)^{2}}{R^2} \IS^0(1,\mathcal Q_{k-1} )
  \left(\frac{1}{(t \lambda)^{2k}} \IS^3(R (\ln R)^{2k-1},\mathcal
    Q_{k-1})   \right)^3
  \\ &\subset&  \frac{1}{(t \lambda)^{6k-2}} \IS^7(R
  (\ln R)^{6k-3},\mathcal Q_{k-1})
  \\ &\subset& \frac{1}{(t \lambda)^{2k}}\IS^7(R
  (\ln R)^{2k-1},\mathcal Q_{k-1})
\end{eqnarray*}
This concludes the analysis of $N_{2k-1}(v_{2k-1})$. We continue with
the terms in $E^t v_{2k-1}$, where we can neglect the $a$ dependence.
Therefore, it suffices to compute
\[
t^2 \partial_t^2 \left( \frac{1}{(t \lambda)^{2k}} \IS^3(R (\ln
  R)^{2k-1}) \right) \subset \frac{1}{(t \lambda)^{2k}} \IS^1(R (\ln
R)^{2k-1})
\]
Finally, we consider the terms in $E^a v_{2k-1}$. For
\[
v_{2k-1}(r,t) = \frac{1}{(t \lambda)^{2k}} w(R,a), \qquad w \in S^3(R
(\ln R)^{2k-1},\mathcal Q_{k-1})
\]
we have
\begin{eqnarray*}
  t^2 E^a v_{2k-1}& =&   \frac{1}{(t \lambda)^{2k}}  \big[ 2k\nu a
  w_a(R,a) - (\nu+1) R a w_{Ra}(R,a) + 2R^{-1} a^{-1} w_{Ra}(R,a) +
  a^{-1} w_a(R,a) \\ &+& (1-a^2)w_{aa}(R,a)  - a w_a(R,a) \big]
\end{eqnarray*}
Since $\mathcal Q_{k-1}$ are even in $a$ we conclude that
\[
a \partial_a, a^{-1} \partial_a, (1-a^2) \partial_a^2: \mathcal
Q_{k-1} \to \mathcal Q'_{k-1}
\]
Also the $R^{-1}$ factor simply removes one order of vanishing at
$R=0$. Hence we easily obtain
\[
t^2 E^a v_{2k-1} \in \frac{1}{(t \lambda)^{2k}} \IS^1(R (\ln
R)^{2k-1},\mathcal Q'_{k-1})
\]
This concludes the proof of \eqref{e2k-1}. We remark that for the
special case of $k=1$, i.e., with $v_1$ as in~\eqref{eq:v1}, these
arguments yield
\begin{equation}
  \label{eq:e1_def}t^2 e_1 \in (t\lambda)^{-2} S^3(R\log R)
\end{equation}

\medskip

{\bf Step 3:} {\em Define $v_{2k}$ so that \eqref{v2k} holds.}

\smallskip\noindent We begin the analysis with $e_{2k-1}$ replaced by
its main asymptotic component $f_{2k-1}$ at $R = \infty$ for $b=0$.
This has the form
\[
t^2 f_{2k-1} = \frac{R}{(t \lambda)^{2k}} \sum_{j=0}^{2k-1} q_j(a)
(\ln R)^j, \qquad q_j \in \mathcal Q'_{k-1}
\]
which we rewrite as
\[
t^2 f_{2k-1} = \frac{1}{(t \lambda)^{2k-1}} \sum_{j=0}^{2k-1} a q_j(a)
(\ln R)^j
\]
We remark that \eqref{eq:e1_def} implies that $t^2 f_1(a)=
(t\lambda)^{-1} a\log R$.  Consider the equation \eqref{vkeven} with
$ f_{2k-1}$ on the right-hand side,
\[
t^2\left(-\partial_t^2+\partial_r^2 +\frac{1}r \partial_r -
  \frac1{r^{2}}\right) w_{2k} = t^2 f_{2k-1}
\]
Homogeneity considerations suggest that we should look for a
solution $w_{2k}$ which has the form
\[
w_{2k} = \frac{1}{(t \lambda)^{2k-1}} \sum_{j=0}^{2k-1} W_{2k}^j(a)
(\ln R)^j
\]
The one-dimensional equations for $ W_{2k}^j$ are obtained by
matching the powers of $\ln R$. This gives the system of equations
\[
\begin{split}
  t^2 \left(-\partial_t^2+\partial_r^2 +\frac{1}r \partial_r -
    \frac1{r^{2}}\right) \left(\frac{1}{(t \lambda)^{2k-1}}
    W_{2k}^{j}(a)\right) = \frac{1}{(t \lambda)^{2k-1}} (a q_{j}(a)
  -F_j(a))
\end{split}
\]
where
\begin{equation}\label{eq:Fj_def}\begin{split}
F_j(a) &= (j+1)\left[((\nu+1)\nu(2k-1)+a^{-2}) W_{2k}^{j+1} +(a^{-1}
-
  (1+\nu) a) \partial_a W_{2k}^{j+1}\right]\\
  &\qquad  + (j+2)(j+1) ((\nu+1)^2 +
a^{-2}) W_{2k}^{j+2}
\end{split}
\end{equation}
Here we make the convention that $W_{2k}^{j} = 0$ for $j \geq 2k$.
Then we solve the equations in this system successively for decreasing
values of $j$ from $2k-1$ to $0$.

Conjugating out the power of $t$ we get
\[
t^2 \Big( -\Big(\partial_t+\frac{(2k-1)\nu}t\Big)^2+
  \partial_r^2 + \frac{1}r \partial_r - \frac{1}{r^2}
\Big)W_{2k}^{j}(a) = a q_{j}(a) -F_j(a)
\]
which we rewrite as an equation in the $a$ variable,
\begin{equation}
  L_{(2k-1)\nu} W_{2k}^{j} = a q_{j}(a) -F_j(a)
  \label{w2kj}
  \end{equation}
  where the one-parameter family of
operators $L_{\beta}$ is defined by
\begin{equation}\label{eq:Lbeta_def}
L_\beta =(1-a^2) \partial_a^2 + (a^{-1} + 2a \beta - 2a)
\partial_a + (-\beta^2 + \beta-a^{-2})
\end{equation}
We claim that solving this system with $0$ Cauchy data at $a=0$ yields
solutions which satisfy
\begin{equation}
  W_{2k}^j \in a^3 \cQ_k, \qquad j=\overline{0,2k-1}
  \label{cqk}\end{equation}

To prove this we need the following

\begin{lemma}
  Let $0 \leq m(j) \lesssim j^2$. Let $f$ be an analytic function in
  $[0,1)$ with an odd expansion at $0$ and an absolutely convergent
  expansion near $a = 1$ of the form
  \begin{equation}\label{eq:f_exp}\begin{split}
  f(a) &= f_0(a) + \sum_{j=1}^\infty \Bigg[ (1-a)^{(2j-1) \nu -\frac12} \sum_{m=0}^{m(2j-1)}
    f_{2j-1,m}(a)\, [\ln(1-a)]^m +\\
    &\qquad \qquad\qquad + (1-a)^{2j \nu}\sum_{m=0}^{m(2j)} f_{2j,m}(a)\, [\ln(1-a)]^m
    \Bigg]
  \end{split}\end{equation}
  with $f_{i,j}$ analytic near $a = 1$.
  Then there is a unique solution $w$ to the equation
  \begin{equation}
    L_{(2k-1)\nu}\, w = f, \qquad w(0) = 0, \ \partial_a w(0) = 0
  \end{equation}
  with the following properties:

  (i) $w$ is an analytic function in $[0,1)$

  (ii) $w$ is cubic at $0$ and has an odd expansion at $0$

  (iii) $w$ has an absolutely convergent expansion near $a = 1$ of the
  form
  \begin{equation}\label{eq:w_exp}
  \begin{split}
  w(a) &= w_0(a) + \sum_{j=1}^\infty \Bigg[ (1-a)^{(2j-1) \nu +\frac12} \sum_{\ell=0}^{\ell(2j-1)}
    w_{2j-1,\ell}(a) [\ln(1-a)]^\ell +\\
    &\qquad\qquad\qquad + (1-a)^{2j \nu+1}\sum_{\ell=0}^{\ell(2j)} w_{2j,\ell}(a)
    [\ln(1-a)]^\ell
  \Bigg]
  \end{split}
  \end{equation}
  with $w_{i,j}$ analytic near $a = 1$ and $\ell(i)=m(i)$ with one
  exception, namely $\ell(2k-1)=m(2k-1)+1$. If however
  $f_{2k-1,m(2k-1)}(1)=0$, then $\ell(2k-1)=m(2k-1)$. In that case
  also $w_{2k-1,\ell}(1)=0$ if $\ell>0$, but not necessarily when $\ell=0$.
  Finally, if $f_{2i-1,j}(1)=0$ for all $i>k$ and all $j$, then also
  $w_{2i-1,\ell}(1)=0$ for all $i>k$ and all $\ell$.
\end{lemma}

\begin{proof}
  Denote $\beta= (2k-1)\nu$. Since $k\ge1$ and $\nu>\frac12$, also
  $\beta>\frac12$.
  Cearly, $L_{\beta}$ is the sum of two
  pieces: that part which is homogeneous of degree $-2$ in~$a$, viz.,
  \[
  L_\beta^0 = \partial_a^2 + a^{-1} \partial_a - a^{-2} = a^{-1}
  \partial_a (a\partial_a) - a^{-2}
  \]
  and the remainder which is homogeneous of degree $0$.
  The equation $L_\beta^0\, y=0$  has fundamental solutions $a$ and $a^{-1}$.
  A standard power-series ansatz then leads to
  fundamental
  solutions of $L_\beta\, y=0$ of the form
  \[
  \phi_1(a) = a(1+O(a^2)), \qquad \phi_2(a) = a^{-1}(1 + O(a^2))
  \]
  where the $O(\cdot)$ terms are analytic functions of $a^2$. Since our right-hand side $f$ has size $O(a)$ at
  $0$, this implies that we can use the equation to write a complete
  Taylor expansion for $w$ near $0$. Matching coefficients in
  $L_\beta \, w=f$ with
  \[
  f(a) = \sum_{j=1}^\infty f_j\, a^{2j-1}, \qquad w(a) = \sum_{j=2}^\infty w_j \, a^{2j-1}
  \]
  yields
  \[
  2j(2j-2) w_j =(2j(2j-1) - (4j-1) \beta +\beta^2) w_{j-1} +
  f_{j-1}
  \]
  where we take $w_1=0$. The coefficient of $w_j$ is always nonzero;
  this allows us to successively compute the coefficients $w_j$. The
  convergence of the series for $w$ easily follows from the
  convergence of the series for $f$.

  \noindent It remains to study the solution $w$ near $a=1$. The behavior of
  $L_{\beta}$ at $1$ is well approximated by
  \[
  L_\beta^1 = 2(1-a) \partial_a^2 + (2\beta -1) \partial_a +
  (\beta -\beta^2 -1) = 2(1-a)^{\beta+\frac12} \partial_a
  \big[(1-a)^{-\beta+\frac12} \partial_a\big] + (\beta -\beta^2 -1)
  \]
in the sense that
\begin{equation}\label{eq:Lzerl}
L_\beta =  L_\beta^1 + (a-1)\big[ (1-a)\partial_a^2 +
(2(\beta-1)-a^{-1}) \partial_a + (a+1)a^{-2}\big] =: L_\beta^1 +
(a-1) L_\beta^2
\end{equation}
The differential operator \begin{equation}   2(1-a)^{\beta+\frac12}
\partial_a
  \big[(1-a)^{-\beta+\frac12} \partial_a\big]\label{eq:hom_diffop}
\end{equation}
annihilates
   $1$ and $(1-a)^{\beta+\frac12}$. Therefore, we seek a fundamental system
    for $L_\beta^1\, y=0$ of the form
  \begin{equation}\label{eq:fundsys}
  \phi_1(a) = 1+\sum_{\ell=1}^\infty \mu_\ell (1-a)^\ell , \qquad \phi_2(a) =
  (1-a)^{\beta+\frac12}\Big[ 1+\sum_{\ell=1}^\infty \tilde \mu_\ell (1-a)^\ell \Big]
  \end{equation}
This leads to the conditions, for $\ell\ge1$,
\begin{align}
  & \mu_1(1-2\beta)+\beta -\beta^2 -1 =0,\quad \mu_{\ell+1}
  (\ell+1)(2\ell+1-2\beta) + (\beta -\beta^2 -1)\mu_\ell=0 \label{eq:phi1def}\\
& \tilde\mu_1(2\beta+3)+\beta -\beta^2 -1 =0,\quad
\tilde\mu_{\ell+1} (\ell+1)(2\ell+3+2\beta) + (\beta -\beta^2
-1)\tilde\mu_\ell =0 \label{eq:phi2def}
\end{align}
Clearly, \eqref{eq:phi2def} always has a solution whereas
\eqref{eq:phi1def} requires $\beta-\half\not\in \Z^+$; in the latter
case, the series in~\eqref{eq:fundsys} define entire functions. If,
on the other hand, $\ell_0:=\beta-\half\in \Z^+$, then
  $\phi_1$ is modified to
  \begin{equation}\label{eq:phi1mod}
  \phi_1(a) = 1+\sum_{\ell=1}^\infty \mu_\ell (1-a)^\ell + c_1\,\phi_2(a) \ln(1-a)
  \end{equation}
with the unique choice
$c_1=-(2\beta+1)^{-1}(\beta-\beta^2-1)\mu_{\ell_0}$. Here
\eqref{eq:phi1def} is unchanged and can be solved for $\mu_\ell$ up
to $\ell\le\ell_0$; for $\ell>\ell_0$ this equation is then modified
by the terms from the $\phi_2$ series (in particular, for
$\ell=\ell_0+1$ the choice of $c_1$ assures the validity of the
equation, whereas for all  $\ell>\ell_0+1$ we can again solve for
$\mu_{\ell}$). Finally, observe that the same process also leads to
a fundamental system for $L_\beta$; indeed, the remainder
$(a-1)L_\beta^2$ in~\eqref{eq:Lzerl} does not change the
coefficients of $\mu_{\ell+1}$ or $\tilde\mu_{\ell+1}$
in~\eqref{eq:phi1def} and~\eqref{eq:phi2def}.  In conclusion, the
preceding power series argument leads to a fundamental system of
$L_\beta\, y=0$, which we again denote by $\phi_1(a)$ and
$\phi_2(a)$.

Modulo a linear combination of $\phi_1, \phi_2$ it suffices to
  find one solution to the inhomogeneous equation $L_\beta\, w=f$ near $a=1$.
At this point, it will be convenient to write $L_\beta$ as a
Sturm-Liouville operator. Thus, we write
\[
L_\beta = q_1^{-1}(a) \partial_a (q_2(a) \partial_a) + q_3(a)
\]
with, cf.~\eqref{eq:Lbeta_def},
\[
q_1^{-1} q_2(a) = 1-a^2, \qquad q_1^{-1} q_2'(a) = a^{-1} +
2a(\beta-1), \qquad q_3(a) = -\beta^2 + \beta-a^{-2}
\]
One checks that for $a$ close to $1$ the first two equations admit
solutions
\[ q_2(a) = (1-a)^{-\beta+\frac12}[1+(1-a)\tilde q_1(a)],\qquad q_1(a) = \frac12(1-a)^{-\beta-
\frac12}[1+(1-a)\tilde q_2(a)]
\]
with $\tilde q_1,\tilde q_2$ analytic near $a=1$. The Wronskian can
now be computed as
  \[
q_2(a) [\phi_1(a)\phi_2'(a)-\phi_1'(a)\phi_2(a)]=\beta+1/2
  \]
Thus, a particular solution of the inhomogeneous problem is given by
\begin{equation}
\label{eq:sol_form}
 w(a) = (\beta+1/2)^{-1}\phi_1(a)\int_a^1
\phi_2(a')q_1(a') f(a')\, da'+(\beta+1/2)^{-1}\phi_2(a)\int_{a_0}^a
\phi_1(a') q_1(a') f(a')\, da'
\end{equation}
where $a_0<1$ is some number close to $1$. For the first integral,
note that $\phi_2(a')q_1(a')$ is an analytic function in the
neighborhood of $a'=1$. Let $\gamma\ne-1$ and  $m$ be a positive
integer. Iterating the relation
\begin{equation}\label{eq:log_int}
 \int_a^1 (1-a')^\gamma [\log(1-a')]^m\,da'= \frac{1}{\gamma+1}
[\log(1-a)]^m (1-a)^{\gamma+1} - \frac{m}{\gamma+1} \int_a^1
(1-a')^\gamma [\log(1-a')]^{m-1}\, da'
\end{equation}
shows that each summand on the right-hand side of \eqref{eq:f_exp},
inserted into the first integral in~\eqref{eq:sol_form}, makes an
admissible contribution to $w$ in the sense of~\eqref{eq:w_exp}
 (for this it does not matter whether $\phi$ takes the form~\eqref{eq:fundsys}
 or~\eqref{eq:phi1mod}).
The analysis of the second integral in~\eqref{eq:sol_form} is again
based on~\eqref{eq:log_int} provided $j\ne k$, since then
$\gamma\ne-1$. If $j=k$, then we encounter
\[
\int_{a_0}^a (1-a')^{-1} [\log(1-a')]^m\, da' = -(m+1)^{-1}
[\log(1-a)]^{m+1} +C
\]
which explains why we obtain an extra log-factor when $j=k$.
Clearly, if $f_{2k-1,m(2k-1)}(1)=0$ then there is no extra
log-factor and the lemma is proved. In that case also we write, with
\[f= (1-a)^{\beta -\frac12} \sum_{m=0}^{m(2k-1)}
    f_{2k-1,m}(a)\, [\ln(1-a)]^m\]
    the second integral in \eqref{eq:sol_form} as
\[
\phi_2(a)\int_{a_0}^a \phi_1(a') q_1(a') f(a')\, da' =
\phi_2(a)\int_{a_0}^1 \phi_1(a') q_1(a') f(a')\, da' -
\phi_2(a)\int_{a}^1 \phi_1(a') q_1(a') f(a')\, da'
\]
The first term on the right-hand side here is just a multiple of
$\phi_2(a)$, whereas the second one possesses the extra vanishing at
$a=1$, as claimed. The final claim of the lemma follows similarly.
\end{proof}

Before turning to the proof of~\eqref{cqk} in full generality, we
first discuss the special case $k=1$. This will also serve to
explain how the algebra $\cQ_k$ arises at all in the iteration. If
$k=1$,  then~\eqref{w2kj} reduces to the system
\[
L_\nu \, W_2^1(a)=a,\qquad L_\nu\, W_2^0(a) =
-(\nu(\nu+1)+a^{-2})W_2^1(a) - (a^{-1}-(\nu+1)a)\partial_a W_2^1(a)
\]
due to $t^2 f_1(a)= (t\lambda)^{-1} a\log R$.  In view of the
solution formula \eqref{eq:sol_form} with $\beta=\nu$, provided
$\nu-\frac12\not\in\Z^+$,
\[\begin{split}
W_2^1(a) &= g_0(a) +  g_1(a)(1-a)^{\nu+\frac12}\\
 W_2^0(a) &= h_0(a)
+ h_1(a)(1-a)^{\nu+\frac12} + h_2(a)(1-a)^{\nu+\frac12}\log(1-a)
\end{split}
\]
where $g_j(a), h_j(a)$ are analytic around $a=1$.  Note that the
term $(1-a)^{\nu+\frac12} \log(1-a)$ appears in~$W_2^0$ due to
$\partial_a W_2^1$. Similarly, if $\nu-\frac12\in\Z^+$, then
\begin{equation}\nn\begin{split}
W_2^1(a) &= g_0(a) + g_1(a)(1-a)^{\nu+\frac12}  +
g_2(a)(1-a)^{\nu+\frac12}\ln(1-a)\\
 W_2^0(a) &= h_0(a) + (1-a)^{\nu+\frac12}
\sum_{\ell=0}^2 h_{\ell+1}(a)[\log(1-a)]^\ell + (1-a)^{2\nu+1}
\sum_{\ell=0}^2 h_{\ell+4}(a)[\log(1-a)]^\ell,
\end{split}
\end{equation}
with analytic $g_j,h_j$.  The terms involving the $(1-a)^{2\nu+1}$
factor in $W_2^0$ are due to the modified $\phi$,
see~\eqref{eq:phi1mod}. Thus, we see that in all cases $W_2^j \in
\cQ_1$ for $j=0,1$ and $a$ near~$1$.

 We now continue with the proof of \eqref{cqk} for general $k$. At
first we consider the easier case when $\nu$ is irrational.  We
apply the lemma in \eqref{w2kj} using for the right-hand side the
fact that $q_{2k-1} \in Q'_{k-1}$. This implies that the coefficient
of $(1-a)^{(2k-1)\nu-\frac12}$ in $q_{2k-1} $ vanishes at $a=1$. The
lemma gives a similar expansion for $W_{2k}^{2k-1}$ with the
required vanishing conditions. Hence $W_{2k}^{2k-1}\in \cQ_{k}$,
with one extra $(1-a)^{(2k-1)\nu+\frac12}$ term (this is the
$w_{2k-1,0}(1)\ne0$ statement of the lemma) -- we refer to this as
the "free term" in what follows.

Next we reiterate the argument for the remaining $W_{2k}^j$ which
solve \eqref{w2kj}.  At each step we have to compute $F_j$,
see~\eqref{eq:Fj_def}. Since $W_{2k}^{j+1}$ and $W_{2k}^{j+2}$ have
an odd Taylor expansion at $0$ beginning with a cubic term, it
follows that $F_j$ has an odd Taylor expansion at $0$ beginning with
a linear term. The expansion of $F_j$ around $a=1$ is similar to the
one for $W_{2k}^{j+1}$ except that one $(1-a)$ factor is lost in the
"free term". For $j=2k-2$ this produces the term
$(1-a)^{(2k-1)\nu+\frac12}\log(1-a)$ in $W_{2k}^j$ etc. At the
conclusion of the iteration we have gained at most $2k-1$ logarithms
in the free term for the $W_{2k}^j$'s. Then \eqref{cqk} follows.

Next we consider the case when $\nu$ is rational. This is more
difficult since now the term $(1-a)^{(2k-1)\nu-\frac12}$ can also
arise in expressions of the form
\[
f_{2j-1,m}(a) (1-a)^{(2j-1) \nu -\frac12} [\ln(1-a)]^m \ \ or \ \
f_{2j,m}(a) (1-a)^{2j \nu} [\ln(1-a)]^m
\]
using the notations of the lemma. This leads to more logarithms than
in the irrational case. The first term above will be of interest if
$2(k-j)\nu$ is an integer, while the second needs to be considered
if $(2k-2j-1)\nu-\frac12$ is an  integer. The worst case is $j=k-1$.
Then we can have $m(2k-2)$ logarithms in the second term above,
while $2k$ more logarithms are produced by the $2k$ applications of
the lemma. Hence we need the relation
\[
m(j) \geq m(j-1) + j+1
\]
which is verified e.g.~by $m(j)= j^2$ (we pick $n_j=2j^2$ because of
$j=1$, see above).

\noindent We cannot use $w_{2k}$ for $v_{2k}$ due to the singularity
of $\ln R$ at $R=0$. However, we define instead
\[
v_{2k} := \frac{1}{(t \lambda)^{2k-1}} \sum_{j=0}^{2k-1} W_{2k}^j(a)
\Big(\frac12 \ln (1+R^2)\Big)^j
\]
In doing this we add an additional component to the error. This is
large near $R=0$, but this is not so important since the aim of this
correction is to improve the error for large $R$.  Since $a^3 = R^3/
(t \lambda)^{3}$, pulling a cubic factor $a^3$ out of the $W$'s it is
easy to see that \eqref{v2k} holds.

{\bf Step 4: } For $v_{2k}$ defined as above we show that the
corresponding error $e_{2k}$ satisfies \eqref{e2k}. We can write
$e_{2k}$ in the form
\[
t^2 e_{2k} = t^2 \left(e_{2k-1} - e_{2k-1}^0 \right) + t^2 \left(
  e_{2k-1}^0- \left(-\partial_t^2+\partial_r^2 +\frac{1}r \partial_r -
    \frac1{r^{2}}\right) v_{2k} \right) + t^2 N_{2k}(v_{2k})
\]
where $N_{2k}$ is defined by \eqref{eeven} and
\[
e_{2k-1}^0= \frac{R}{(t \lambda)^{2k}} \sum_{j=0}^{2k-1} q_j(a)
\Big(\frac12 \ln (1+R^2) \Big)^j
\]
We begin with the first term in $e_{2k}$, which has the form
\[
t^2(e_{2k-1} - e_{2k-1}^0) \in \frac{1}{(t \lambda)^{2k}}
\big[\IS^1(R^{-1} (\ln R)^{2k},\mathcal Q'_{k-1}) + b \IS^1(R (\ln
R)^{2k-1},\mathcal Q'_{k-1})\big]
\]
The second term is contained in the second term of \eqref{e2k}. It
remains to show that
\begin{equation}
  \IS^1(R^{-1} (\ln R)^{2k},\mathcal Q'_{k-1}) \subset \IS^1(R^{-1}
  (\ln R)^{2k},\mathcal Q_{k-1}) + b \IS^1(R (\ln R)^{2k-1},\mathcal
  Q'_{k-1}) \label{qqp}\end{equation} For $w \in \IS^1(R^{-1} (\ln
R)^{2k},\mathcal Q'_{k-1})$ we write
\[
w = (1-a^2)w + \frac{1}{(t \lambda)^2} R^{2} w
\]
Then
\[
(1-a^2)w \in \IS^1(R^{-1} (\ln R)^{2k},\mathcal Q_{k-1}), \qquad
\frac{1}{(t \lambda)^2} R^{2} w \in b \IS^1(R (\ln R)^{2k-1},\mathcal
Q'_{k-1})
\]
as desired. The second term in $e_{2k}$ would equal $0$ if we were
to replace $\frac12\ln(1+R^2)$ by $\ln R$ in both $e_{2k}^0$ and
$v_{2k}$. Hence the difference is obtained when we replace the
derivatives of $ \frac12\ln (1+R^2)$ by derivatives of $\ln R$ in
the expression
\[
t^2 \left(-\partial_t^2+\partial_r^2 +\frac{1}r \partial_r \right)
v_{2k} = t^2 \left(-\partial_t^2+\partial_r^2 +\frac{1}r \partial_r
\right) \left( \frac{1}{(t \lambda)^{2k-1}} \sum_{j=0}^{2k-1}
  W_{2k}^j(a) \Big(\frac12\ln (1+R^2)\Big)^j \right)
\]
Computing these differences one finds that the second term in $e_{2k}$
is a sum of expressions of the form
\[
\frac{1}{(t \lambda)^{2k-1}} \sum_{j=0}^{2k-1} \frac{
  W_{2k}^j(a)}{a^{2}} \Big[S(R^{-2}) (\ln (1+R^2))^{j-1} + S(R^{-2}) (\ln
(1+R^2))^{j-2}\Big] + \frac{\partial_a W_{2k}^j(a)}{a} S(R^{-2})
(\ln (1+R^2))^{j-1}
\]
Since $W_{2k}^j$ are cubic at $0$ it follows that we can pull out an
$a$ factor and see that this part of the error is in
\[
\frac{1}{(t \lambda)^{2k}} \IS^1(R^{-1} (\ln R)^{2k-2},\mathcal Q'_k)
\]
which is admissible by \eqref{qqp}.

Finally, we consider the nonlinear terms in $N_{2k}$.  Again the
$a,b$ dependence is uninteresting since $\mathcal Q_k$ is an
algebra. We shall use that
\[
u_{2k-1}-u_0 \in \frac{1}{(t\lambda)^2}IS^3 (R\log R, \cQ_k)
\]
 By Lemma~\ref{sincosv}, for the linear
term we therefore have
\begin{eqnarray*}
  && t^2 \frac{1 -\cos(2 u_{2k-1})}{r^2} v_{2k} \\
  &\in & \frac{(t
    \lambda)^2}{R^2} \left( IS^1(R^{-1},\cQ_k) + \frac{1}{(t \lambda)^{2}}
    \IS^3(R\ln
    R,\mathcal Q_k)\right)^2 \frac{1}{(t \lambda)^{2k+2}} \IS^3(R^3 (\ln R)^{2k-1},\mathcal Q_k)
  \\ &\subset&  \frac{1}{(t \lambda)^{2k}}
  \left( \IS^3(R^{-1} (\ln R)^{2k-1} ,\mathcal Q_k) +  \frac{1}{(t
      \lambda)^{2}} \IS^5(R
    (\ln R)^{2k},\mathcal Q_k) +  \frac{1}{(t \lambda)^{4}} \IS^7(R^3
    (\ln R)^{2k+1},\mathcal Q_k)\right)
  \\ &\subset&  \frac{1}{(t \lambda)^{2k}}
  \left( \IS^3(R^{-1} (\ln R)^{2k-1} ,\mathcal Q_k) +  \frac{b}{(t
      \lambda)^{2}} \IS^5(R
    (\ln R)^{2k-1},\mathcal Q_k) \right)
\end{eqnarray*}
For the quadratic term we obtain
\begin{eqnarray*}
  &&t^2 \frac{\sin(2u_{2k-1})}{2r^{2}} (1-\cos(2v_{2k})) \\
  &\in & \frac{(t \lambda)^2}{R^2} \left( \IS^1(R^{-1} ,\mathcal Q_k) + \frac{1}{(t \lambda)^{2}} \IS^3(R\ln
    R ,\mathcal Q_k)\right) \left( \frac{1}{(t \lambda)^{2k+2}} \IS^3(R^3 (\ln
    R)^{2k-1}) ,\mathcal Q_k)\right)^2
  \\ &\subset&  \frac{1}{(t \lambda)^{2k}}
  \left( \frac{1}{(t \lambda)^{2k+2}}  \IS^5(R^3  (\ln R)^{4k-2}
    ,\mathcal Q_k) +  \frac{1}{(t \lambda)^{2k+4}}  \IS^7(R^5  (\ln
    R)^{4k-1} ,\mathcal Q_k) \right)
  \\ &\subset&  \frac{1}{(t \lambda)^{2k}}
  \left(  \IS^1(R^{-1}  (\ln R)^{2k} ,\mathcal Q_k) +  b  \IS^3(R
    (\ln R)^{2k-1} ,\mathcal Q_k) \right)
\end{eqnarray*}
Finally, the cubic term is
\begin{eqnarray*}
  t^2 \frac{  \cos(2u_{2k-1})}{r^{2}}(2
  v_{2k} - \sin( 2 v_{2k})) &\in& \frac{(t \lambda)^2}{R^2} \left( \frac{1}{(t \lambda)^{2k+2}} \IS^3(R^3 (\ln
    R)^{2k-1}),\mathcal Q_k)\right)^3
  \\ &\subset&\frac{1}{(t \lambda)^{2k}} \frac{1}{(t \lambda)^{4k+4}} IS^7(R^7(\ln
  R)^{6k-3},\cQ_k))
  \\ &\subset&\frac{a^6 b^{4k-2}}{(t \lambda)^{2k}}  IS^1(R(\ln
  R)^{2k-1},\cQ_k))\\
&\subset&\frac{b}{(t \lambda)^{2k}}  IS^1(R(\ln
  R)^{2k-1},\cQ_k'))
\end{eqnarray*}
This concludes the proof of Theorem~\ref{thm:sec2}.
\end{proof}

\section{The perturbed equation}
\label{sec:perturb}

We now need to complement the approximate solution found in the first
section to an actual solution. The mechanism for achieving this will
rely on the construction of an approximate parametrix for a suitable
wave-type equation. We now set about deriving this equation: we seek
an exact solution of the form
\begin{equation}\nonumber
  u(t,r)=u_{2k-1}(t,r)+\epsilon(t,r)
\end{equation}
where $u_{2k-1}$ is as in the previous section and $\epsilon$ will
be obtained by means of an iteration procedure. To motivate this
procedure, note that we need to solve the following equation,
see~\eqref{eq:gleich},
\begin{equation}\label{big}\begin{split}
    &-\epsilon_{tt}+\epsilon_{rr}+\frac{1}{r}\epsilon_{r}-\frac{\cos(2Q(\lambda
      r))}{r^{2}}\epsilon= N_{2k-1}(\epsilon) + e_{2k-1}
  \end{split}\end{equation}
where $N_{2k-1}$ is defined in \eqref{eodd} but with $u_{2k-2}$
replaced by $u_{2k-1}$.

In order to remove the time dependence of the potential in
\eqref{big}, we now introduce new coordinates: first, the new time is
to satisfy the relation
\begin{equation}\nonumber
  \frac{\partial}{\partial\tau}=\frac{1}{\lambda(t)}\frac{\partial}{\partial t}
\end{equation}
Specifically, we may put
$\tau=-\int_{t}^{1}\lambda(s)ds+\frac{1}{\nu}=\frac{1}{\nu}t^{-\nu}$.
Thus, the singularity now corresponds to $\tau=\infty$.  Next,
introduce the new dependent variable $v(\tau, R):=\epsilon(t(\tau),
\lambda^{-1}R)$, where we now understand $\lambda$ as a function of
$\tau$. Then we have
\begin{equation}\nonumber
  \frac{\partial}{\partial\tau}v=t'(\tau)\epsilon_{t}(t(\tau), \lambda^{-1}R)-\frac{\lambda_{\tau}}{\lambda^{2}}
  R\epsilon_{r}(t(\tau), \lambda^{-1}R),\,\frac{\partial}{\partial R}v=\lambda^{-1}\epsilon_{r}(t(\tau), \lambda^{-1}R)
\end{equation}
This entails that
\begin{equation}\nonumber
  (\partial_{\tau}+\frac{\lambda_{\tau}}{\lambda}R\partial_{R})v=\lambda^{-1}\epsilon_{t}(t(\tau), \lambda^{-1}R)
\end{equation}
From here we get
\begin{equation}\nonumber
  (\partial_{\tau}+\frac{\lambda_{\tau}}{\lambda}R\partial_{R})^{2}v=(\partial_{\tau}+\frac{\lambda_{\tau}}{\lambda}R\partial_{R})
  [\lambda^{-1}\epsilon_{t}]=\lambda^{-2}\epsilon_{tt}-\frac{\lambda_{\tau}}{\lambda^{2}}\epsilon_{t}
  =\lambda^{-2}\epsilon_{tt}-\frac{\lambda_{\tau}}{\lambda}\partial_{\tau}v-[\frac{\lambda_{\tau}}{\lambda}]^{2}R\partial_{R}v
\end{equation}
We conclude that we may recast the wave equation~\eqref{big} in the
following way:
\begin{equation}\nonumber
  -[(\partial_{\tau}+\frac{\lambda_{\tau}}{\lambda}R\partial_{R})^{2}+\frac{\lambda_{\tau}}{\lambda}(\partial_{\tau}
  +\frac{\lambda_{\tau}}{\lambda}R\partial_{R})]v+\Big(\partial_{R}^{2}+\frac{1}{R}\partial_{R}-\frac{\cos[2Q(R)]}{R^{2}}\Big)v
  =\frac{1}{\lambda^{2}}[N_{2k-1}(\epsilon) + e_{2k-1}](t(\tau), \lambda^{-1}R)
\end{equation}

In order to turn the above second order elliptic operator in $R$ into
a selfadjoint operator we  introduce the new variable $\tilde{\epsilon}(\tau,
R):=R^{\frac{1}{2}}v(\tau, R)$. This leads to
\begin{equation}\nonumber
  (\partial_{\tau}+\frac{\lambda_{\tau}}{\lambda}R\partial_{R})\tileps=R^{\frac{1}{2}}(\partial_{\tau}
  +\frac{\lambda_{\tau}}{\lambda}R\partial_{R})v+\frac{1}{2}R^{\frac{1}{2}}\frac{\lambda_{\tau}}{\lambda}v(\tau, R)
\end{equation}
\begin{equation}\nonumber
  (\partial_{\tau}+\frac{\lambda_{\tau}}{\lambda}R\partial_{R})^{2}\tilde{\epsilon}=
  R^{\frac{1}{2}}(\partial_{\tau}+\frac{\lambda_{\tau}}{\lambda}R\partial_{R})^{2}v+R^{\frac{1}{2}}
  \frac{\lambda_{\tau}}{\lambda}(\partial_{\tau}+\frac{\lambda_{\tau}}{\lambda}R\partial_{R})v+\frac{1}{2}R^{\frac{1}{2}}
  \partial_{\tau}(\frac{\lambda_{\tau}}{\lambda})v+\frac{1}{4}R^{\frac{1}{2}}(\frac{\lambda_{\tau}}{\lambda})^{2}v
\end{equation}
One checks that
\begin{equation}\nonumber
  R^{\frac{1}{2}}(\partial_{\tau}+\frac{\lambda_{\tau}}{\lambda}R\partial_{R})^{2}v+R^{\frac{1}{2}}
  \frac{\lambda_{\tau}}{\lambda}(\partial_{\tau}+\frac{\lambda_{\tau}}{\lambda}R\partial_{R})v
  =(\partial_{\tau}+\frac{\lambda_{\tau}}{\lambda}R\partial_{R})^{2}\tilde{\epsilon}-
  \frac{1}{4}(\frac{\lambda_{\tau}}{\lambda})^{2}\tilde{\epsilon}
  -\frac{1}{2}\partial_{\tau}(\frac{\lambda_{\tau}}{\lambda})\tilde{\epsilon}
\end{equation}
as well as
\begin{equation}\nonumber
  R^{\frac{1}{2}}\Big(\partial_{R}^{2}+\frac{1}{R}\partial_{R}-\frac{\cos[2Q(R)]}{R^{2}}\Big)v=
  (\partial_{R}^{2}-\frac{3}{4
    R^{2}}+\frac{8}{(1+R^{2})^{2}})\tilde{\epsilon}
\end{equation}
Combining these observations with \eqref{big}, we now obtain the wave
equation
\begin{equation}\label{equation2}
   \left(-(\partial_{\tau}+\frac{\lambda_{\tau}}{\lambda}R\partial_{R})^{2}
    +\frac{1}{4}(\frac{\lambda_{\tau}}{\lambda})^{2}
+\frac{1}{2}\partial_{\tau}(\frac{\lambda_{\tau}}{\lambda})\right) \tilde{\epsilon}
  - \calL \tilde{\epsilon} = \lambda^{-2} R^\frac12
    \left(  N_{2k-1}(R^{-\frac12} \tileps) + e_{2k-1}\right)
\end{equation}
where
\begin{equation}\label{eq:L_def}
  \calL:=-\partial_{R}^{2}+\frac{3}{4R^{2}}-\frac{8}{(1+R^{2})^{2}}
\end{equation}
Equation~\eqref{equation2} is the main equation which we need to solve
in this paper. As a first step,  in the following section we will carefully analyze the spectral
properties of the underlying linear operator $\calL$.

\section{The analysis of the underlying strongly singular Sturm-Liouville
  operator}
\label{sec:spec}

The goal of this section is to develop the scattering theory of
$\cL$ from~\eqref{eq:L_def}. We start with the basic\footnote{In
this section we use the variable $r>0$
  for the independent variable. The reader should note that this now
  plays the role of $R$ in the previous section. }
\begin{defi}
  Let
  \[ \cL_0 := -\frac{d^2}{\,dr^2} + \frac{3}{4r^2}, \qquad \cL :=
  \cL_0 - \frac{8}{(1+r^2)^2} =: \cL_0 + V\] be half-line operators on
  $L^2(0,\infty)$. They are self-adjoint with the same domain, namely
  \[ \dom(\cL)=\dom(\cL_0)=\{ f\in L^2((0,\infty))\::\: f,f'\in
  AC_{{\rm loc}}((0,\infty)),\; \cL_0 f \in L^2((0,\infty)) \}
  \]
\end{defi}

It is important to realize that because of the strong singularity of
the potential at $r=0$ no boundary condition is needed there to
insure self-adjointness. Technically speaking, this means that
$\cL_0$ and $\cL$ are in the {\em limit point case} at $r=0$, see
Gesztesy, Zinchenko~\cite{GZ}. It is worth noting that the potential
$\frac{3}{4r^2}$ is critical with respect to this property --- any
number smaller than $\frac34$ leads to an operator which is in the
limit circle case at $r=0$. We remark that $\cL_0$ and $\cL$ are in
the limit point case at $r=\infty$ by a standard criterion
(sub-quadratic growth of the potential).

\begin{lemma}
  \label{lem:spec} The spectrum of $\cL$ is purely
  absolutely continuous and equals $\spec(\cL)=[0,\infty)$.
\end{lemma}
\begin{proof}
  That $\cL$ has no
  negative spectrum follows from
\begin{equation}
  \label{eq:reson}
  \cL \phi_0 =0, \qquad \phi_0(r)= \frac{r^{3/2}}{1+r^2}
\end{equation}
with $\phi_0$ positive (by the Sturm oscillation theorem,
see~\cite{DS}). And since $\phi_0\not\in L^2((0,\infty))$, zero is
not an eigenvalue. The pure absolute continuity of the spectrum of
$\cL$ is an immediate consequence of the fact that the potential of
$\cL$ is integrable at infinity.
\end{proof}

Since $\phi_0 \not\in L^2((0,\infty))$, one refers to zero energy as
a {\em resonance}. Heuristically speaking, this notion can be
thought of as follows: by inspection, $\cL_0 r^{-\frac12}=0$ and
$\cL_0 r^{\frac32}=0$. A "generic" perturbation
$\tilde\cL=\cL_0+\tilde V$ with $\tilde V$ bounded, smooth, and
nicely decaying, will have zero energy solutions that behave just
like $r^{-\frac12}$ and $r^{\frac32}$, respectively. However, in
some cases $\tilde V$ is such that these two $\cL_0$ solutions will
be "in resonance" and produce a globally bounded zero energy
solution of~$\tilde\cL$ which behaves like $r^{\frac32}$ around zero
and $r^{-\frac12}$ around infinity
--- just like $\phi_0$.

For the parametrix construction in the following sections the
relevance of the zero energy resonance lies with the singularity of
the spectral measure of~$\cL$ at zero energy. Indeed, for $\cL_0$
the density of the spectral measure behaves like $\xi$ as $\xi\to0$,
whereas for $\cL$ we will show that it behaves like
$(\xi\log^2\xi)^{-1}$ as $\xi\to0$. We now briefly summarize the
results from~\cite{GZ} relevant for our purposes, see Section~3 in
their paper, in particular Example~3.10.

\begin{theorem}\label{thm:GZ}
a) For each $z \in \Compl$ there exists a fundamental system
$\tilphi(r,z)$, $\tiltheta(r,z)$ for $\cL-z$ which is analytic in
$z$ for each $r > 0$ and has the asymptotic behavior
\begin{equation}\label{eq:phitheta}
\tilphi(r,z) \sim r^\frac32, \qquad \tiltheta(r,z) \sim \frac12
r^{-\frac12} \text{\ \ as\ \ }r\to0
\end{equation}
In particular, their Wronskian is
$W(\theta(\cdot,z),\phi(\cdot,z))=1$ for all $z\in\Compl$. We remark
that $\phi(\cdot,z)$ is the Weyl-Titchmarsh solution\footnote{Our
$\phi(\cdot,z)$ is the $\tilde\phi(z,\cdot)$ function from~\cite{GZ}
where the analyticity is only required in a strip around~$\R$ -- but
here there is no need for this restriction.} of $\cL-z$ at $r=0$. By
convention, $\phi(\cdot,z), \theta(\cdot,z)$ are real-valued for
$z\in\R$.

b) For each  $ z\in\Compl$, $\Im z>0$, let $\psi^+(r,z)$ denote the
Weyl-Titchmarsh solution of $\cL-z$ at $r=\infty$ normalized so that
\[
\psi^+(r,z) \sim z^{-\frac14}\,e^{i z^{\frac12} r} \text{\ \ as\ \
}r\to \infty,\;\Im z^{\frac12}>0
\]
If $\xi>0$, then the limit $\psi^+(r,\xi+i0)$ exists point-wise for
all $r>0$ and we denote it by $\psi^+(r,\xi)$. Moreover, define
$\psi^-(\cdot,\xi):=\overline {\psi^+(\cdot,\xi)}$. Then
$\psi^+(r,\xi)$, $\psi^-(r,\xi)$ form a fundamental system of
$\cL-\xi$ with asymptotic behavior $\psi^\pm(r,\xi)\sim
\xi^{-\frac14}\,e^{\pm i\xi^{\frac12} r}$ as $r\to\infty$.

c) The spectral measure of $\cL$ is absolutely continuous and its
density is given by
\begin{equation} \label{sm}
\rho(\xi )  = \frac{1}{\pi} \Im\; \tilm(\xi +i0)\,\chi_{[\xi>0]}
\end{equation}
with the ``generalized Weyl-Titchmarsh" function
\begin{equation}\label{mw}
  \tilm(\xi) =\frac{W({\theta}(.,\xi),\,\psi^{+}(.,\xi))}{W(\psi^{+}(.,\xi),
    {\phi}(.,\xi))}
\end{equation}

d) The distorted Fourier transform defined as
  \begin{equation}\nonumber
    \calF: f\longrightarrow \hat{f}(\xi)=\lim_{b\rightarrow
      \infty}\int_{0}^{b}\tilphi(r,\xi)f(r)\,dr
  \end{equation}
  is a unitary operator from $L^2 (\R^+)$ to
  $L^{2}({{\R^+}},{\rho})$ and its inverse is given by
 \begin{equation}\nonumber
    \calF^{-1}: \hat{f}\longrightarrow f(r)=\lim_{\mu \rightarrow
      \infty}\int_{0}^{\mu} \tilphi(r,
    \xi)\hat{f}(\xi)\,{\rho}(\xi)\,d\xi
  \end{equation}
Here $\lim$ refers to the $L^{2}({{\R^+}},{\rho})$, respectively the
$L^2(\R^+)$, limit.
\end{theorem}

Needless to say, part b) above has nothing to do with~\cite{GZ} and
is standard. Most relevant for our computations are~\eqref{mw}
(which is  formula (3.22) in~\cite{GZ}), as well as the Fourier
inversion theorem in this context (see Theorem~3.5 in~\cite{GZ}).

Theorem~\ref{thm:GZ} of course also holds for $\cL_0$ instead
of~$\cL$.  In that case we have a Bessel equation with solutions
 \begin{align}
\phi(r;z) &=  2z^{-1/2} r^{1/2} J_1(z^{1/2} r) \label{eq:J1}\\
\theta(r;z) &= \frac{\pi}{4} z^{1/2} r^{1/2} [-Y_1(z^{1/2} r) +
\pi^{-1} \log(z)
J_1(z^{1/2} r) ] \nn \\
\psi(r;z) &= z^{1/2}r^{1/2} [-Y_1(z^{1/2} r) + iJ_1(z^{1/2} r)] = z^{1/2} r^{1/2} i H^{(1)}_1(z^{1/2}r)\nn\\
&= \theta(r;z) + m(z) \phi(r;z) \nn \\
m(z) &= \frac{\pi}{4} z [i-\pi^{-1}\log(z)], \quad
z\in\Compl\setminus\R^+ \nn
\end{align}
The last formula shows that for strongly singular potentials the
Weyl-Titchmarsh function ceases to be Herglotz, see~\cite{GZ} for
further discussion. Although we shall make no use of these formulas
for $\cL_0$, the reader should note the similarities between the
asymptotic expansions on $\tilphi$, $\tiltheta$ and $\psi^+$ we
derive below and the classical ones for the Bessel functions,
cf.~\cite{Wat}.

\subsection{Asymptotic behavior of $\tilphi$ and $\tiltheta$}

Beginning with two explicit solutions for $\cL f=0$, namely
\[
\phi_0(r) = \frac{r^\frac32}{1+r^2}, \qquad \theta_0(r) =
\frac{1-4r^2 \ln r - r^4}{2r^\frac12(1+r^2)} = r^{-\frac12}
(1-r^2)/2 - 2\phi_0(r)\log r
\]
we shall construct power series expansions for $\tilphi$ and
$\tiltheta$ from~\eqref{eq:phitheta} in $z \in\Compl$ when $r>0$ is
fixed.

\begin{proposition}  \label{pphitheta} For any $z\in\Compl$
the fundamental system $\tilphi(r,z)$, $\tiltheta(r,z)$ from
Theorem~\ref{thm:GZ} admits absolutely convergent asymptotic
expansions
\begin{align*}
\tilphi(r,z) &= \phi_0(r) + r^{-\frac12} \sum_{j =1}^\infty (r^2
z)^{j} \phi_j(r^2)\\
 \tiltheta(r,z) &=  r^{-\frac12}
\Big( 1-r^2 - \sum_{j =1}^\infty (r^2 z)^{j}  \theta_j(r^2) \Big)/2
- (2+z/4) \tilphi(r,z) \ln r
\end{align*}
where the functions $\phi_j$, $\theta_j$ are holomorphic in $U = \{
\Re u > -\frac12\} $ and satisfy the bounds
\begin{align*}
| \phi_j(u)| &\leq \frac{3C^j}{(j-1)!}  {\ln(1+|u|)}, \quad |\phi_1(u)| > \frac{1}{2}\log u \text{\ \ if\ \ }u\gg1 \\
 |\theta_1(u)|
&\leq C|u|,\qquad | \theta_j(u)| \leq \frac{C^j}{(j-1)!} {(1+|u|)},
\qquad u \in U
\end{align*}
Furthermore, \begin{equation}\label{eq:phi1} \phi_1(u) = \left\{
\begin{array}{ll} -\frac14\log u + \frac12 +
O(u^{-1}\log^2 u) & \text{\ \ as\ \ }u\to\infty \\
-\frac{u}{8} + \frac{u^2}{12} + O(u^3) & \text{\ \ as\ \ }u\to 0
\end{array}\right.
\end{equation}
\end{proposition}

\begin{proof}
We begin with $\tilphi$. We formally write
\[
\tilphi(r,z) = r^{-\frac12} \sum_{j=0}^\infty z^j  f_j(r)
\]
This becomes rigorous once we verify the convergence of the series
in any reasonable sense. The functions $f_j$ should solve
\[
\cL (r^{-\frac12} f_j) = r^{-\frac12} f_{j-1}, \qquad f_{-1} = 0
\]
The forward fundamental solution for $\cL$ is
\[
H(r,s) = \frac12(\phi_0(r) \theta_0(s) -\phi_0(s) \theta_0(r))1_{[r
> s]}
\]
Hence we have the iterative relation
\[
f_j(r) = \frac12\int_0^s  r^{\frac12} s^{-\frac12}  (\phi_0(r)
\theta_0(s) -\phi_0(s) \theta_0(r)) f_{j-1}(s)\, ds, \qquad f_0(r) =
\frac{r^2}{1+r^2}
\]
Using the expressions for $\phi_0$, $\theta_0$ we rewrite this as
\[
f_j(r) =  \int_0^r  \frac{r^2(-1+4s^2\log s+s^4) -s^2(-1+ 4r^2 \log
r
  + r^4)}{2s(1+r^2)(1+s^2)} f_{j-1}(s) ds
\]
We claim that all functions $f_j$ extend to even holomorphic
functions in any even simply connected domain not containing $\pm
i$, vanishing at $0$. Indeed, we now suppose that $f_{j-1}$ has
these properties and we shall prove them for $f_j$. Clearly, $f_j$
extends to a holomorphic function in any even simply connected
domain not containing $\pm i$ and $0$. We first show that at $0$
there is at most an isolated singularity.  For this we consider a
branch of the logarithm which is holomorphic in $\Compl \setminus
\R^+$ and show that $f_j(r+i0) = f_j(r-i0)$ for $r < 0$.
Disregarding the terms not involving logarithms, we need to show
that for any holomorphic function $g$ we have
\[
\int_{0}^{r+i0} (\log s - \log (r+i0)) g(s) ds =
\int_{0}^{r-i0} (\log s - \log (r-i0)) g(s) ds
\]
This is obvious since for $s < 0$ we have
\[
\log (s+i0) - \log (r+i0) = \log (s-i0) - \log (r-i0)
\]

The singularity at $0$ is a removable singularity. Indeed, for $s$
close to $0$ we have $|f_{j-1}(s)| \lesssim |s|$ which by a crude
bound on the denominator in the above integral leads to
$ |f_{j}(r)| \lesssim |r|$ (again with $r$ close to $0$). This also
shows that $f_j$ vanishes at $0$.

The fact that $f_j$ is even is obvious if we substitute $2 \log s$ and
$2 \log r$ by $\log s^2$ respectively $\log s^2$ in the integral. This
is allowed since due to the above discussion we can use any branch of
the logarithm. Indeed, denoting $\tilde f_{j-1}(s^2) = f_{j-1}(s)$
the  change of variable $s^2 = v$ yields
the iterative relation
\begin{equation}\label{eq:tilfj}
\tilde f_j(u) =  \int_0^{u}  \frac{u(-1+2 v \log v +v^2) -v(-1+ 2u  \log u
  + u^2)}{4v(1+u)(1+v)} \tilde f_{j-1}(v) dv, \qquad \tilde f_0 (u) = \frac{u}{1+u}
\end{equation}
Next, we obtain bounds on the functions $\tilde f_j$. To avoid the
singularity at $-1$ we restrict ourselves to the region $U = \{ \Re
u > -\frac12\}$.  We claim that the $\tilde f_j$ satisfy the bound
\[
|\tilde f_j(u)| \leq \frac{3C^{j}}{(j-1)!}  |u|^j \ln (1+|u|)
\]
The kernel above can be estimated by
\[
\left| \frac{u(-1+2 v \log v +v^2) -v(-1+ 2u \log u +
    u^2)}{2v(1+u)(1+v)}\right| \leq C \frac{|u|}{|v|}
\]
We have
\[
|\tilde f_0(u)| \leq 3 \frac{|u|}{1+|u|}
\]
which yields
\[
|\tilde f_1(u)| \leq 3C|u| \int_{0}^{|u|} \frac{1}{1+x} dx = 3C |u| \ln(1+|u|)
\]
From here on we use induction, noting that
\[
\int_{0}^{|u|} x^{j-1} \ln(1+x) dx \leq \frac{1}j |u|^j \ln(1+|u|)
\]
To conclude the proof, we note that the functions $\phi_j$ are given
by $\phi_j (u) = u^{-j} \tilde f_j(u)$ and satisfy the desired
pointwise bound. Finally, \eqref{eq:phi1} follows by an asymptotic
evaluation of the explicit integral~\eqref{eq:tilfj} with $j=1$,
which we leave to the reader.

\noindent The argument for the function $\theta$ is similar. The
ansatz
\begin{align*}
\tiltheta(r,z) &= r^{-\frac12} \Big( 1-r^2  - \sum_{j =1}^\infty z^j
g_j(r)\Big)/2 - (2+z/4)
  \phi(r,z) \ln r \\
   &= r^{-\frac12} \Big( 1-r^2 - \sum_{j =1}^\infty
z^j g_j(r)\Big)/2 + (2+z/4)
  \big(\phi_0(r) - \sum_{j=1}^\infty z^j r^{-\frac12} f_j(r)\big) \ln r
\end{align*}
yields a recurrence relation for the $g_j$ via $(\cL-z)\theta=0$.
Indeed, for $j =1$,
\begin{align*}
  \cL(r^{-\frac12} g_1(r)) &= \theta_0(r) - \cL\big(\frac12\phi_0(r) \log r +
  4r^{-\frac12} f_1(r)\log r\big)\\
  &= r^{-\frac12} \Big[r^2-  \frac{r^2(3+r^2)}{(1+r^2)^2} - \frac{8}{r^2} f_1(r) +
  \frac{8}{r} f_1'(r)\Big]
\end{align*}
where the important fact is that the quantity in brackets is even
analytic around $0$ and vanishes at $0$. A similar computation
yields for $j\ge2$
\[
\cL (r^{-\frac12} g_j(r)) = r^{-\frac12}\Big[ g_{j-1}(r) - r^{-2}
f_{j-1}(r)
  +  r^{-1}f'_{j-1}(r) - 8r^{-2} f_j(r) + 8r^{-1} f'_j(r)\Big]
\]
The same considerations as in the case of $f_j$ show that each $g_j$
is an even holomorphic function in any even simply connected domain
not containing $\pm i$. Also, the same bound for the fundamental
solution for $\cL$ leads to $|g_1(r)| \leq C r^4$ and more
generally, for $j\ge 2$,
\[
|g_j(r)|   \leq \frac{C^j}{(j-1)!} r^{2j}(1+r^2)
\]
The proof of the proposition is concluded.
\end{proof}

\begin{remark}\label{rem:log}
The logarithmic behavior of $\phi_1(u)$ for large $u$ is inherited
by $\phi(r,\xi)$; indeed, suppose that $1\gg\xi>0$ and
$r=\delta\xi^{-\frac12}$ where $\delta>0$ is small. Then the
proposition shows that
\[
\phi(r,\xi) \gtrsim r^{-\frac12} \log r
\]
The size of $\delta$ here only depends on various constants in the
expansion of $\phi$ and is thus itself an absolute constant. We
remark that the appearance of the $\log r$ term is a specific
feature of $\cL$ --- it does not occur for $\cL_0$,
see~\eqref{eq:J1} --- indicative of the fact that $\cL$ is a {\em
long range} perturbation of $\cL_0$. We shall see later that the
logarithm in $\phi$ produces crucial logarithmic factors in the
small $\xi$ asymptotics of the spectral density of~$\cL$, see
Proposition~\ref{atilphi} below.
\end{remark}

 We note that although the above series for $\tilphi$
converges for all $r,z$, we can only use it to obtain various
estimates for $\tilphi$ in the region $|z| r^2 \lesssim 1$. On the
other hand, in the region $\xi r^2 \gtrsim 1$ where $z=\xi>0$,  we
will represent $\tilphi$ in terms of $\psi^+$ and use the $\psi^+$
asymptotic expansion, described in what follows.

\subsection{The asymptotic behavior of $\psi^+$}

The following result provides good asymptotics for $\psi^+$ in the
region $r^2 \xi \gtrsim 1$.

\begin{proposition} For any $\xi>0$, the solution  $\psi^+(\cdot,\xi)$ from Theorem~\ref{thm:GZ}  is of
the form
\[
\psi^+(r,\xi) =  \xi^{-\frac14}e^{ir \xi^\frac12}
\sigma(r\xi^\frac12,r),\qquad r^2\xi\gtrsim 1
\]
where $\sigma$ admits the asymptotic series approximation
\[
\sigma(q,r) \approx \sum_{j=0}^\infty q^{-j} \psi^+_j(r), \qquad
\psi^+_0 = 1, \qquad \psi_1^+ = \frac{3i}{8} + O(\frac{1}{1+r^2})
\]
with zero order symbols  $\psi^+_j(r)$ that are analytic at
infinity,
\[\sup_{r>0}
|(r \partial_r)^k \psi^+_j(r)| <\infty
\]
in the sense that for all large integers $j_0$, and all indices
$\alpha$, $\beta$, we have
\[
\sup_{r>0}\Bigl|(r \partial_r)^\alpha (q \partial_q)^\beta
\Big[\sigma(q,r)
  - \sum_{j=0}^{j_0} q^{-j} \psi_j^+(r)\Big]\Bigr| \leq
c_{\alpha,\beta,j_0}  q^{-j_0-1}
\]
for all $q>1$.
 \label{ppsipsi}\end{proposition}

\begin{proof}
With the notation
\[
\sigma(q,r) =  \xi^{\frac14}\psi^+(r,\xi)  e^{-ir \xi^\frac12}
\]
we need to solve the conjugated equation
\begin{equation}
\left(-\partial_r^2 - 2 i \xi^\frac12 \partial_r + \frac{3}{4r^2} -
  \frac{8}{(1+r^2)^2}\right) \sigma(r\xi^{\frac12}, r) = 0
\label{conjug}\end{equation}
We look for a formal power series solving this equation,
\begin{equation}
\sum_{j=0}^\infty \xi^{-\frac{j}2} f_j(r)
\label{formal} \end{equation}
This yields a recurrence relation for the $f_j$'s,
\[
2i \partial_r f_j = \left(-\partial_r^2 + \frac{3}{4r^2} -
  \frac{8}{(1+r^2)^2}\right) f_{j-1}, \qquad f_0 = 1
\]
which is solved by
\[
f_j = \frac{i}{2} \partial_r f_{j-1} + \frac{i}{2} \int_{r}^\infty
\left( \frac{3}{4s^2} -
  \frac{8}{(1+s^2)^2}\right) f_{j-1}(s)\, ds
\]
Extending this into the complex domain, it is easy to see that the
functions $f_j$ are holomorphic in $\Compl \setminus [-i,i]$. They
are also holomorphic at $\infty$, and the leading term in the Taylor
series at $\infty$ is $r^{-j}$. At $0$, on the other hand, $f_j$ are
singular. The worst singularity is of power type, namely $r^{-j}$;
however, weaker terms contain logarithms and powers of logarithms so
it is not easy to obtain a complete expansion. Instead we contend
ourselves with a weaker estimate, namely
\[
|(r \partial_r)^k f_j| \leq c_{jk}\, r^{-j} \qquad \forall r>0
\]
which is easy to establish inductively. The functions
\[
\psi_j^+(r) := r^j f_j(r)
\]
now satisfy the desired bounds due to the bounds above on $f_j$.

Unlike in the expansion for small $r$, here we make no effort to
obtain a uniform estimate on the size of the derivatives of
$\psi_j^+$.  This is because we do not expect the formal series
\eqref{formal} to converge, on account of the fact that derivatives
are lost in the iterative construction of the $f_j$'s. Instead we
can construct an approximate sum, i.e., a function
$\sigma_{ap}(q,r)$ with the property that for each $j_0 \geq 0$ we
have
\begin{equation}\label{eq:sigma_schr}
\Big|(r \partial_r)^\alpha (q \partial_q)^\beta
\big[\sigma_{ap}(q,r)
  - \sum_{j=0}^{j_0} q^{-j} \psi_j^+(r)\big]\Big| \leq
c_{\alpha,\beta,j_0}\,  q^{-j_0-1}
\end{equation}
The construction of $\sigma_{ap}(q,r)$ is standard in symbol
calculus; indeed, we can set
\[
\sigma_{ap}(q,r) := \sum_{j=0}^\infty q^{-j} \psi_j^+(r)
\chi(q\delta_j)
\]
where $\delta_j\to0$ sufficiently fast and $\chi$ is a cut-off
function which vanishes around zero and is equal to one for large
arguments. The bound~\eqref{eq:sigma_schr} implies that
$\sigma_{ap}(r\xi^\frac12,r)$ is a good approximate solution for
\eqref{conjug} at infinity, namely the error
\[
e(r \xi^\frac12,r) = \left(-\partial_r^2 -2 i \xi^\frac12 \partial_r
+
  \frac{3}{4r^2} - \frac{8}{(1+r^2)^2}\right) \sigma_{ap}(r,\xi)
\]
satisfies for all indices $\alpha,\beta,j$
\[
| (r \partial_r)^\alpha (q \partial_q)^\beta e(q,r) | \leq
c_{\alpha,\beta,j}\,  r^{-2} q^{-j}
\]
To conclude the proof it remains to solve the equation for the
difference $\sigma_1 = -\sigma + \sigma_{ap}$,
\[
 \left(-\partial_r^2 - 2i \xi^\frac12 \partial_r +
  \frac{3}{4r^2} - \frac{8}{(1+r^2)^2}\right) \sigma_1(r
\xi^\frac12,r) = e(r \xi^\frac12,r)
\]
with zero Cauchy data at infinity. We claim that the solution
$\sigma_1$ satisfies
\[
| (r \partial_r)^\alpha (q \partial_q)^\beta \sigma_1(q,r) | \leq
c_{\alpha,\beta,j}\,  q^{-j}, \qquad j \geq 2
\]
Note that this finishes the proof by defining $\sigma=\sigma_{ap}
-\sigma_1$.  A change of variable allows us to switch from the pair
of operators $(r\partial_r, q\partial_q)$ to $(r\partial_r, \xi
\partial_\xi)$ with comparable bounds. We rewrite the above equation as
a first order system for $(v_1,v_2) = (\sigma_1,r \partial_r \sigma_1)$:
\[
\partial_r \left( \begin{array}{c} v_1 \cr v_2 \end{array}\right) -
\left( \begin{array}{cc} 0 & r^{-1} \cr
 \frac{3}{4r} - \frac{8r}{(1+r^2)^2}
 &r^{-1}  -2 i\xi^\frac12
 \end{array}\right)
 \left( \begin{array}{c} v_1 \cr v_2 \end{array}\right) =
 \left( \begin{array}{c} 0 \cr re \end{array}\right)
\]
Then we have
\[
\frac{d}{dr} |v|^2 \gtrsim  -r^{-1} |v|^2 - r |v||e|
\]
which gives
\[
\frac{d}{dr} |v| \geq -C(r^{-1}|v| + r |e|)
\]
and by Gronwall
\[
|v(r)| \leq \int_{r}^\infty \left(\frac{s}{r}\right)^C s |e(s)| ds
\]
Then for large $j$ we have
\begin{equation}
|e| \lesssim \xi^{-\frac{j}{2}} r^{-j-2} \implies |v| \lesssim
\xi^{-\frac{j}{2}} r^{-j} = q^{-j}
 \label{ev}\end{equation}
 To estimate derivatives of $v$ we commute them with the operator. For
derivatives with respect to $r$ we have
\[
\partial_r (r \partial_r) \left( \begin{array}{c} v_1 \cr v_2 \end{array}\right) -
\left( \begin{array}{cc} 0 & \frac{1}r \cr
 \frac{3}{4r} - \frac{8r}{(1+r^2)^2}
 &\frac1r  -2 i\xi^\frac12
 \end{array}\right)
 (r \partial_r) \left( \begin{array}{c} v_1 \cr v_2 \end{array}\right)
=
\left( \begin{array}{cc} 0 & \frac1r \cr
 \frac{3}{4r} - \frac{8r(3r^2-1)}{(1+r^2)^3}
 &\frac1r
 \end{array}\right)
 \left( \begin{array}{c} v_1 \cr v_2 \end{array}\right) +
 \left( \begin{array}{c} 0 \cr r \partial_r (re) \end{array}\right)
\]
But the right-hand side is bounded by $r^{-j-1}$ from the previous
step and the hypothesis on $e$, therefore as above $r \partial_r v$
is bounded by $r^{-j}$.

\noindent We argue similarly for the $\xi$ derivatives. We have
\[
\partial_r (\xi \partial_\xi) \left( \begin{array}{c} v_1 \cr v_2 \end{array}\right) -
\left( \begin{array}{cc} 0 & \frac{1}r \cr
 \frac{3}{4r} - \frac{8r}{(1+r^2)^2}
 &\frac1r  -2 i\xi^\frac12
 \end{array}\right)
 (\xi \partial_\xi) \left( \begin{array}{c} v_1 \cr v_2 \end{array}\right)
=
\left( \begin{array}{cc} 0 & 0 \cr
 0 &   i\xi^\frac12
 \end{array}\right)
 \left( \begin{array}{c} v_1 \cr v_2 \end{array}\right) +
 \left( \begin{array}{c} 0 \cr \xi \partial_\xi (re) \end{array}\right)
\]
The only difference is in the first term on the right, for which we
write $\xi^\frac12 = r^{-1} q$ and we use the decay property of $v$
with $j$ replaced by $j+1$:
\[\begin{split}
|\xi^{\frac12} v_2| \less \xi^{\frac12} q^{-j-1}\less r^{-1}
q^{-j},\qquad  |\xi \partial_\xi (re)| \less r^{-1} q^{-j}
\end{split}
\]
as desired.  Finally, higher order derivatives are estimated by
induction using the above arguments at each step.
\end{proof}

\subsection{Structure of the spectral measure of $\cL$}

We begin by relating the functions $\tilphi$, $\tiltheta$ and
$\psi^\pm$. By examining the asymptotics at $r = 0$ we see that
\begin{equation}
W(\tiltheta,\tilphi) = 1\label{phitheta}\end{equation} Also by
examining the asymptotics as $r \to \infty$ we obtain
\begin{equation}
W(\psi^+,\psi^-) = -2 i \label{psipsi}\end{equation} Hence we can
express the $\cL -\xi$ solutions in either the  $\tilphi$,
$\tiltheta$ basis or the $\psi^\pm$  basis. On the other hand,
$\tilphi$, $\tiltheta$ are real-valued while the real and imaginary
parts of $\psi^\pm$ are equally strong. Hence the two bases are
quite separated. These are the main ingredients of the next result.

\begin{proposition} \label{atilphi}
a) We have
\begin{equation}\label{eq:phiapsi}
\tilphi(r,\xi) =  a(\xi) \psi^+(r,\xi) + \overline{ a(\xi) \psi^+(r,\xi)}
\end{equation}
where $a$ is smooth, always nonzero, and has size\footnote{$a\asymp
b$ means that for some constant $C$ one has $C^{-1} a<b<Ca$}
\[
|a(\xi)| \asymp \left\{ \begin{array}{cc} -
   \xi^\frac12  \log \xi & \xi \ll 1 \cr  \cr \xi^{-\frac12} &  \xi \gtrsim 1
\end{array}\right.
\]
Moreover, it satisfies the symbol type bounds
\[
| (\xi \partial_\xi)^k a(\xi) | \leq c_k |a(\xi)|\quad
\forall\;\xi>0
\]

b) The spectral measure $\rho(\xi)d\xi$ has density
\[
\rho(\xi) = \frac{1}\pi |a(\xi)|^{-2}
\]
and therefore satisfies
\[
\rho(\xi)\asymp \left\{ \begin{array}{cc}
   \frac{1}{\xi  (\log \xi)^2} & \xi \ll 1 \cr  \cr \xi &  \xi \gtrsim 1
\end{array}\right.
\]
\end{proposition}

\begin{proof}
a) Since $\tilphi$ is real-valued, due to \eqref{psipsi}, the
relation~\eqref{eq:phiapsi} above holds with
\[
a(\xi) = -\frac{i}{2} W(\tilphi(\cdot,\xi),{\psi^-}(\cdot,\xi))
\]
We evaluate the Wronskian in the region where both the
$\psi^+(r,\xi)$ and $\tilphi(r,\xi)$ asymptotics are useful, i.e.,
where $r^2 \xi \approx 1$. By Proposition~\ref{pphitheta} we obtain
that both $ \tilphi(\xi^{-\frac12},\xi)$ and $(r
\partial_r\tilphi)(\xi^{-\frac12},\xi)$ can be expressed in the form
$\xi^{\frac14} f(\xi^{-1})$ with $f$ holomorphic and satisfying
\[
|f(u)| \lesssim \log(1+|u|), \qquad \Re u > \frac14
\]
On the other hand, it follows from Proposition~\ref{ppsipsi}  that
both $ \psi^+(\xi^{-\frac12},\xi)$ and $(r
\partial_r\psi^+)(\xi^{-\frac12},\xi)$ can be expressed in the form
$\xi^{-\frac14} h(\xi^{-\frac12})$ with $h$ satisfying symbol type
bounds
\[
|(r \partial_r)^k h(r)| \leq c_k
\]
Combining the two expressions above, it follows that $a$ is a sum of
terms of the form $\xi^\frac12 f(\xi^{-1}) h(\xi^{-\frac12})$ with $f,h$ as
above. The bounds from above on $a$ and its derivatives follow.

It remains to prove the bound from below on $a$, which is more
delicate. By \eqref{psipsi} we have
\[
\Im (\psi^{+}(r,\theta) \partial_r \psi^-(r,\theta))  = -1
\]
Since $\tilphi$ is real-valued,  this gives
\[
\Im\big[ \partial_r \psi^+(r,\xi)
W(\tilphi(\cdot,\xi),{\psi^-}(\cdot,\xi))\big] =
-\partial_r\tilphi(r,\xi)
\]
which implies that for all $r$ we have
\[
|a(\xi)| \geq \frac{| \partial_r\tilphi(r,\xi)|}{2|\partial_r
\psi^+(r,\xi)|}
\]
We use this relation for $r = \delta \xi^{-\frac12}$ with a small constant
$\delta$. Then by Proposition~\ref{pphitheta} we have
\[
|\partial_r \tilphi(r,\xi)| \gtrsim r^{-\frac32} \ln(1+r^2)
\]
while by  Proposition~\ref{ppsipsi}
\[
{|\partial_r \psi^+(r,\xi)|} \lesssim \xi^{\frac14} (r^2 \xi)^{-j_0}
\]
This give the desired bound from below on $a$.

b) By \eqref{phitheta} we can express $\psi^+$ in terms of $\theta$ and
$\phi$ by
\[
\psi^+ = -\tilphi W(\psi^+,\tiltheta) + \tiltheta W( \psi^+,\tilphi)
\]
Since both $ \tilphi $ and $\tiltheta$ are real-valued, by inserting
this into \eqref{psipsi} we obtain the relation
\[
\Im ( W(\psi^+,\tiltheta)   {W( \psi^-,\tilphi)}) = -1
\]
Inserting this in the expression for the spectral measure \eqref{sm}
and taking \eqref{mw} into account we obtain
\[
\rho (\xi) = \frac1\pi \frac{\Im ( W(\psi^+,\tiltheta)
   {W( \psi^-,\tilphi)})}{| W(\psi^+,\tilphi)|^2}
= \frac1\pi |W(\psi^+,\tilphi)|^{-2} = \frac{1}{\pi |a(\xi)|^2}
\]
as desired.
\end{proof}

 \section{The transference identity}
 \label{sec:transference}

 Returning to the radiation part $\tilde{\epsilon}$ in
 \eqref{equation2}, the idea is to expand it in terms of the
 generalized Fourier basis\footnote{We now return to the
variable $R$ as the independent spatial variable instead of $r$ as
in the previous section.} $\tilphi(R, \xi)$ associated with the
 operator $\calL=-\partial_{R}^{2}+\frac{3}{4
   R^{2}}-\frac{8}{(1+R^{2})^{2}}$, i.e., write
 \begin{equation}\nonumber
   \tileps(\tau, R)=\int_{0}^{\infty}  x(\tau,\xi)
   \tilphi(R,\xi)\rho(\xi)\,d\xi
 \end{equation}
 and deduce a transport equation for the Fourier coefficients $x(\tau,
 \xi)$.  The main difficulty in doing this is caused by the operator
 $R \partial_R$ which is not diagonal in the Fourier basis.  Our
 strategy for dealing with this is to replace it with $2 \xi
 \partial_\xi$ modulo an error which we treat perturbatively.  The
 operator $R\partial_{R}-2\xi\partial_{\xi}$ is natural since it
 annihilates the expression $e^{i\xi^\frac12 R}$ arising in the
 asymptotic expansion of $\tilphi(R, \xi)$ for large $R$.
 Consequently, we define the error operator $\calK$ by
 \begin{equation}\label{eq:transfer}
 \widehat{R \partial_R u} = - 2 \xi \partial_\xi \hat u + \calK \hat u
 \end{equation}
 where $\hat{f}=\calF f$ is the ``distorted Fourier transform" from
 Theorem~\ref{thm:GZ}.
 Using the expressions for the direct and inverse Fourier transform in
 that theorem we obtain
 \[
 \calK f(\eta)= \Bla
 \int_{0}^{\infty}f(\xi) R\partial_{R} \tilphi(R,
 \xi)\rho(\xi)\,d\xi\,,\, \tilphi(R, \eta)\Bra_{L^2_R} + \Bla
 \int_{0}^{\infty}2\xi\partial_{\xi} f(\xi) \tilphi(R,
 \xi)\rho(\xi)\,d\xi\,,\, \tilphi(R, \eta)\Bra_{L^2_R}
 \]
 Integrating by parts with respect to $\xi$ in the second expression
 we obtain
 \begin{equation}\label{calk}
   \calK f(\eta)= \Bla
\int_{0}^{\infty}f(\xi)[R\partial_{R}-2\xi\partial_{\xi}]\tilphi(R,\xi)\rho(\xi)\,d\xi\,,\,
\tilphi(R, \eta)\Bra_{L^2_R} - 2 \left(1+ \frac{\eta
  \rho'(\eta)}{\rho(\eta)}\right) f(\eta)
\end{equation}
 where the scalar product is to be interpreted in the principal value
 sense with $f\in C_0^\infty((0,\infty))$.
 Apriori we have
 \[
 \calK: C_0^\infty((0,\infty)) \to C^\infty((0,\infty))
 \]
 therefore we can write
 \begin{equation}\nonumber
   \calK f(\eta)=\int_{0}^{\infty}K(\eta, \xi)f(\xi)\,d\xi
 \end{equation}
 for a distribution valued
   function $\eta\rightarrow K(\eta, \xi)$.
 We refer to \eqref{eq:transfer} as the {\it{transference identity}} to indicate that we
 are  transferring derivatives from $R$ to~$\xi$.  To asses
 its usefulness we need to understand the boundedness properties of the operator
 $\calK$. We begin with a description of the kernel $K(\eta,\xi)$.

 \begin{theorem}\label{tp}
   The operator $\calK$ can be written as
   \begin{equation}\label{eq:kern}
   \calK = -\Big(\frac32 +\frac{\eta \rho'(\eta)}{\rho(\eta)}\Big)\delta(\xi-\eta) + \calK_0
   \end{equation}
   where the operator $\calK_0$ has a kernel $K_0(\eta, \xi)$ of the
   form\footnote{The kernel below is interpreted in the principal
     value sense}
   \begin{equation} \label{ketaxi} K_0(\eta,
     \xi)=\frac{\rho(\xi)}{\xi-\eta} F(\xi,\eta)
   \end{equation}
   with a symmetric function $F(\eta, \xi)$ of class $C^2$ in
   $(0,\infty) \times (0,\infty)$ satisfying the bounds
   \[\begin{split}
   | F(\xi,\eta)| &\lesssim \left\{ \begin{array}{cc} \xi+\eta &
       \xi+\eta \leq 1 \cr (\xi+\eta)^{-\frac32} (1+|\xi^\frac12
       -\eta^\frac12|)^{-N} & \xi+\eta \geq 1
     \end{array} \right.\\
   | \partial_{\xi} F(\xi,\eta)|+| \partial_{\eta} F(\xi,\eta)| &\lesssim \left\{
     \begin{array}{cc} 1 & \xi+\eta \leq 1 \cr (\xi+\eta)^{-2}
       (1+|\xi^\frac12 -\eta^\frac12|)^{-N} & \xi+\eta \geq 1
     \end{array} \right.\\
  \sup_{j+k=2} | \partial^j_{\xi}\partial^k_{\eta} F(\xi,\eta)| &\lesssim \left\{
     \begin{array}{cc} |\log(\xi+\eta)|^3 & \xi+\eta \leq 1 \cr
       (\xi+\eta)^{-\frac52} (1+|\xi^\frac12 -\eta^\frac12|)^{-N} &
       \xi+\eta \geq 1
     \end{array} \right.
     \end{split}
   \]
   where $N$ an arbitrary large integer.
 \end{theorem}

 \begin{proof}
   We first establish the off-diagonal behavior of $K$, and later return to the issue
   of identifying the $\delta$-measure that sits on the diagonal. We begin with \eqref{calk}
   with $f \in C_0^\infty((0,\infty))$. The integral
   \[
   u(R) =
   \int_{0}^{\infty}f(\xi)[R\partial_{R}-2\xi\partial_{\xi}]\tilphi(R,
   \xi)\rho(\xi)\,d\xi
   \]
   behaves like $R^{\frac32}$ at $0$ and is a Schwartz function at
   infinity. The second factor $\tilphi(R,\eta)$ in \eqref{calk} also
   decays like $R^{\frac32}$ at $0$ but at infinity it is only bounded
   with bounded derivatives. Then the following integration by parts
   is justified:
   \[
   \eta \calK f(\eta) = \Bla u, \cL \tilphi(R, \eta)\Bra_{L^2_R} = \Bla
   \cL u, \tilphi(R, \eta)\Bra_{L^2_R}
   \]
   Moreover,
   \[
   \begin{split}
     \cL u =& \int_{0}^{\infty}f(\xi)[\cL,R\partial_{R}] \tilphi(R,
     \xi)\rho(\xi)\,d\xi +
     \int_{0}^{\infty}f(\xi)(R\partial_{R}-2\xi\partial_{\xi}) \xi
     \tilphi(R,\xi)\rho(\xi)\,d\xi \\ = & \int_{0}^{\infty}f(\xi) [\cL,R\partial_{R}]
     \tilphi(R, \xi)\rho(\xi)\,d\xi + \int_{0}^{\infty}\xi
     f(\xi)(R\partial_{R}-2\xi\partial_{\xi})
     \tilphi(R,\xi)\rho(\xi)\,d\xi - 2\int_{0}^{\infty}\xi
     f(\xi)\tilphi(R,\xi)\rho(\xi)\,d\xi
   \end{split}
   \]
   with the commutator
   \[
   [\cL,R\partial_{R}] = 2\cL + \frac{16}{(1+R^2)^2} -
   \frac{32R^2}{(1+R^2)^3} =: 2\cL + W(R)
   \]
    Thus,
    \[
\cL u =\int_{0}^{\infty}f(\xi) W(R)
     \tilphi(R, \xi)\rho(\xi)\,d\xi + \int_{0}^{\infty}\xi
     f(\xi)(R\partial_{R}-2\xi\partial_{\xi})
     \tilphi(R,\xi)\rho(\xi)\,d\xi
    \]
   Hence we obtain
   \[
   \eta \calK f(\eta) - \calK (\xi f)(\eta) = \Bla
   \int_{0}^{\infty}f(\xi) W(R) \tilphi(R, \xi)\rho(\xi)\,d\xi ,
   \tilphi(R, \eta)\Bra_{L^2_R}
   \]
   The double integral on the right-hand side is absolutely
   convergent, therefore we can change the order of integration to
   obtain
   \[
   (\eta -\xi) K(\eta,\xi) = \rho(\xi) \Bla W(R) \tilphi(R, \xi),
   \tilphi(R, \eta)\Bra_{L^2_R}
   \]
   This leads to the representation in \eqref{ketaxi} when $\xi\ne\eta$ with
   \[
   F(\xi,\eta) = \Bla W(R) \tilphi(R, \xi), \tilphi(R,
   \eta)\Bra_{L^2_R}
   \]
   It remains to study its size and regularity. First, due
   to our pointwise bound from the previous section,
   \begin{equation}\begin{split}
   \sup_{R\ge0}|\tilphi(R, \xi)|
   &\less \la\xi\ra^{-\frac34},\\
    |R\partial_R \phi(R,\xi)|& \less \min(R\xi^{-\frac14},R^{\frac32}) \qquad \forall\; \xi>1 \\
   |\partial_\xi \tilphi(R, \xi)|
   &\less \min(R\xi^{-\frac54}, R^{\frac72}) \qquad \forall\;
   \xi>1/2
   \\
|\partial_\xi \tilphi(R, \xi)|
   &\less \min(R^{\frac32}\log(1+R^2),\xi^{-\frac14}|\log\xi|R) \qquad \forall\;
   0<\xi<1/2\\
|\partial^2_\xi \tilphi(R, \xi)|
   &\less \min(R^2\xi^{-\frac74}, R^{\frac{11}{2}}) \qquad \forall\;
   \xi>1/2
   \\
|\partial^2_\xi \tilphi(R, \xi)|
   &\less \min(R^{\frac72}\log(1+R^2),\xi^{-\frac34}|\log\xi|\,R^2) \qquad \forall\;
   0<\xi<1/2
   \end{split}\label{eq:phi_est}
   \end{equation}
   we always have the estimates
   \begin{equation}\label{eq:Fbd}\begin{split}
   |F(\xi,\eta)| &\less \la\xi\ra^{-\frac34}\la\eta\ra^{-\frac34},\\
 |\partial_\xi F(\xi,\eta)| &\less
 \la\xi\ra^{-\frac54}\la\eta\ra^{-\frac34}, \quad
 |\partial_\eta F(\xi,\eta)| \less
 \la\xi\ra^{-\frac34}\la\eta\ra^{-\frac54}, \\
 |\partial^2_{\xi\eta} F(\xi,\eta)| &\less
\la\xi\ra^{-\frac54}\la\eta\ra^{-\frac54} \qquad \forall\;
\xi+\eta\gtrsim 1\\
|\partial_\xi^2 F(\xi,\eta)| &\less \xi^{-\frac74}\eta^{-\frac34} \qquad \forall\;\xi>1,\,\eta>1\\
|\partial_\eta^2 F(\xi,\eta)| &\less \xi^{-\frac34}\eta^{-\frac74}
\qquad \forall\;\xi>1,\,\eta>1
\end{split}
\end{equation}
   They are only useful when $\xi$ and $\eta$ are very close. To improve on them, we consider two
   cases:

   {\bf Case 1: $1 \lesssim \xi+\eta $.}
   To capture the cancelations when $\xi$ and $\eta$ are separated we
   resort to another integration by parts,
   \[
   \eta F(\xi,\eta) = \Bla W(R) \tilphi(R, \xi), \cL \tilphi(R,
   \eta)\Bra = \Bla [\cL,W(R)] \tilphi(R, \xi), \tilphi(R, \eta)\Bra +
   \xi F(\xi,\eta)
   \]
   Hence, evaluating the commutator,
   \begin{equation} \label{Kk1} (\eta -\xi) F(\xi,\eta) = -\Bla (2 W_R
     \partial_R + W_{RR}) \tilphi(R, \xi), \tilphi(R, \eta)\Bra
   \end{equation}
   Since $W_R(0)=0$ it follows that $ (2 W_R \partial_R + W_{RR})
   \tilphi(R, \xi)$ has the same behavior as $\tilphi(R, \xi)$ in the
   first region. Then we can repeat the argument above to obtain
   \[
   (\eta -\xi)^2 F(\xi,\eta) = -\Bla [\cL,2 W_R \partial_R + W_{RR}]
   \tilphi(R, \xi), \tilphi(R, \eta)\Bra
   \]
   The second commutator has the form, with $V(R):= -
   8(1+R^2)^{-2}$,
   \[ [\cL,2 W_R \partial_R + W_{RR}] = 4 W_{RR} \cL - 4W_{RRR}
   \partial_R -W_{RRRR}
   + 3R^{-2} (R^{-1} W_R -W_{RR}) - 2 W_R V_R -4W_{RR} V
   \]
   Since $R^{-1} W_R(R) -W_{RR}(R)=O(R^2)$ this leads to
   \[
   (\eta -\xi)^2 F(\xi,\eta) = \Bla ( W^{odd}(R) \partial_R + W^{even}(R) +
   \xi W^{even}(R) ) \tilphi(R, \xi), \tilphi(R, \eta)\Bra
   \]
   where by $ W^{odd}$, respectively $W^{even}$, we have generically
   denoted odd, respectively even, nonsingular rational functions with
   good decay at infinity.  Inductively, one now verifies the identity
        \begin{equation} \begin{split}
     (\eta -\xi)^{2k} F(\xi,\eta) &=  \Bla \Big( \sum_{j=0}^{k-1} \xi^{j}\, W_{kj}^{odd}(R)\, \partial_R +
     \sum_{\ell=0}^k \xi^\ell W_{k\ell}^{even}(R) \Big) \tilphi(R, \xi), \tilphi(R, \eta)\Bra
     \label{Kk} \\
\la R\ra|W_{kj}^{odd}(R)| + |W_{k\ell}^{even}(R)| &\less \la
R\ra^{-4-2k} \qquad \forall\; j,\ell
     \end{split}\end{equation}
   By means of the pointwise bounds on $\tilphi$ and $\partial_R
   \tilphi$ from~\eqref{eq:phi_est} we infer from this that
   \[
   |F(\xi,\eta)| \lesssim \frac{\xi^{k -\frac34}
     \la\eta\ra^{-\frac34}} {(\eta -\xi)^{2k}} \qquad\forall \;
     \xi\gtrsim 1, \;\eta>0
   \]
Combining this estimate with \eqref{eq:Fbd} yields, for arbitrary
$N$,
\[ |F(\xi,\eta)|\less (\xi+\eta)^{-\frac32} (1+|\xi^\frac12
       -\eta^\frac12|)^{-N} \text{\ \ provided\ \ } \xi+\eta \gtrsim 1,
       \]
as claimed.
   For the derivatives of $F$ we follow a similar procedure. If $\xi$
   and $\eta$ are comparable,  then from~\eqref{eq:Fbd},
   $|\partial_\eta F(\xi,\eta)|\less \la \xi\ra^{-2}$.
   Otherwise we differentiate with respect to $\eta$ in \eqref{Kk}.
   This yields
   \[
   (\eta -\xi)^{2k} \partial_\eta F(\xi,\eta) = \Bla \Big( \sum_{j=0}^{k-1} \xi^{j}\, W_{kj}^{odd}(R)\, \partial_R +
     \sum_{\ell=0}^k \xi^\ell W_{k\ell}^{even}(R) \Big)\tilphi(R, \xi),
   \partial_\eta \tilphi(R, \eta)\Bra - 2k (\eta -\xi)^{2k-1}
   F(\xi,\eta)
   \]
   Using also the bound on $F$ from above we obtain
   \[
   |\partial_\eta F(\xi,\eta)| \lesssim \frac{\xi^{k-\frac34}\eta^{-\frac54}} {(\eta -\xi)^{2k}}, \qquad 1
   \lesssim \xi,\eta
   \]
   respectively
   \[
   |\partial_\eta F(\xi,\eta)| \lesssim\frac{
     \eta^{-\frac54}} {(\eta -\xi)^{2k}} \qquad \xi \ll 1 \lesssim
   \eta
   \]
   and
   \[
   |\partial_\eta F(\xi,\eta)| \lesssim\frac{\xi^{k-\frac34}}{(\eta -\xi)^{2k}} \qquad \eta \ll 1 \lesssim
   \xi
   \]
   which again yield the desired bounds.
   Finally, we consider the second order derivatives with respect to
   $\xi$ and $\eta$. For $\xi$ and $\eta$ close we again use the
   bound from~\eqref{eq:Fbd}.
   Otherwise we differentiate twice in \eqref{Kk} and continue as
   before.  We note that it is important here that the decay of
   $W_{kj}^{odd}$ and $W_{k\ell}^{even}$ improves with $k$. This is because the
   second order derivative bound at $0$ has a sizeable growth at
   infinity which has to be canceled,
   \[
   | \partial_\xi^2 \tilphi(R, 0)| \approx R^\frac72 \ln R
   \]

   {\bf Case 2: $ \xi,\eta \ll 1$}. Our first observation is that $F(0,0) =
   0$. This is easy to verify by direct integration, and is
   heuristically justified by the fact that $W = [L,R\partial_R]$. The
   pointwise bound
   \[
   | \partial_\xi F(\xi,\eta)| \lesssim 1
   \]
   follows by a direct computation. The second order derivative bound
   is, however, more delicate. We have at our disposal the pointwise
   bounds
   \[
   |\partial^j_\xi \tilphi(R, \xi)| \lesssim \left\{ \begin{array}{cc}
       R^{-\frac12+2j} \ln (1+R^2) & R < \xi^{-\frac12} \\ \xi^{\frac14-j/2}
       |\log \xi|\, R^j & R \geq \xi^{-\frac12}
     \end{array} \right., \qquad j = 0,1,2
   \]
   If $\eta < \xi<1/2$, then these bounds imply that
   \[\begin{split}
   | \partial^2_{\xi\eta} F(\xi,\eta)| &\lesssim
   \int_{0}^{\xi^{-\frac12}} \la R\ra^{-4} R^3 (\log(1+R^2))^2\, dR +
   \int_{\xi^{-\frac12}}^{\eta^{-\frac12}} \la R\ra^{-4} R^{\frac52} \xi^{-\frac14}|\log\xi|\log(1+R^2)
   \,dR\\
   &\quad
   +\int_{\eta^{-\frac12}}^\infty \la R\ra^{-2} \xi^{-\frac14}
   \eta^{-\frac14} |\log\xi||\log\eta|\, dR
   \end{split}
   \]
   The main contribution comes from the first term. Integrating we
   obtain
   \[
   | \partial^2_{\xi\eta} F(\xi,\eta)| \lesssim |\log \xi|^3
   \]
   A similar computation yields, again when $\eta<\xi<1/2$,
\[\begin{split}
   | \partial^2_{\xi} F(\xi,\eta)| &\lesssim
   \int_{0}^{\xi^{-\frac12}} \la R\ra^{-4} R^3 (\log(1+R^2))^2\, dR +
   \int_{\xi^{-\frac12}}^{\eta^{-\frac12}} \la R\ra^{-4} R^{\frac32} \xi^{-\frac34}|\log\xi|\log(1+R^2)
   \,dR\\
   &\quad
   +\int_{\eta^{-\frac12}}^\infty \la R\ra^{-2} \xi^{-\frac34}
   \eta^{\frac14} |\log\xi||\log\eta|\, dR \less |\log \xi|^3
   \end{split}
   \]
   It remains to consider the expression $\partial_\xi^2 F(\xi,\eta)$
   for $ \xi \ll \eta<1/2$.  Differentiating in \eqref{Kk1} we
   obtain
   \[
   (\eta-\xi) \partial_\xi^2 F(\xi,\eta) = 2 \partial_\xi
   F(\xi,\eta) - \Bla \partial_\xi^2 \tilphi(R, \xi), (2 W_R \partial_R + W_{RR})
    \tilphi(R, \eta)\Bra
   \]
   We differentiate and integrate with respect to $\eta$ to obtain
   \begin{equation}\label{eq:letzt}
   (\eta-\xi) \partial_\xi^2 F(\xi,\eta) = \int_\xi^\eta \Big[2
   \partial^2_{\xi\zeta} F(\xi,\zeta) - \Bla  \partial^2_\xi \tilphi(R, \xi),
   (2 W_R \partial_R
   + W_{RR}) \partial_\zeta
   \tilphi(R, \zeta)\Bra\Big] \, d\zeta
   \end{equation}
   Using also the bound
   \[
   | \partial_R \partial_\zeta \tilphi(R,\zeta) | \lesssim \left\{
     \begin{array}{cc} R^{\frac12} \ln(1+ R^2) & R < \zeta^{-\frac12} \\
       \zeta^{-\frac14} |\ln \zeta| & R \geq \zeta^{-\frac12}
     \end{array}\right.
   \]
   we can evaluate the inner product in \eqref{eq:letzt} as follows:
\begin{align*}
&\Big|\Bla \partial_\xi^2 \tilphi(R, \xi), (2 W_R \partial_R +
W_{RR})
    \tilphi(R, \eta)\Bra\Big| \less \int_0^{\zeta^{-\frac12} }\la
    R\ra^{-6} R^{\frac72} \log(1+R^2) R^{\frac32} \log(1+R^2)\, dR
    \\
    &+\int_{\zeta^{-\frac12}}^{\xi^{-\frac12}} \la
    R\ra^{-6} R^{\frac72} \log(1+R^2) \zeta^{-\frac14} |\log\zeta|\,R\, dR
    +\int_{\xi^{-\frac12}}^\infty \la
    R\ra^{-6} \xi^{-\frac34}|\log\xi|R^2 \zeta^{-\frac14}|\log\zeta|R\,
    dR
    \less |\log\zeta|^3
\end{align*}
Thus, \eqref{eq:letzt} is controlled by
   \[
   |(\eta-\xi) \partial_\eta^2 F(\xi,\eta)| \lesssim \Big|\int_\xi^\eta
   (\ln \zeta)^3 \,d\zeta\Big| \lesssim \eta |\log \eta|^3
   \]
   Since $\eta \ll \xi$ this yields
   \[
   |\partial_\eta^2 F(\xi,\eta)| \lesssim |\log \eta|^3
   \]
   which concludes the analysis of the off-diagonal part of the kernel.

   Next, we extract the $\delta$ measure that sits on the diagonal of the
   kernel $K$ from the representation formula~\eqref{calk}, see also~\eqref{eq:kern}.
   To do so, we can restrict
    $\xi,\eta$  to a compact subset of
   $(0,\infty)$.  This is convenient, as we then have the following
   asymptotics of $\tilphi(R, \xi)$ for $R \xi^\frac12 \gg 1$:
   \[ \begin{split}
   \tilphi(R,\xi) &=
   \Re \left[a(\xi)\xi^{-\frac{1}{4}}
     e^{iR\xi^{\frac{1}{2}}} \Big(1+
       \frac{3i}{8R\xi^\frac12}\Big)\right] + O(R^{-2})\\
   (R\partial_R -2\xi \partial_\xi) \tilphi(R,\xi) &= -2\Re \left[ \xi
     \partial_\xi(a(\xi)\xi^{-\frac{1}{4}})
     e^{i R\xi^{\frac{1}{2}}} \Big(1+
       \frac{3i}{8R\xi^\frac12}\Big)\right] + O(R^{-2})
       \end{split}
   \]
   where the $O(\cdot)$ terms depend on the choice of the compact
   subset.
   The $R^{-2}$ terms are integrable so they contribute a bounded
   kernel
    to the inner product in~\eqref{calk}. The same applies to the contribution of a
   bounded $R$ region.  Using the above expansions, we
   conclude that the $\delta$-measure contribution of the inner
   product in~\eqref{calk} can only come from one of the following
   integrals:
   \begin{align}
     &  -\int_0^\infty \int_0^\infty f(\xi) \chi(R) \Re
     \left[ \xi \partial_\xi({a}(\xi)\xi^{-\frac{1}{4}})
       a(\eta)\eta^{-\frac{1}{4}}
       e^{iR(\xi^{\frac{1}{2}}+\eta^\frac12)}\Big(1+
       \frac{3i}{8R\xi^\frac12}\Big)\Big(1+
       \frac{3i}{8R\eta^\frac12}\Big)\right]\rho(\xi)\, d\xi dR \label{eq:I1}\\
       &  -\frac12\int_0^\infty \int_0^\infty f(\xi) \chi(R)\,
      \xi \partial_\xi({a}(\xi)\xi^{-\frac{1}{4}})
       \bar{a}(\eta)\eta^{-\frac{1}{4}}
       e^{iR(\xi^{\frac{1}{2}}-\eta^\frac12)}\Big(1+
       \frac{3i}{8R\xi^\frac12}\Big)\Big(1-
       \frac{3i}{8R\eta^\frac12}\Big)\;\rho(\xi)\, d\xi dR \label{eq:I2}\\
     & -\frac12\int_0^\infty \int_0^\infty f(\xi) \chi(R)\,
     \xi \partial_\xi(\bar{a}(\xi)\xi^{-\frac{1}{4}})
       a(\eta)\eta^{-\frac{1}{4}}
       e^{-iR(\xi^{\frac{1}{2}}-\eta^\frac12)}\Big(1-
       \frac{3i}{8R\xi^\frac12}\Big)\Big(1+
       \frac{3i}{8R\eta^\frac12}\Big)\;\rho(\xi)\, d\xi
       dR\label{eq:I3}
   \end{align}
   where $\chi$ is a smooth cutoff function which equals $0$ near
   $R=0$ and $1$ near $R=\infty$.  In all of the above integrals we can
   argue as in the proof of the classical Fourier inversion formula to
   change the order of integration. Integrating by parts in the first
   integral~\eqref{eq:I1} reveals that it cannot
   contribute a $\delta$-measure. Discarding the $R^{-2}$ terms from
   \eqref{eq:I2} and~\eqref{eq:I3} reduces us further to the
   expressions
\begin{align}
     &  -\int_0^\infty \int_0^\infty f(\xi) \chi(R) \Re
     \left[ \xi \partial_\xi({a}(\xi)\xi^{-\frac{1}{4}})
       \bar{a}(\eta)\eta^{-\frac{1}{4}}
       e^{iR(\xi^{\frac{1}{2}}-\eta^\frac12)}\right]\rho(\xi)\, d\xi dR \label{eq:I4}\\
& + \frac38\int_0^\infty \int_0^\infty f(\xi) \chi(R) \Im
     \left[ \xi \partial_\xi({a}(\xi)\xi^{-\frac{1}{4}})
       \bar{a}(\eta)\eta^{-\frac{1}{4}}
       e^{iR(\xi^{\frac{1}{2}}-\eta^\frac12)}\right]R^{-1} (\xi^{-\frac12}-\eta^{-\frac12})\rho(\xi)\, d\xi dR \label{eq:I5}
   \end{align}
   The second integral \eqref{eq:I5} has both an
   $R^{-1}$ and a $(\xi^{-\frac12} -\eta^{-\frac12})$ factor so its
   contribution to $K$ is bounded. The first integral~\eqref{eq:I4}
contributes both a Hilbert transform type kernel  as well as a
$\delta$-measure to $K$. By inspection, the $\delta$ contribution is
\[\begin{split}
&-\frac12 \int_{-\infty}^\infty \Re
     \left[ \xi \partial_\xi({a}(\xi)\xi^{-\frac{1}{4}})
       \bar{a}(\eta)\eta^{-\frac{1}{4}}
       e^{iR(\xi^{\frac{1}{2}}-\eta^\frac12)}\right]\rho(\xi)\, dR
      \\
     &  = -\pi \Re
     \left[ \xi \partial_\xi({a}(\xi)\xi^{-\frac{1}{4}})
       \bar{a}(\eta)\eta^{-\frac{1}{4}}
       \right]\rho(\xi)\delta(\xi^{\frac12}-\eta^{\frac12})\\
       &= -2\pi \xi^{\frac12} \rho(\xi) \Re
     \left[ \xi \partial_\xi({a}(\xi)\xi^{-\frac{1}{4}})
       \bar{a}(\xi)\xi^{-\frac{1}{4}}
       \right]\delta(\xi-\eta)\\
&= -2\pi \xi^{\frac12} \rho(\xi) \Re
     \left[ -\frac14 \xi^{-\frac12}|a(\xi)|^2 + \xi^{\frac12} a(\xi)\bar{a}'(\xi)
       \right]\delta(\xi-\eta)\\
       &= \Big[\frac12 +
       \frac{\xi\rho'(\xi)}{\rho(\xi)}\Big]\delta(\xi-\eta)
\end{split}\]
where we used that $\rho(\xi)^{-1} = \pi |a|^{2}$ in the final step.
Combining this with the $\delta$-measure in \eqref{calk}
yields~\eqref{eq:kern}.
 \end{proof}

 Next we consider the $L^2$ mapping properties for $\calK$. We
 introduce the weighted $L^2$ spaces $L^{2,\alpha}_\rho$ with norm
 \begin{equation}\label{eq:L2weight}
 \| f\|_{L^{2,\alpha}_\rho} := \Big(\int_0^\infty |f(\xi)|^2
 \la \xi\ra^{2\alpha} \rho(\xi)\,d\xi\Big)^{\frac12}
 \end{equation}
 Then we have

\begin{proposition}
 a)  The operator $\calK_0$ from \eqref{eq:kern} maps
  \[
  \calK_0\::\: L^{2,\alpha}_\rho \to L^{2,\alpha+1/2}_\rho
  \]
b) In addition, we have the commutator bound
\[
[ \calK_0,\xi \partial_\xi] \::\:  L^{2,\alpha}_\rho \to
L^{2,\alpha}_\rho
\]
Both statements hold for all $\alpha\in\R$.
\label{l2k}\end{proposition}
\begin{proof}
a) This is equivalent to showing that the kernel
\[
\rho^{\frac12}(\eta)\la \eta\ra^{\alpha+1/2} K_0(\eta,\xi) \la
\xi\ra^{-\alpha} \rho^{-\frac12}(\xi) \::\: L^2(\R^+)\to L^2(\R^+)
\]
With the notation of the previous theorem, the kernel on the
left-hand side is
\[
\tilde K_0(\eta,\xi):= \la\eta\ra^{\alpha+1/2}\la\xi\ra^{-\alpha}
\frac{\sqrt{\rho(\xi)\rho(\eta)}}{\xi-\eta} F(\xi,\eta)
\]
We first separate the diagonal and off-diagonal behavior of $\tilde
K_0$, considering several cases.

  { \bf Case 1: $(\xi,\eta) \in Q := [0,4] \times [0,4]$.}

  We cover the unit interval with dyadic subintervals $I_j =
  [2^{j-1},2^{j+1}]$. We cover the diagonal with the union of squares
  \[
  A = \bigcup_{j=-\infty}^2 I_j \times I_j
  \]
  and divide the kernel $\tilde K_0$ into
  \[
  1_Q\tilde K_0 = 1_{A\cap Q} \tilde K_0 + 1_{Q \setminus A} \tilde K_0
  \]

  { \bf Case 1(a)}: Here we show that the diagonal part $1_{A\cap Q} \tilde K_0$ of
  $\tilde K_0$ maps $L^2$ to $L^2$. By orthogonality it suffices to restrict
  ourselves to a single square $I_j \times I_j$. We recall the $T1$
  theorem for Calderon-Zygmund operators, see page~293 in~\cite{Stein}: suppose the kernel
  $K(\eta,\xi)$ on $\R^2$ defines an operator $T:\calS\to\calS'$ and has the
  following pointwise properties with some $\gamma\in(0,1]$ and a constant $C_0$:
  \begin{enumerate}
  \item[(i)] $|K(\eta,\xi)|\le C_0|\xi-\eta|^{-1}$
  \item[(ii)] $|K(\eta,\xi)-K(\eta',\xi)|\le
  C_0|\eta-\eta'|^\gamma|\xi-\eta|^{-1-\gamma}$ for all
  $|\eta-\eta'|<|\xi-\eta|/2$
\item[(iii)] $|K(\eta,\xi)-K(\eta,\xi')|\le
  C_0|\xi-\xi'|^\gamma|\xi-\eta|^{-1-\gamma}$ for all
  $|\xi-\xi'|<|\xi-\eta|/2$
\end{enumerate}
If in addition $T$ has the restricted $L^2$ boundedness property,
i.e., for all $r>0$ and $\xi_0,\eta_0\in\R$,
$\|T(\omega^{r,\xi_0})\|_2 \le C_0r^{\frac12}$ and
$\|T^*(\omega^{r,\eta_0})\|_2 \le C_0r^{\frac12}$ where
$\omega^{r,\xi_0}(\xi)= \omega((\xi-\xi_0)/r)$ with a fixed
bump-function $\omega$, then $T$ and $T^*$ are $L^2(\R)$ bounded
with an operator norm that only depends on~$C_0$.

Within the
  square $I_j\times I_j$,
  Theorem~\ref{tp} shows that the kernel of $\tilde K_0$ satisfies these properties with $\gamma=1$,
  and is thus bounded on $L^2$.

  { \bf Case 1(b)}: Consider now the off-diagonal part $1_{Q \setminus
    A} \tilde K_0$. In this region, by Theorem~\ref{tp},
    \[
|\tilde K_0(\eta,\xi)| \less
\frac{1}{\sqrt{\xi\eta}|\log\xi||\log\eta|}
    \]
which is a Hilbert-Schmidt kernel on $Q$ and thus $L^2$ bounded.

  {\bf Case 2:} $(\xi,\eta) \in Q^c $. We cover the diagonal with the union of squares
  \[
  B = \bigcup_{j=1}^\infty I_j \times I_j
  \]
  and divide the kernel $\tilde K_0$ into
  \[
  1_{Q^c}\tilde K_0 = 1_{B\cap Q^c} \tilde K_0 + 1_{Q^c \setminus B} \tilde K_0
  \]

  {\bf Case 2a:} Here we consider the estimate on $B$. As in case 1a)
  above, we use Calderon-Zygmund theory. Evidently, $|\tilde
  K_0(\eta,\xi)|\less |\xi-\eta|^{-1}$ on $B$ by Theorem~\ref{tp}.
  To check (ii) and (iii), we differentiate $\tilde K_0$. It will suffice to
  consider the case where the $\partial_\xi$ derivative falls on $F(\xi,\eta)$.
  We distinguish two cases:
if $|\xi^{\frac12}-\eta^{\frac12}|\le 1$, then $|\xi-\eta|\less
\xi^{\frac12}$ which implies that
\[
\frac{\xi^{-\frac12}|\xi-\xi'|}{|\xi-\eta|} \less
\frac{|\xi-\xi'|^{\frac12}}{|\xi-\eta|^{\frac32}} \qquad \forall\;
|\xi-\xi'|<|\xi-\eta|/2
\]
if, on the other hand, $|\xi^{\frac12}-\eta^{\frac12}|> 1$, then
\[
\frac{\xi^{-\frac12}|\xi-\xi'|}{|\xi-\eta||\xi^{\frac12}-\eta^{\frac12}|}
\less \frac{|\xi-\xi'|}{|\xi-\eta|^2} \qquad \forall\;
|\xi-\xi'|<|\xi-\eta|/2
\]
which proves property (iii) on $B$ with $\gamma=\frac12$, and by
symmetry also (ii). The restricted $L^2$ property follows form the
cancelation in the kernel and the previous bounds on the kernel.
Hence, $\tilde K_0$ is $L^2$ bounded on~$B$.

  {\bf Case 2b:} Finally, in the exterior region $Q_c \setminus B$ we
  have the bound, with arbitrarily large $N$,
  \[
  |\tilde K_0(\eta,\xi)| \lesssim  (1+\xi)^{-N} (1+\eta)^{-N}
  \]
  which is $L^2$ bounded by Schur's lemma.

b) A direct computation shows that the kernel $K_0^{com}$
of the commutator $[\xi \partial_\xi, K_0]$  is given by
\[
K_0^{com}(\eta,\xi) = (\eta \partial_\eta + \xi \partial_\xi)
K_0(\eta,\xi) + K_0(\eta,\xi) =
\frac{\rho(\xi)}{\xi-\eta}F^{com}(\xi,\eta)
\]
interpreted in the principal value sense and with $F^{com}$ given by
\[
F^{com}(\xi,\eta) = \frac{\xi \rho'(\xi)}{\rho(\xi)} F(\xi,\eta) +
(\eta
\partial_\eta + \xi \partial_\xi) F(\xi,\eta)
\]
By Theorem~\ref{tp} this satisfies the same pointwise off-diagonal
bounds as $F$. Near the diagonal the bounds for $F^{com}$ and its
derivatives are worse\footnote{The one derivative loss can be
avoided by a more careful analysis, but this does not seem necessary
here.} than those for $F$ by a factor of $(1+\xi)^\frac12$. Then the
proof of the $L^2$ commutator bound is similar to the argument in
part~(a).
\end{proof}




\section{ The final equation}

To rewrite the equation \eqref{equation2} in a final form, we begin
by expressing the operator  $R \partial_R$ in terms of the kernel
$\calK$ in the transference identity~\eqref{eq:transfer}. We have,
with $\calF$ as in Theorem~\ref{thm:GZ},
\[
 \calF \Big(\partial_{\tau}+\frac{\lambda_{\tau}}{\lambda}R\partial_{R}\Big) =
\Big(\partial_{\tau} + \frac{\lambda_{\tau}}{\lambda}(- 2\xi
\partial_{\xi}
 + \calK) \Big) \calF
\]
which gives
\[
\begin{split}
\calF
\Big(\partial_{\tau}+\frac{\lambda_{\tau}}{\lambda}R\partial_{R}\Big)^2
=  &
 \Big(\partial_{\tau} + \frac{\lambda_{\tau}}{\lambda}(-2\xi \partial_{\xi}
 + \calK) \Big)^2 \calF
\\ = &  \Big(\partial_{\tau} - \frac{\lambda_{\tau}}{\lambda}2\xi
\partial_{\xi}\Big)^2  \calF+
2\frac{\lambda_{\tau}}{\lambda} \calK  \Big(\partial_{\tau} -
\frac{\lambda_{\tau}}{\lambda}2\xi
\partial_{\xi}\Big)  \calF+
\frac{\lambda_{\tau}^2}{\lambda^2}(\calK^2 + 2[\xi \partial_\xi,\calK]) \calF
\end{split}
\]
This leads to a transport type equation for the Fourier transform
$x(\tau,\xi)$ of $\tileps$ by applying $\calF$ to~\eqref{equation2}.
Indeed, in view of the preceding
\begin{equation}\label{final}
\begin{split}
 -\Big(\partial_{\tau} - \frac{\lambda_{\tau}}{\lambda}2\xi
\partial_{\xi}\Big)^2  x  - \xi x  =& 2\frac{\lambda_{\tau}}{\lambda}
\calK  \Big(\partial_{\tau} - \frac{\lambda_{\tau}}{\lambda}2\xi
\partial_{\xi}\Big) x + \frac{\lambda_{\tau}^2}{\lambda^2}(\calK^2 +
2[\xi \partial_\xi,\calK]) x \\
&- \Big(\frac{1}{4}(\frac{\lambda_{\tau}}{\lambda})^{2}
+\frac{1}{2}\partial_{\tau}(\frac{\lambda_{\tau}}{\lambda})\Big) x
 + \lambda^{-2} \calF R^{\frac12} (N_{2k-1} (R^{-\frac12}
\calF^{-1} x) +  e_{2k-1})
\end{split}
\end{equation}

We want to obtain solutions to \eqref{final} which decay as $\tau
\to \infty$, which means we need to solve the equation backward in
time, i.e., with zero Cauchy data at $\tau = \infty$. We treat this
problem iteratively, as a small perturbation of the linear equation
governed by the operator on the left--hand side. For this we need to
solve the following {\bf{transport equation}}
\begin{equation}\label{transport}
  -\Big[\Big(\partial_{\tau}-
2\frac{\lambda_{\tau}}{\lambda}\xi\partial_{\xi}\Big)^{2}+\xi\Big]x(\tau,
\xi)=b(\tau, \xi),
\end{equation}
We denote by $H$ the backward fundamental solution for
the operator
\[
 \Big(\partial_{\tau} - 2\frac{\lambda_{\tau}}{\lambda}\xi
\partial_{\xi}\Big)^2   +\xi
\]
and by $H(\tau,\sigma)$ its kernel,
\[
x(\tau) = \int_{\tau}^\infty H(\tau,\sigma) f(\sigma)\, d\sigma
\]
The mapping properties of $H$ are described in the following result,
which will be proven in the next section.
\begin{proposition} For any $\alpha\ge 0$
there exists some (large) constant $C=C(\alpha)$ so that the
operator $H(\tau,\sigma)$ satisfies the bounds
\begin{align}
\|H(\tau,\sigma) \|_{ L^{2,\alpha}_\rho \to L^{2,\alpha+1/2}_\rho}
&\lesssim \tau \Big(\frac{\sigma}{\tau}\Big)^{C} \label{hbd}\\
\Big\| \Big(\partial_{\tau} - \frac{\lambda_{\tau}}{\lambda}2\xi
\partial_{\xi}\Big) H(\tau,\sigma) \Big\|_{ L^{2,\alpha}_\rho \to L^{2,\alpha}_\rho}
&\lesssim  \Big(\frac{\sigma}{\tau}\Big)^{C}
\label{htaubd}\end{align} uniformly in $\sigma\ge\tau$.
\label{proph}\end{proposition}

\noindent This leads us to introduce the spaces $L^{\infty,N}
L^{2,\alpha}_\rho$ with norm
\[
\| f\|_{L^{\infty,N} L^{2,\alpha}_\rho} := \sup_{\tau\ge1} \tau^N
\|f(\tau)\|_{L^{2,\alpha}_\rho}
\]
Then the above proposition immediately allows us to draw the
following conclusions:
\begin{cor} Given $\alpha\ge0$, let $N$ be large enough. Then
\[
\| H b\|_{L^{\infty,N-2} L^{2,\alpha+1/2}_\rho}  + \Big\|
\Big(\partial_{\tau} - 2\frac{\lambda_{\tau}}{\lambda}\xi
\partial_{\xi}\Big)H
b\Big\|_{L^{\infty,N-1} L^{2,\alpha}_\rho} \le C_0\,\frac{1}N
\|b\|_{L^{\infty,N} L^{2,\alpha}_\rho}
\]
with a constant $C_0$ that depends on $\alpha$ but does not depend
on $N$.
\end{cor}

\noindent The small factor $N^{-1}$ is crucial here for our argument
to work. On the other hand, the nonlinear operator $N_{2k-1}$
from~\eqref{final} has the following mapping properties:

\begin{proposition}
Assume that $N$ is large enough and $\frac{\nu}{2}+\frac34 > \alpha
> \frac14$. Then the map
\[
x \to \lambda^{-2} \calF R^\frac12 (N_{2k-1} (R^{-\frac12}
\calF^{-1} x))
\]
is locally Lipschitz from $ L^{\infty,N-2} L^{2,\alpha+1/2}_\rho$ to
$L^{\infty,N} L^{2,\alpha}_\rho$. \label{propn}
\end{proposition}

The above two results, combined with  Proposition~\ref{l2k}, allow
us to use a contraction argument to solve equation \eqref{final}.
The next two sections are devoted to proving
Propositions~\ref{proph}, \ref{propn}. Finally, in the last section
we close the argument.

\section{The transport equation}

Here we consider the backward fundamental solution $H$ for
\eqref{transport} and prove Proposition~\ref{proph}.  Observe that
\eqref{transport} implies
\begin{equation}\nonumber
  [\partial_{\tau}^{2}+\lambda^{-2}(\tau)\xi]x(\tau,
  \lambda^{-2}(\tau)\xi)=b(\tau,\lambda^{-2}(\tau)\xi)
\end{equation}
We introduce the operator
\[
L_{\xi,\tau}:=\partial_{\tau}^{2}+\lambda^{-2}(\tau)\xi
\]
and the fundamental solutions $S(\tau, \sigma ,\xi)$, $U(\tau, \sigma,
\xi)$, which satisfy
\begin{equation}\nonumber
  L_{\xi,\tau}S(\tau, \sigma, \xi)=0,\,S(\tau, \tau,
  \xi)=0,\,\partial_{\tau}S(\tau, \sigma, \xi)|_{\tau=\sigma}=-1
\end{equation}
\begin{equation}\nonumber
  L_{\xi,\tau}U(\tau, \sigma, \xi)=0,\,U(\tau, \tau,
  \xi)=1,\,\partial_{\tau}U(\tau, \sigma, \xi)|_{\tau=\sigma}=0
\end{equation}
Then \eqref{transport} may be solved via
\begin{equation}\nonumber
  x(\tau, \lambda^{-2}(\tau) \xi)=-\int_{\tau}^{\infty}
  S(\tau, \sigma, \xi)b(\sigma,
  \lambda^{-2}(\sigma) \xi)\,d\sigma
\end{equation}
Given this representation, we note that the index $\alpha$ plays no
role in \eqref{hbd} and \eqref{htaubd} since
\[
\frac{(1+ \lambda^{-2}(\tau) \xi)^\alpha}{(1+ \lambda^{-2}(\sigma)
  \xi)^\alpha} \lesssim \left(\frac{\sigma}{\tau}\right)^{C}
\]
Hence without loss of generality we set $\alpha = 0$.  Similarly, we
can neglect the measure of integration $\rho(\xi)d\xi$ which also
has a polynomial behavior both at $0$ and at infinity,
\[
\frac{\rho(\lambda^{-2}(\tau) \xi)}{\rho(\lambda^{-2}(\sigma) \xi)}
\lesssim \Big(\frac{\sigma}{\tau}\Big)^C
\]
Then the bounds \eqref{hbd} and \eqref{htaubd} reduce to proving that
\[
|S(\tau, \sigma, \xi)| \lesssim \sigma
\left(\frac{\sigma}{\tau}\right)^{C} (1+ \lambda^{-2}(\tau)
\xi)^{-\frac12}, \qquad |\partial_\tau S(\tau, \sigma, \xi)| \lesssim
\left(\frac{\sigma}{\tau}\right)^{C}, \qquad 1 \lesssim \tau < \sigma
\]
Recalling that $\lambda(\tau) = \tau^{1+ \nu^{-1}}$ (we are ignoring
a multiplicative constant here), we strengthen the first bound and
prove instead that
\begin{equation}
  |S(\tau, \sigma, \xi)| \lesssim \sigma
  \left(\frac{\sigma}{\tau}\right)^{C}   (1+ \tau^{-\frac{2}{\nu}}
  \xi)^{-\frac12}, \qquad
  |\partial_\tau S(\tau, \sigma, \xi)| \lesssim
  \left(\frac{\sigma}{\tau}\right)^{C}   \qquad 0 < \tau < \sigma
  \label{S}\end{equation}
The advantage of doing this is that the last bound is scale invariant.
Precisely, one verifies directly the scaling relation
\[
S(\tau,\sigma,\xi) = \xi^{\frac{\nu}2} S(\tau
\xi^{-\frac{\nu}2},\sigma \xi^{-\frac{\nu}2},1)
\]
which leaves \eqref{S} unchanged. Hence in what follows it suffices
to prove \eqref{S} in the case $\xi =1$. We begin by constructing
two special solutions for the operator $L_{1,\tau}$. For
small\footnote{The reader should bear in mind that by this $\tau$ we
mean the rescaled one, i.e, $\xi^{-\frac{\nu}{2}}\tau$, which can be
arbitrarily close to zero.} $\tau$ we use a standard WKB ansatz.

\begin{lemma}

  a) (Large $\tau$ solutions) If $\nu$ is not an even integer
  then there exist two analytic solutions $\phi_0$ and $\phi_1$ of $
  L_{1,\tau}\phi_j=0$  with a series representation
  \[
  \phi_j (\tau) = \sum_{k =0}^\infty c_{jk}\, \tau^{j-\frac{2k}{\nu}}, \qquad
  c_{j0} = 1
  \]
  which is convergent for all $\tau > 0$. If $\nu$ is an even integer
then the result still holds with a modification in the expression
for $\phi_1$, namely
\[
 \phi_1 (\tau) = c_1 \phi_0(\tau) \ln \tau+
 \sum_{k =0}^\infty c_{1k}\, \tau^{1-\frac{2k}\nu }, \qquad
  c_{10} = 1
  \]

  b) (Small $\tau$ solutions) There is a  solution $\phi_2$
  for $ L_{1,\tau}$ of the form
  \[
  \phi_2(\tau) = \tau^{\frac12 +\frac{1}{2\nu}} e^{i \nu
    \tau^{-\frac{1}{\nu}}} [1+ a(\tau^\frac1{\nu})]
  \]
  with $a$ smooth and satisfying  $a(0) = 0$.
\end{lemma}

\begin{proof}
a) We substitute the formal series in the equation
\[
(\partial_\tau^2 + \tau^{-2 - \frac{2}{\nu}}) \phi_j = 0
\]
in the equation and identify the coefficients of the similar terms.
This yields
\[
c_{j,k} (j-\frac{2 k}\nu)(j-1- \frac{2 k}\nu) + c_{j,k-1} = 0 \qquad k \geq 1
\]
Hence the coefficients $c_{jk}$ can be iteratively computed and
satisfy a bound of the type
\[
|c_{j,k}| \le  \frac{C^k}{(k!)^2}
\]
which implies that the series converges for all $\tau$.

If $j=0$ then the argument works for all $\nu > 0$. If $j =1$ then
there is an obstruction if $\nu$ is an even integer; indeed, this
happens  precisely when $2k = \nu$. As usual, this is compensated
for by adding in the logarithmic term, since
\[
L_{1,\tau} (\phi_0(\tau) \log \tau) = -\tau^{-2} \phi_0 + \tau^{-1}
2\partial_\tau \phi_0
\]
has a nonzero coefficient on the $\tau^{-2}$ term.

b) In this case, we use the usual WKB-ansatz which we now recall in
a more general setting: we wish to solve the equation
$(\partial_\tau^2 + Q)\psi=0$ where $Q(\tau)$ is a smooth potential
for $\tau>0$. Fix some (small) $\tau_0>0$. WKB means that we seek a
solution of the form $\psi(\tau)=\psi_0(\tau)[1+a(\tau)]$ with
\[ \psi_0(\tau)= Q^{-\frac14}(\tau)e^{iS(\tau)},\qquad
S(\tau)=\int_{\tau_0}^\tau Q^{\frac12}(\sigma)\,d\sigma\] Since
\[ \partial_\tau^2\psi_0 + Q\psi_0 = V\psi_0, \qquad
V= -\frac14 \frac{Q''}{Q} +\frac{5}{16} \Big(\frac{Q'}{Q}\Big)^2 \]
we obtain the following equation for $a(\tau)$:
\[
(a'\psi_0^2)'(\tau) = -V\psi_0^2(\tau) [1+a(\tau)]
\]
which we solve in the form
\[\begin{split}
a(\tau) &= -\int_0^\tau \int_0^{\tau'} \psi_0^{-2}(\tau')
\psi_0^2(\sigma)  V(\sigma) [1+a(\sigma)]\,d\sigma d\tau'\\
&= \frac{i}{2} \int_0^\tau Q^{-\frac12}(\sigma)
\big[1-e^{2i(S(\sigma)-S(\tau))}\big]\, V(\sigma)
[1+a(\sigma)]\,d\sigma
\end{split}\]
provided these integrals converge at zero. They do in our case: in
fact, $Q(\tau)=\lambda^{-2}(\tau)$ which implies that
\[\begin{split}
\psi_0(\tau) &= \tau^{\frac12+\frac{1}{2\nu}} e^{i\nu
\tau^{-\frac{1}{\nu}}}\\
 a(\tau) &= c\, i\int_0^\tau
\sigma^{-1+\frac{1}{\nu}}
\big[1-e^{2i\nu(\sigma^{-\frac{1}{\nu}}-\tau^{-\frac{1}{\nu}})}\big]\,
[1+a(\sigma)]\,d\sigma\end{split}
\]
or, after changing variables to $a(\tau^\nu)= \tilde a(\tau)$,
\begin{equation}\label{eq:atilde}
\tilde a(\tau) = ic\nu\int_0^\tau
\big[1-e^{2i\nu(\sigma^{-1}-\tau^{-1})}\big]\, [1+\tilde
a(\sigma)]\,d\sigma
\end{equation}
By the boundedness of the kernel,  this Volterra equation has a
solution $\tilde a\in C([0,\infty))$ which is clearly then also
smooth for all $\tau>0$. We now claim that in fact $\tilde a\in
C^\infty([0,\infty))$.  Indeed, the zero order iterate here is a
smooth function at $\tau=0$:
\begin{align*}
\int_0^\tau \big[1-e^{2i\nu(\sigma^{-1}-\tau^{-1})}\big]\,d\sigma &=
\int_{\tau^{-1}}^\infty
\big[1-e^{2i\nu(u-\tau^{-1})}\big]\,\frac{du}{u^2}   \\
& = \tau - \int_{\tau^{-1}}^\infty
e^{2i\nu(u-\tau^{-1})}\,\frac{du}{u^2} = \sum_{j=1}^m c_j\, \tau^j +
O(\tau^{m+1})
\end{align*}
for any positive integer $m$ by repeated integration by parts. One
now proceeds to show the same for the higher Volterra iterates;
alternatively, we insert the ansatz
\[\tilde a(\tau) = \sum_{j=1}^m d_j\, \tau^j +
O(\tau^{m+1})
\]
into \eqref{eq:atilde}  and solve for the coefficients $d_j$. In
either case, the conclusion is that \eqref{eq:atilde} has a smooth
solution, as claimed.
\end{proof}

We now use this lemma to prove \eqref{S}, which will then conclude
the proof of Proposition~\ref{proph}. Considering the limits at
infinity, respectively at $0$, one finds that
\[
W(\phi_0,\phi_1) = 1, \qquad W(\phi_2,\overline\phi_2) = -2i
\]
This allows us to express the backward fundamental solution
$S(\tau,\sigma)$ in terms of these bases. Note that we suppress the
$\xi$ variable as $\xi=1$ is fixed. We consider two cases.

{\bf Case 1:} $\sigma > 1$. Then we have
\[
S(\tau,\sigma) = \phi_1(\sigma) \phi_0(\tau) - \phi_0(\sigma)
\phi_1(\tau)
\]
If $1 \leq \tau \leq \sigma$, then \eqref{S} follows directly from
the properties of $\phi_0$ and $\phi_1$. If $\tau < 1$ then we
express $\phi_0(\tau)$ and $\phi_1(\tau)$ in terms of the
$\{\phi_2,\overline\phi_2\}$ basis to obtain
\[
S(\tau,\sigma) = \Re (c(\sigma) \phi_2(\tau)), \qquad |c(\sigma)|
\lesssim \sigma
\]
This gives
\[
|S(\tau,\sigma)| \lesssim \sigma \tau^{\frac12+\frac1{2\nu}},
\qquad  |\partial_\tau S(\tau,\sigma)| \lesssim
\sigma \tau^{-\frac12-\frac1{2\nu}}
\]
Again \eqref{S} follows.

{\bf Case 2:} $\sigma < 1$. Then we express $S(\tau,\sigma)$ in the
$\{\phi_2$, $\overline \phi_2\}$ basis to obtain
\[
S(\tau,\sigma) = \Im (\phi_2(\sigma) \overline\phi_2(\tau))
\]
This gives the bounds
\[
|S(\tau,\sigma)| \lesssim \sigma^{\frac12+\frac1{2\nu}}
 \tau^{\frac12+\frac1{2\nu}},
\qquad  |\partial_\tau S(\tau,\sigma)| \lesssim
\sigma^{\frac12+\frac1{2\nu}}  \tau^{-\frac12-\frac1{2\nu}}
\]
which imply \eqref{S}.

\section{The nonlinear terms}

  In this section we consider the {\it{nonlinear source
      terms}}, i.e., those given by the right-hand side of
  \eqref{equation2}, and prove Proposition~\ref{propn}.
Recalling that $R=r\lambda$,  we write
\begin{equation}\label{eq:N_split}\begin{split} \lambda^{-2} R^\frac12 N_{2k-1}(
R^{-\frac12} \tileps) &= \frac{\cos(2u_{2k-1}) -\cos( 2Q)}{R^2}
    2 {\tileps}+ \frac{\sin(2u_{2k-1})}{2R} \frac{\cos
(2 \tileps{R^{-\frac{1}{2}}})-1}{R^{\frac{1}{2}}}\\
& \quad + \cos (2u_{2k-1}) \frac{\sin
  (2{\tileps}{R^{-\frac{1}{2}}})-2{\tileps}{R^{-\frac{1}{2}}}}{2R^{\frac{3}{2}}}
  \end{split}\end{equation}
where the regularity of the coefficients above is computed as in
Step~2 of the proof of Theorem~\ref{thm:sec2},
\begin{align}
 \frac{\cos(2u_{2k-1}) -\cos (2Q)}{R^2} &\in \tau^{-2} \IS^2(R^{-2} (\log
 R)^2,\calQ_{k-1})\label{eq:nicht1}\\
 \frac{\sin(2u_{2k-1})}{2R} &\in \IS^0(R^{-2} \log
 R,\calQ_{k-1})\label{eq:nicht2}\\
 \cos (2u_{2k-1}) &\in  \IS^0(1,\calQ_{k-1})\label{eq:nicht3}
\end{align}
where we used here that $t\lambda(t)\asymp \tau$ and also that
$R\lesssim \tau$ (recall the algebras $\cQ$ and $\cQ_{k}$ from
Definition~\ref{def:Q}). Proposition~\ref{propn} amounts to proving
multiplicative estimates in the context of the classical Sobolev
spaces. Here we use Sobolev spaces adapted to the operator $\calL$,
namely
\[
\| u\|_{H^\alpha_\rho} := \|\hat{u}\|_{L^{2,\alpha}_\rho}
\]
Restating Proposition~\ref{propn} with this notation shows that we
need to prove that the map
\[
\tileps \mapsto \lambda^{-2}R^\frac12 N_{2k-1}( R^{-\frac12}
\tileps)
\]
is locally Lipschitz from $L^{\infty,N-2} H^{\alpha+1/2}_\rho$ to
$L^{\infty,N} H^{\alpha}_\rho$. Note that~\eqref{eq:nicht1} has an
explicit gain of $\tau^{-2}$ which explains why we can improve the
time-decay of the first (linear) term on the right-hand side
of~\eqref{eq:N_split} from $N-2$ to~$N$. On the other hand, there is
no such gain in~\eqref{eq:nicht2} and~\eqref{eq:nicht3}. What saves
us here is that both the second and third terms on the right-hand
side of~\eqref{eq:N_split} are truly {\em nonlinear terms}
in~$\tilde\eps$.

\noindent As a technical tool we introduce an inhomogeneous
Littlewood-Paley decomposition
\[
f = \sum _{\lambda =1}^\infty P_\lambda f = \sum_\lambda
\int_0^\infty p_\lambda(\xi) \phi(R,\xi)\hat{f}(\xi)\rho(\xi)\,d\xi
\]
corresponding to a smooth partition of unity $\{p_\lambda\}$ in the
Fourier space. Here $\lambda\in\{2^j\}_{j=0}^\infty$ and $p_\lambda$
is adapted to frequencies of size~$\lambda$. Our first result is
\begin{lemma}
Let $q \in S(1,\calQ)$ and $|\alpha| < \frac{\nu}{2}+\frac34$. Then
\[
\| q f\|_{H^\alpha_\rho} \lesssim \|f\|_{H_\rho^\alpha}
\]
\end{lemma}
\begin{proof}
We decompose the multiplication operator into its Littlewood-Paley
pieces:
\[
q = \sum_{\lambda,\mu} P_\lambda\, q P_\mu
\]
The diagonal sum corresponding to $\lambda \asymp \mu$ is estimated
using only the $L^\infty$ bound on~$q$.  For the off-diagonal
component it suffices to show rapid decay. In fact, we claim that
\[
\|  P_\lambda\, q P_\mu\|_{L^2 \to L^2} \lesssim
(\mu+\lambda)^{-\frac14-\frac{\nu}2} [\ln(\mu+\lambda)]^m
\]
where $m$ is some large integer. The Fourier kernel of $ P_\lambda
\, q P_\mu$ is
\[
K_{\lambda,\mu}(\eta,\xi)=  \sqrt{\rho(\xi)\rho(\eta)}\;
p_\lambda(\xi) p_\mu(\eta)
  \int q(R) \tilphi(\xi,R) \tilphi(\eta,R)\, dR\]
  in the sense that
  \[ \sqrt{\rho(\eta)}\;\calF(P_\lambda
\, q P_\mu f)(\eta)=\int
K_{\lambda,\mu}(\eta,\xi)\,\hat{f}(\xi)\sqrt{\rho(\xi)}\,d\xi
\]
Therefore, the above $L^2$ bound would follow from the pointwise
estimate (recall $\rho(\xi)\asymp\xi$ for $\xi>1$)
\[
\left|  \int q(R) \tilphi(\xi,R) \tilphi(\eta,R) dR \right| \lesssim
\la\xi\ra^{-1} \la\eta\ra^{-1} \la\xi+\eta\ra^{
-\frac14-\frac{\nu}2} [\ln(2+ \xi+\eta)]^m
\]
The function $q$ has a symbol type behavior with respect to $R$
except near $R = \tau$, where it has a power type singularity
$(1-a)^{\nu+\frac12}$, $a = R/\tau$, possibly involving also
logarithms\footnote{Strictly speaking, there is a multiplicative
constant in $a=cR/\tau$,  but we ignore it}.  To separate this
singularity from the behavior at $0$ we use a smooth cutoff to split
$q$ into (recall that $\tau$ is a large parameter)
\[
q = q_{<\tau/2} + q_{>\tau/2}
\]
The first term is a  symbol of order $0$ with respect to $R$. To
proceed, we recall the calculations leading up to~\eqref{Kk}. The
main tool  there is the following double commutator identity: if
$\xi\ne\eta$ and $U$ is a zero order symbol, then
\begin{equation}\label{eq:doppel}\begin{split} &(\xi-\eta)^{2} \la
U(R)\phi(R,\xi),\phi(R,\eta)\ra =
\Bla \big[[U,\cL],\cL\big] \phi(R,\xi),\phi(R,\eta)\Bra\\
& =\Bla \big(-4U_{RR}\xi + 3R^{-2}(U_{RR}-R^{-1}U_R) + 4U_{RR} V +
U_{RRRR} + 2U_RV_R + 4U_{RRR}\partial_R\big)
\phi(R,\xi),\phi(R,\eta)\Bra
\end{split}\end{equation}
where the inner products exist in the principal value sense (recall
that $V(R)=-8(1+R^2)^{-2}$). Iterating this identity $k$ times
yields
\[
(\xi-\eta)^{2k} \Bla q_{<\tau/2}(R)  \tilphi(R,\xi), \tilphi(R,\eta)
\Bra =  \Bla \Big[ \sum_{j=0}^{k-1} \xi^{j }q_j^{odd}(R)
\partial_R +\sum_{\ell=0}^k \xi^\ell q_\ell^{even}(R)\Big] \tilphi(R,\xi), \tilphi(R,\eta) \Bra
\]
where $q_j^{odd}$ and  $q_\ell^{even}$ are symbols of order at most
$-2k$ with odd, respectively even, expansions around $R=0$. For
$1+\xi \not\asymp 1+\eta$ this gives
\[
|\la q_{<\tau/2}(R)  \tilphi(R,\xi), \tilphi(R,\eta) \ra| \lesssim
\la\xi+\eta\ra^{-k}
\]
for all $k$ which is more than we need.

\noindent The second term $q_{>\tau/2}$ can be thought of as a
function of $a$,
\[
q_{>\tau/2}(R)  = q_1(a), \qquad a = \frac{R}{\tau}
\]
where $q_1$ is supported in $[\frac12,2]$ and has a $\calQ$ type
singularity\footnote{$q_1$ also has a nonsingular part, which by a
  slight abuse of notation we include in $ q_{<\tau/2} $} at $a=1$. We
divide it into a singular and a nonsingular component,
\[
q_{>\tau/2} = q_{>\tau/2}^{{\rm s}} + q_{>\tau/2}^{{\rm ns}}, \qquad
q_{>\tau/2}^{{\rm s}} := q_{>\tau/2} \chi_{[|R-\tau| <
\la\xi+\eta\ra^{-\frac12}]}, \qquad q_{>\tau/2}^{{\rm ns}} :=
q_{>\tau/2} \chi_{[|R-\tau| > \la\xi+\eta\ra^{-\frac12}]},
\]
where the  $\chi$'s define a smooth partition of unity relative to
the indicated sets. For the singular component we bound the integral
directly using the pointwise bounds on $\tilphi(R,\xi)$ to obtain
\[
\begin{split}
\left| \int q_{>\tau/2}^{{\rm s}} (R) \tilphi(R,\xi)
\tilphi(R,\eta)\, dR \right|  & \lesssim \int_{\frac{\tau}{2}}^\tau
(1-R/\tau)^{\nu+\frac12} |\log(1-R/\tau)|^m \;1_{[|R-\tau| <
\la\xi+\eta\ra^{-\frac12}]}
\la\xi\ra^{-\frac34}\la\eta\ra^{-\frac34}\, dR\\
&\lesssim \la \xi\ra^{-\frac34}\la\eta\ra^{-\frac34}
\tau^{-\nu-\frac12} \la\xi+\eta\ra^{-\frac{\nu}{2}-\frac34}[\log(2+\xi+\eta)]^m\\
& \lesssim \la \xi\ra^{-1}\la\eta\ra^{-1}
\la\xi+\eta\ra^{-\frac{\nu}{2}-\frac14} [\log(2+\xi+\eta)]^m
\end{split}
\]
For the nonsingular component, a $k$-fold iteration
of~\eqref{eq:doppel} yields
\begin{equation}
(\xi-\eta)^{2k} \Bla q_{>\tau/2}^{{\rm ns}}(R)  \tilphi(R,\xi),
\tilphi(R,\eta) \Bra =  \Bla\Big[ \sum_{j=0}^{k-1} \xi^j
q_{k,j}^{{\rm odd}}(R)\partial_R + \sum_{\ell=0}^k \xi^\ell
 q^{{\rm even}}_{k,\ell}(R) \Big ]\tilphi(R,\xi), \tilphi(R,\eta) \Bra
\label{Kk2}\end{equation} with
\[
q^{{\rm odd}}_{k,j}(R) = \sum_{i=0}^{2k-j-1} r_{k,j,i}^{{\rm
odd}}(R)\;
\partial_R^{2i+1} q_{>\tau/2}^{{\rm ns}}(R),\qquad q^{{\rm
even}}_{k,\ell}(R) = \sum_{i=1}^{2k-\ell} r_{k,\ell,i}^{{\rm
even}}(R)\;
\partial_R^{2i} q_{>\tau/2}^{{\rm ns}}(R)
\]
where the coefficients are rational functions, smooth for all
$R\ge0$, decaying at rates
\[
|r_{k,j,i}^{{\rm odd}}(R)|\less R^{-2-(4k-2j-2i)}, \qquad
|r_{k,\ell,i}^{{\rm even}}(R)|\less R^{-4-(4k-2\ell-2i)}
\]
The logic behind the numerology here is simple: a factor $\xi^j$
consumes $2j$ derivatives, so the remaining $4k$ derivatives need to
hit either the symbol $q_{>\tau/2}^{{\rm ns}}(R)$ or the weight $V$
(the latter leading to the rational functions).

We show how to apply these formulas for the case of the even
weights, the odd ones being analogous. As for the derivatives
\[
\partial_R^{2i} q_{>\tau/2}^{{\rm ns}}(R) = \partial_R^{2i}\Big(q_{>\tau/2} \chi_{[|R-\tau| >
\la\xi+\eta\ra^{-\frac12}]}\Big)
\]
it will suffice to consider two extreme cases: when all derivatives
fall on the symbol, or all fall on the cut-off function. The
contribution by the latter to $|\la q_{>\tau/2}^{{\rm ns}}(R)
\tilphi(R,\xi), \tilphi(R,\eta) \ra|$ is bounded by (ignoring logs)
\begin{align*}
&(\xi+\eta)^{-2k} \il_{[|R-\tau|\asymp\la \xi+\eta\ra^{-\frac12}]}
R^{-4-(4k-2\ell-2i)} (1-a)^{\nu+\frac12} \la \xi+\eta\ra^i \xi^\ell
\la\xi\ra^{-\frac34}\la\eta\ra^{-\frac34}\, dR \\
& \less (\xi+\eta)^{-2k}\tau^{-3-(4k-2\ell-2i)}\tau^{-\nu-\frac32}
\la\xi+\eta\ra^{-\frac{\nu}{2}-\frac34 +i}
\la\xi\ra^{\ell-\frac34}\la\eta\ra^{-\frac34}\\
&\lesssim \la\xi\ra^{-1}\la\eta\ra^{-1}
\la\xi+\eta\ra^{-\frac{\nu}{2}-\frac14}
\end{align*}
as desired. The other cases are checked similarly and we skip them.
\end{proof}

This allows us to deal with the coefficients in front of the
$\tileps$ terms. As remarked above,  the $\tau$ decay for the first
term in $N_{2k-1}$ comes from the $\tau^{-2}$ factor in the
coefficient and from the quadratic (respectively, cubic) expressions
in $\tileps$ for the remaining terms. It remains to prove the
following:

\begin{proposition}\label{prop:uff}
Let $\alpha > \frac14$. Then the  maps \begin{align} \tileps
&\mapsto R^{-\frac12} (\cos (2 \tileps{R^{-\frac{1}{2}}})-1)\label{eq:epsmap1}\\
\tileps &\mapsto R^{-\frac32}(\sin
(2{\tileps}{R^{-\frac{1}{2}}})-2{\tileps}{R^{-\frac{1}{2}}})\label{eq:epsmap2}
\end{align}
are locally Lipschitz from $H^{\alpha+1/2}_\rho$ to
$H^\alpha_\rho$.
\end{proposition}

The proof will be split up into the following four lemmas.  We first
obtain a pointwise bound for frequency localized $L^2$ functions:

\begin{lemma}
For dyadic $\lambda \ge 1$ we have
\[
| P_\lambda f(R)| \lesssim \lambda
\min\{R^\frac32,\lambda^{-\frac34}\} \|f\|_{L^2}
\]
for all $f\in L^2(\R^+)$. \label{linf}\end{lemma}

\begin{proof}
Using the inversion formula we write
\[
P_\lambda f(R) = \int_0^\infty p_\lambda(\xi) \hat f(\xi)
\tilphi(R,\xi) \rho(\xi)\, d\xi
\]
The  pointwise bounds for $ \tilphi$,
\[
|\tilphi(R,\xi)| \lesssim \min\{R^\frac32,\xi^{-\frac34}\}
\]
and the Cauchy-Schwarz inequality finish the proof.
\end{proof}

\noindent We also have estimates for the derivative:

\begin{lemma}
For dyadic $\lambda \ge 1$ we have
\[
| \partial_R P_\lambda f(R)| \lesssim \lambda
\min\{R^\frac12,\lambda^{-\frac14}\} \|f\|_{L^2}
\]
and
\[
\| \partial_R P_\lambda f\|_{L^2} \lesssim \lambda^\frac12
\|f\|_{L^2}
\]
for all $f\in L^2(\R^+)$. \label{dpl}\end{lemma}
\begin{proof}
  The first estimate follows from the pointwise bounds on $\partial_R
  \tilphi$. For the second bound we can integrate by parts (justified
by the first bound) to obtain
\[
\lambda \|P_\lambda f(R)\|_{L^2}^2 \gtrsim\Bla \cL f,f\Bra \ge \|
\partial_R f\|_{L^2}^2 + \frac34 \|R^{-1} f\|^2_{L^2} - C\|f\|^2_{L^2}
\]
which leads to the desired conclusion.
\end{proof}

Next we consider bilinear estimates but with a weight that is
singular at $0$. This suffices in order to estimate the quadratic
and the cubic terms in the proposition. The logic behind the
Lemma~\ref{lem:8.5} is the following: dividing by $R^{\frac32}$
should amount to a loss of $\xi^{\frac34}$ on the Fourier side
(since the scaling relation is $R\xi^{\frac12}=1$). Inspection of
the following estimates shows that we do indeed lose a combined
$\frac34$ weight in $\xi$ on the right-hand side.

\begin{lemma}\label{lem:8.5}
Let $\alpha > \frac14$. Then
\[
\| R^{-\frac32} f g\|_{H^{\alpha+\frac14}_\rho} \lesssim
\|f\|_{H^{\alpha+\frac12}_\rho} \|g\|_{H^{\alpha+\frac12}_\rho}
\]
respectively
\[
\| R^{-\frac32} f g\|_{H^{\alpha}_\rho} \lesssim
\|f\|_{H^{\alpha+\frac14}_\rho} \|g\|_{H^{\alpha+\frac12}_\rho}
\]
for all $f,g$ so that the right-hand sides are finite.
\end{lemma}

\begin{proof}
We first use the above pointwise bound to obtain an $L^2$ estimate,
\[
\|R^{-\frac32} P_{\lambda_1} f P_{\lambda_2} g\|_{L^2} \lesssim \min\{ \lambda_1,\lambda_2\}
\|  P_{\lambda_1} f\|_{L^2} \|P_{\lambda_2} g\|_{L^2}
\]
This suffices for both of the above estimates provided that the
output is measured at frequency $\sigma \lesssim
\max\{\lambda_1,\lambda_2\}$. Indeed, in that case
\begin{align*}
  &\sum_{\lambda_1,\lambda_2}\;\sum_{\sigma<\max(\lambda_1,\lambda_2)}
  \sigma^{\alpha+\frac14} \big\|P_\sigma [R^{-\frac32} P_{\lambda_1} f
  P_{\lambda_2} g]\big\|_2 \\
  &\lesssim \sum_{\lambda_1>\lambda_2}
  \lambda_1^{\alpha+\frac14}\lambda_2 \|P_{\lambda_1} f\|_2
  \|P_{\lambda_2}g\|_2 + \sum_{\lambda_1\le \lambda_2}\lambda_2^{\alpha+\frac14}\lambda_1 \|P_{\lambda_1} f\|_2
  \|P_{\lambda_2}g\|_2 \\
  &\lesssim \sum_{\lambda_1>\lambda_2} \lambda_1^{-\frac14}
  \lambda_2^{\frac12-\alpha}
  \|f\|_{H_\rho^{\alpha+\frac12}}\|g\|_{H_\rho^{\alpha+\frac12}}
\end{align*}
which gives the desired bound since $\alpha>\frac14$.

\noindent For larger $\sigma$, however, we need some additional
decay. For this we compute using integration by parts
\[
\begin{split}
\Bla  R^{-\frac32} P_{\lambda_1} f P_{\lambda_2} g, P_\sigma h \Bra
&= \Bla  R^{-\frac32} P_{\lambda_1} f P_{\lambda_2} g, \cL^k
\cL^{-k} P_\sigma h \Bra
\\ &= \Bla  \cL^{k} (R^{-\frac32} P_{\lambda_1} f P_{\lambda_2} g),  \cL^{-k} P_\sigma h \Bra
\end{split}
\]
To justify the integration by parts we observe that near $R=0$ we have
\[
 P_{\lambda_1} f (R) = R^\frac32 q(R^2), \qquad q \text{ analytic }
\]
Then the bilinear form is given by
\[
R^{-\frac32} P_{\lambda_1} f P_{\lambda_2} g = R^\frac32 q(R^2), \qquad q \text{ analytic }
\]
which successively implies that (recall $\cL_0 R^{\frac32}=0$)
\[
 \cL^{k} (R^{-\frac32} P_{\lambda_1} f P_{\lambda_2} g) = R^\frac32
 q(R^2), \qquad q \text{ analytic }
\]
We claim that  we can estimate the left-hand side here in $L^2$ by
\begin{equation}
\|   \cL^{k} (R^{-\frac32} P_{\lambda_1} f P_{\lambda_2} g)\|_{L^2}
\lesssim \min\{ \lambda_1,\lambda_2\} \max\{ \lambda_1,\lambda_2\}^k
\label{biest}\|  P_{\lambda_1} f\|_{L^2} \|P_{\lambda_2} g\|_{L^2}
\end{equation}
Given the above integration by parts, this implies that
\[
\left| \Bla  R^{-\frac32} P_{\lambda_1} f P_{\lambda_2} g, P_\sigma h
  \Bra \right| \lesssim \min\{ \lambda_1,\lambda_2\} \max\{ \lambda_1,\lambda_2\}^k \sigma^{-k}
\|  P_{\lambda_1} f\|_{L^2} \|P_{\lambda_2} g\|_{L^2} \| P_\sigma h\|_{L^2}
\]
and further
\[
\|  P_\sigma (R^{-\frac32} P_{\lambda_1} f P_{\lambda_2} g)\|_{L^2}
\lesssim \min\{ \lambda_1,\lambda_2\} \max\{ \lambda_1,\lambda_2\}^k
\sigma^{-k} \|  P_{\lambda_1} f\|_{L^2} \|P_{\lambda_2} g\|_{L^2}
\]
thus providing the additional  decay for large $\sigma$.

It remains to prove \eqref{biest}.  We assume that $\lambda_1 <
\lambda_2$ and use different bounds depending on whether $R$ is small
or large.  Assume first that $R < \lambda_2^{-\frac12}$. Then we start
from
\begin{equation}\label{eq:cLkiden}
\cL^k(R^{-\frac32} P_{\lambda_1}f(R) \, P_{\lambda_2} g(R)) =
\il_0^\infty\il_0^\infty p_{\lambda_1}(\xi)p_{\lambda_2}(\eta) \cL^k
[R^{-\frac32} \phi(R,\xi)\phi(R,\eta)]
\,\hat{f}(\xi)\hat{g}(\eta)\,\rho(\xi)\rho(\eta)\,d\xi d\eta
\end{equation}
Next, we claim that
\begin{equation}\label{eq:cLkdoppel}
 \|\cL^k [R^{-\frac32} \phi(R,\xi)\phi(R,\eta)]
\|_{L^2(0,\lambda_2^{-\frac12})} \less \lambda_2^{k-1}
\end{equation}
If true, then combining \eqref{eq:cLkiden} and~\eqref{eq:cLkdoppel}
via Minkowski and Cauchy-Schwarz yields
\[
\| \cL^k(R^{-\frac32} P_{\lambda_1}f(R) \, P_{\lambda_2}
g(R))\|_{L^2(0,\lambda_2^{-\frac12})} \lesssim \lambda_2^k \lambda_1
\|P_{\lambda_1} f\|_2 \|P_{\lambda_2} g\|_2
\]
as desired. To prove~\eqref{eq:cLkdoppel}, consider first $k=0$.
Then
\[
\|R^{-\frac32} \phi(R,\xi)\phi(R,\eta)
\|_{L^2(0,\lambda_2^{-\frac12})} \less
\left(\int_0^{\lambda_2^{-\frac12}} R^3\,dR\right)^{\frac12}\less
\lambda_2^{-1}
\]
The higher $k$ cases now follow from Proposition~\ref{pphitheta},
which allows us to write
\[
R^{-\frac32} \phi(R,\xi)\phi(R,\eta) =  R^\frac32 q(R^2,\xi R^2,\eta
R^2),
\qquad \text{ $q$ analytic}
\]
Then, following our previous discussion concerning applications of
$\cL^k$, we obtain
\[
\cL^k(R^{-\frac32} \phi(R,\xi)\phi(R,\eta)) = \sum_{\ell+m \leq k}
R^\frac32  \xi^\ell \eta^m q_{\ell m}(R^2,\xi R^2,\eta R^2), \qquad
\text{ $q_{\ell m}$ analytic}
\]
which implies ~\eqref{eq:cLkdoppel}.

\noindent For large $R$ we use the product rule to write
\[\begin{split}
 \cL^{k} (R^{-\frac32} P_{\lambda_1} f P_{\lambda_2} g) &=
\sum_{\substack{2i+2j \leq 2k-\ell-m\\\ell,m=0,1}}   W_{ij}^{\ell m}
(R)\,
\partial_R^\ell \cL^i P_{\lambda_1} f \cdot \partial_R^m \cL^j
P_{\lambda_2} g\\
|W_{ij}^{\ell m} (R)| &\lesssim R^{-2(k-i-j) +\ell+m-\frac32}
\end{split}
\]
Then we have
\[
\|  \cL^{k} (R^{-\frac32} P_{\lambda_1} f P_{\lambda_2}
g)\|_{{L^2(\lambda_2^{-\frac12},\infty)}} \lesssim
\sum_{\substack{2i+2j \leq 2k-\ell-m\\\ell,m=0,1}}
\lambda_2^{(k-i-j)
  -\frac{\ell+m}2}   \|R^{-\frac32} \partial_R^\ell \cL^i
P_{\lambda_1} f \cdot \partial_R^m \cL^j P_{\lambda_2} g\|_{L^2}
\]
We use Lemmas~\ref{dpl},~\ref{linf} to bound the first factor in
$L^\infty$ and the second in $L^2$. This gives
\[
\|  \cL^{k} (R^{-\frac32} P_{\lambda_1} f P_{\lambda_2}
g)\|_{{L^2(\lambda_2^{-\frac12},\infty)}} \lesssim \lambda_2^k
\lambda_1 \| P_{\lambda_1} f\|_{L^2} \|P_{\lambda_2} g\|_{L^2}
\]
as desired.
\end{proof}

Finally, in order to estimate the higher order terms in the Taylor
expansion of the $\sin$ and $\cos$ functions in the proposition we
also prove a trilinear estimate:

\begin{lemma}
Let $\alpha > 0$. Then
\[
\| R^{-1} f g h\|_{H_\rho^{\alpha}} \lesssim
\|f\|_{H^{\alpha+\frac12}_\rho} \|g\|_{H_\rho^{\alpha+\frac12}}
\|h\|_{H_\rho^{\alpha}}
\]
for all $f,g,h$ so that the right-hand side is finite.
\label{tril}\end{lemma}

\begin{proof}
The pointwise bounds above imply the following $L^2$ estimate,
\begin{equation}\label{eq:3easy}
\|R^{-1} P_{\lambda_1} f P_{\lambda_2} g P_{\lambda_3} h
\|_{L^2} \lesssim \min_{i \neq j} \{ \lambda_i^\frac14 \lambda_j^\frac34\}
\| P_{\lambda_1} f \|_{L^2} \| P_{\lambda_2} g \|_{L^2}
\|P_{\lambda_3} h\|_{L^2}
\end{equation}
which again suffices to estimate the output at frequency $ \sigma
\le \lambda :=\max\{\lambda_1,\lambda_2,\lambda_3\}$. To see this,
we write, with the summation variables
$\sigma,M,\lambda_1,\lambda_2,\lambda_3\in\{2^j\}_{j=0}^\infty$,
\[
\sum_{\lambda_1,\lambda_2,\lambda_3} \sum_{\sigma \le \lambda}
P_\sigma (R^{-1} P_{\lambda_1} f P_{\lambda_2} g P_{\lambda_3} h) =
\sum_{M} \sum_{\substack{\lambda_1,\lambda_2,\lambda_3\\\lambda
  \geq M}}
P_{\lambda/M} \big(R^{-1} P_{\lambda_1} f P_{\lambda_2} g
P_{\lambda_3} h\big)
\]
On the right-hand side we distinguish the cases
$\lambda_1\le\lambda_2<\lambda_3$, $\lambda_1\le \lambda_3\le
\lambda_2$, $\lambda_3<\lambda_1\le\lambda_2$. We only treat the
first case, the other two being similar and easier. Thus, we
estimate the right-hand side for fixed $M$ as follows:
\[
\begin{split}
\Big\|  \sum_{\substack{\lambda_1\le\lambda_2<\lambda_3\\\lambda_3
  \geq M}}
P_{\lambda_3/M} \big(R^{-1} P_{\lambda_1} f &P_{\lambda_2} g
P_{\lambda_3} h\big)\Big\|_{H^\alpha_\rho}^2  \lesssim
\sum_{\lambda_3
> M} \Big(\frac{\lambda_3}{M}\Big)^{2\alpha}\|P_{\lambda_3} h\|^2_{L^2}\Big (
\sum_{\lambda_1\le\lambda_2} \lambda_1^\frac34 \lambda_2^\frac14 \|
P_{\lambda_1} f \|_{L^2} \| P_{\lambda_2} g \|_{L^2} \Big)^2
\\& \lesssim \sum_{\lambda_3 > M}
\Big(\frac{\lambda_3}{M}\Big)^{2\alpha} \|P_{\lambda_3} h\|^2_{L^2}
\Big( \sum_{\lambda_1<\lambda_2} \lambda_1^{\frac14-\alpha}
\lambda_2^{-\frac14-\alpha}  \| f \|_{H^{\alpha+\frac12}_\rho} \| g
\|_{H^{\alpha+\frac12}_\rho} \Big)^2 \\
&\lesssim M^{-2\alpha}
 \| f \|_{H^{\alpha+\frac12}_\rho}^2 \| g \|_{H^{\alpha+\frac12}_\rho}^2
\| h\|_{H^{\alpha}_\rho}^2
\end{split}
\]
The summation with respect to $M$ is trivial.

\noindent  For higher frequency outputs we need some additional
decay,
\[
\|P_\sigma (R^{-1} P_{\lambda_1} f P_{\lambda_2} g P_{\lambda_3} h)
\|_{L^2}  \lesssim  \min_{i \neq j} \{ \lambda_i^\frac14
\lambda_j^\frac34\} \max\{ \lambda_1,\lambda_2,\lambda_3\}^k
\sigma^{-k} \| P_{\lambda_1} f \|_{L^2} \| P_{\lambda_2} g \|_{L^2}
\|P_{\lambda_3} h\|_{L^2}
\]
This in turn is a consequence of  the estimate
\[
\| \cL^k (R^{-1} P_{\lambda_1} f P_{\lambda_2} g P_{\lambda_3} h)
\|_{L^2}  \lesssim  \min_{i \neq j} \{ \lambda_i^\frac14
\lambda_j^\frac34\} \max\{ \lambda_1,\lambda_2,\lambda_3\}^k \|
P_{\lambda_1} f \|_{L^2} \| P_{\lambda_2} g \|_{L^2} \|P_{\lambda_3}
h\|_{L^2}
\]
which is proved in the same manner as \eqref{biest}.
\end{proof}

These lemmas now imply Proposition~\ref{prop:uff}. Indeed,
we express the cosine-map in~\eqref{eq:epsmap1} in the form
\[
R^{-\frac12} (\cos (2 \tileps{R^{-\frac{1}{2}}})-1) =  R^{-\frac32} \tileps^2
q( R^{-1} \tileps^2) \qquad  \text{ $q$ entire}
\]
The first factor is bounded by
\[
\|R^{-\frac32}\tileps^2 \|_{H^\alpha_\rho} \less
\|\tileps\|_{H^{\alpha+\frac12}_\rho}^2
\]
while for $q$ we use its Taylor series together with Lemma~\ref{tril},
which shows that as a multiplication operator the factor $ R^{-1}
\tileps^2$ can be bounded by
\begin{equation}
\| R^{-1} \tileps^2 \|_{H^\alpha_\rho \to H^\alpha_\rho} \lesssim
\|\tileps\|_{H^{\alpha+\frac12}_\rho}^2
\label{linop}\end{equation}

Similarly, we write the sine-map from~\eqref{eq:epsmap2} in the form
\[
 R^{-\frac32}(\sin
(2{\tileps}{R^{-\frac{1}{2}}})-2{\tileps}{R^{-\frac{1}{2}}})
= R^{-3} \tileps^3 q(R^{-1} \tileps^2)
\]
For the first factor we apply Lemma~\ref{lem:8.5} twice to estimate
 \[
\| R^{-3} \tileps^3\|_{H^\alpha_\rho} \less  \|R^{-\frac32}
\tileps^2\|_{H^{\alpha+\frac14}_\rho}
\|\tileps\|_{H^{\alpha+\frac12}_\rho} \less
\|\tileps\|_{H^{\alpha+\frac12}_\rho}^3
 \]
while for the $q$ factor we use again \eqref{linop}.

\section{Proof of the main theorem.}
\label{sec:iter}

\noindent Here we summarize how to assemble together the elements of
the proof. Fixing $\nu > \frac12$ we begin with the approximate
solution $u_{2k-1}$ given by Theorem~\ref{thm:sec2} and with the
corresponding error $e_{2k-1}$. The index $k$ is chosen sufficiently
large, depending on $\nu$. Apriori both $u_{2k-1}$ and $e_{2k-1}$
are defined only inside the cone $\{ r \leq t\}$. We can extend them
to functions with similar regularity supported in a double cone $\{
r \leq 2t\}$. This extension is done crudely, without any reference
to the equation but insuring the matching on the cone for all
derivatives which are meaningful.

\noindent With these choices for $u_{2k-1}$ and $e_{2k-1}$ we seek
to solve \eqref{equation2} backward in $\tau$ and find a solution
$\tileps$ so that
\begin{equation}
\| \tileps(\tau) \|_{H^{\alpha+\frac12}_\rho} \lesssim \tau^{2-N},
\qquad \big\|\big(
\partial_\tau +\frac{\lambda_\tau}{\lambda} \big) \tileps(\tau)
\big\|_{H_\rho^{\alpha}} \lesssim \tau^{1-N}, \qquad N \leq 2k
\label{fds}\end{equation} Here the exponent $\alpha$ is chosen so
that
\[
\frac14 < \alpha < \frac{\nu}{2}
\]
The first bound is solely dictated by estimates for the cubic term
in the nonlinearity. The second one is a consequence of the
regularity of $e_{2k-1}$; namely, $e_{2k-1}$ has a singularity of
type $(1-a)^{\nu -\frac12} \log^m (1-a)$ on the cone. Therefore, if
$\alpha < \nu/2$, then for $e_{2k-1}$ we have the bound
\[
\| \lambda^{-2}R^{\frac12} e_{2k-1}(t(\tau),\lambda^{-1}R) \|_{H^\alpha_\rho} \lesssim \tau^{-2k+2}
\]
Using the transference identity we recast \eqref{equation2} for
$\tileps$ in the form \eqref{final} with $x = \calF \tileps$. By
virtue of Propositions~\ref{proph}, \ref{propn}, \ref{l2k} we can
solve \eqref{final}  using the contraction principle with respect to
the norm
\[
\| x\|_{L^{\infty,N-2}L^{2,\alpha+\frac12}_\rho} + \|(\partial_\tau -
\frac{\lambda_\tau}{\lambda})x\|_{L^{\infty,N-1}L^{2,\alpha}_\rho}
\]
Using again  the transference identity and Proposition~\ref{l2k} we
return back to $\tileps$, which has the regularity \eqref{fds}. Now
eventually we have to return to the original coordinates $(t,r)$ as
well as the function $\epsilon(t,r)$.  For this we define the map
\[
 u(R) \to T u(R,\theta) = e^{i\theta} R^{-\frac12} u(R)
\]
where the right hand side is interpreted as a function in $\R^2$
expressed in polar coordinates $(R,\theta)$. It is easy to see
that this is an isometry
\[
T: L^2(\R^+) \to L^2(\R^2)
\]
Then for the corresponding Sobolev spaces we have

\begin{lemma} \label{lem:sobolev}
For any $\alpha \geq 0$ we have
\[
\|u\|_{H_{\rho}^{\alpha/2}(\R^+)} \asymp \| T u\|_{H^{\alpha}(\R^2)}
\]
in the sense that if one side is finite then the other is finite
and they have comparable sizes.
\end{lemma}

\begin{proof}
  The spaces $H_{\rho}^{\beta}(\R^+)$ are defined using fractional
  powers of the operator $\cL$. However, we can also define them using
  fractional powers of the operator $\cL_0$ since the difference $\cL
  -\cL_0$ is bounded in $L^2$ and also in any $H^\beta_\rho$. This is
  easily seen if $\beta$ is an integer, and for noninteger values
  it follows by interpolation.

Then the conclusion of the lemma follows from the identity
\[
\Delta T u = T \cL_0 u
\]
which is valid whenever $u \in L^2$ and $\cL_0 u \in L^2$.
\end{proof}

To pass from $u(\tau, R)$, or alternatively $u(t,r)$,  to the
co-rotational wave map in terms of the ambient coordinates of
$\R^{3}\supset S^{2}$,  observe that these coordinates are given by
$\phi\circ T(u)$, where $\phi:\R^{2}\longrightarrow S^{2}\subset
\R^{3}$ is given by
\begin{equation}\nonumber
\phi(\rho e^{i\theta})=(\cos \rho, \sin\rho\cos\theta, \sin\rho \sin\theta)
\end{equation}
It is then easily seen that $\phi\circ T(u)\in
H^{2\alpha+1}(\R^{2})$, interpreted component-wise. We have now
constructed a wave map on the cone $r\leq t$, $0<t<t_{0}$, which is
of class $H^{1+\nu-}$ on the closure of the cone. To get a solution
on all of $\R^{2+1}$, extend the solution $\partial_{t}u(t_{0},.),
u(t_{0},.)$ at time $t=t_{0}$ to all of $\R^{2}$ within the same
smoothness and equivariance class. Call the corresponding wave map
$\tilde{u}(t,r)$. We claim that this wave map extends to
$(0,t_{0}]\times \R^{2}$ and is of class $H^{1+\nu-}$ until
breakdown at time $t=0$. Indeed, by finite propagation speed
$\tilde{u}(t,r)$ is given by $u(t,r)$ on the light cone $r\leq t,
0<t\leq t_{0}$. Furthermore, the $\tilde{u}$ does not develop
singularities on the interval $0<t<t_{0}$, as this could only happen
outside the light cone, where energy concentration is precluded by
the equivariance condition. The fact that singularity formation is
tantamount to an energy concentration scenario is a consequence of
\cite{Sh-Tah}, \cite{Tao}, \cite{Tat}. This concludes the proof of
Theorem~\ref{Main}.


\begin{thebibliography}{99}

\bibitem[AbrSte]{AS} Abramowitz, M., Stegun, I. {\em Handbook of
    mathematical functions with formulas, graphs, and mathematical
    tables.} National Bureau of Standards Applied Mathematics Series,
  55, U.S.~Government Printing Office, Washington, D.C.~1964

\bibitem[Bi] {Bi} Bizon, P., Tabor, Z., {\em Formation of
    singularities for equivariant $2+1$-dimensional wave maps into the
    $2$-sphere}, Nonlinearity 14(2001), no. 5, 1041-1053
\bibitem[Chr-Tah]{Chr-Tah} Christoulou, D., Tahvildar-Zadeh, Sh., {\em
    On the regularity of spherically symetric wave maps}, Comm. Pure
  Appl. Math. 46(1993), 1041-1091

\bibitem[Co]{Co} Cote, R., {\em Instability of non-constant harmonic
    maps for the $2+1$-dimensional equivariant wave map system},
  preprint

\bibitem[DS]{DS} Dunford, N., Schwartz, J. {\em Linear operators. Part II.}
 Wiley Classics Library. John Wiley \& Sons, Inc., New York, 1988.

\bibitem[GesZin]{GZ} Gesztesy, F., Zinchenko, M. {\em On spectral
    theory for Schr\"odinger operators with strongly singular
    potentials}, Math.\ Nachr.~279 (2006), 1041-1082.


\bibitem[KrSch1]{KS1} Krieger, J., Schlag, W., {\em Stable manifolds
    for all monic supercritical focusing nonlinear Schr\"odinger
    equations in one dimension.}  J.\ Amer.~Math.~Soc.\ 19 (2006),
  no.~4, 815--920.

\bibitem[KrSch2]{KS2} Krieger, J., Schlag, W., {\em On the focusing
    critical semi-linear wave equation}, to appear in Amer.\ Journal
  of Math.

\bibitem[KrSch3]{KrSch} Krieger, J., Schlag, W., {\em Non-generic
    blow-up solutions for the critical focusing NLS in 1-d}, preprint

\bibitem[Kri]{Kri} Krieger, J., {\em Global regularity of wave maps
    from $\R^{2+1}$ to ${\mathbf{H}}^{2}$}, CMP 250(2004), 507-580


\bibitem[Li]{Li} Isenberg, J., Liebling, S., {\em Singularity
    formation for $2+1$ wave maps}, J. Math. Phys. 43(2002), no. 1,
  678-683

\bibitem[Ro-St]{Ro-St} Rodnianski, I., Sterbenz, J., {\em On the
    formation of singularities in the critical $O(3)$ $\sigma$-model},
  preprint
\bibitem[Sch]{Sch} Schlag, W., {\em Stable manifolds for an orbitally
    unstable NLS}, to appear in Annals of Math.

\bibitem[Sha1]{Sha1} Shatah, J., {\em Weak solutions and development
    of singularities of the $SU(2)$ $\sigma$-model}, Comm. Pure Appl.
  Math. 41(1988), no. 4, 459-469.

\bibitem[Sha2]{Sha2} Cazenave, Th., Shatah, J., Tahvildar-Zadeh, A.
  Shadi, {\em Harmonic maps of the hyperbolic space and development of
    singularities in wave maps and Yang-Mills fields}

\bibitem[Sh-Tah]{Sh-Tah} Shatah, J., Tahvildar-Zadeh, Sh., {\em On the
    Cauchy problem for equivariant wave maps}, Comm. Pure Appl. Math.
  47(1994), no. 5, 719-754

\bibitem[Ste]{Stein} Stein, E. {\em Harmonic Analysis}, Princeton
University Press, 1993.

\bibitem[Str1]{Str1} Struwe, M., {\em Radially symmetric wave maps
    from $(1+2)$-dimensional Minkowski space to general targets},
  Calc. Var. Partial Differential Equations 16(2003), no. 4, 431-437.

\bibitem[Str2]{Str2} Struwe, M., {\em Equivariant wave maps in two
    space dimensions}, Comm. Pure Appl. Math. 56(2003), no. 7,
  815-823.
\bibitem[Tao]{Tao} Tao, T., {\em Global regularity of wave maps II},
  CMP 224(2001), 443-544

\bibitem[Tat]{Tat} Tataru, D., {\em Rough solutions for the wave maps
    equation}, Amer. J. Math 127(2005), no. 2, 293-377.

\bibitem[Wat]{Wat} Watson, G. {\em A treatise on the theory of Bessel
    functions}, Cambridge, 1944.

\end{thebibliography}
\end{document}